%% file: main.tex
\documentclass{article}
\usepackage[a4paper, total={6.5in, 10in}]{geometry}
\input{preamble}

\usepackage{lscape}
\usepackage{authblk}
\usepackage{tabularray}
\usepackage{caption}
\usepackage{subcaption}
\usepackage{siunitx}

\title{What is the optimal way to lie? \\From microscopic to kinetic descriptions of consensus control}

\author[1,*]{Sasha Glendinning}
\author[1]{Susana N. Gomes}
\author[1]{Marie-Therese Wolfram}
\affil[1]{Mathematics Institute, University of Warwick}
\affil[*]{Corresponding author}

\begin{document}

\maketitle

\vspace{0.5cm}

\begin{abstract}
We establish an approach for consensus control of opinion dynamics by introducing a liar to the classical system. The liar's aim is to steer the population towards consensus at their goal opinion by showing `apparent opinions', or `lies', to members of the population. We analyse this as an optimal control problem for how best to lie to a population in order to guarantee the consensus that the liar desires. We consider a range of regularisations, each motivated by some social convention, such as the liar wanting to present an opinion close to their true opinion. For each regularisation, we demonstrate the effect of instantaneous controls. Furthermore, we introduce a Boltzmann-type description for the corresponding kinetic system and present analysis and numerical results for the resulting Boltzmann and Fokker-Planck equations.
\end{abstract}

\section{Introduction}\label{sec:intro}

The field of opinion dynamics is concerned with modelling the results of interactions between people, typically leading to either consensus, when the whole population has the same opinion, or polarisation, characterised by opinion clusters. The fundamental models of Hegselmann-Krause \cite{krause2000discrete} and DeGroot \cite{degroot1974reaching} have inspired a rich literature of analysis for both microscopic and mean-field descriptions of opinion games. In particular, the field of optimal consensus control has gained a lot of interest over recent years. Much literature \cite{albi2014boltzmann,albi2014kinetic,nugent2024steering,albi2017mean,albi2017recent} exists for considering the problem of how to drive a population toward a given target consensus. For many of these formulations \cite{albi2014boltzmann,albi2014kinetic}, an extension of the model to a large-population limit is carried out by constructing a kinetic description, expressed by a Boltzmann-type equation as first introduced by Toscani \cite{toscani2006kinetic,pareschi2013interacting}.

Kinetic descriptions of opinion dynamics steer away from considering the individual opinions of agents and instead describe the evolution of the density of opinions on a macroscopic level. For large population sizes, this formulation induces a lower computational cost than microscopic models while remaining highly informative.

We introduce here an alternative model for optimal consensus control. Where previous formulations have added controls to their model by affecting opinions directly \cite{albi2014boltzmann,albi2014kinetic,albi2017mean}, we introduce controls that can be exactly interpreted as `lies' told by a single specified agent, referred to as `the liar'. From this clear interpretation, we can introduce regularisations (terms in the cost functional that limit or penalise the liar in some way) that are intuitive from a social standpoint. We also derive the Boltzmann-type description for the large-population limit of our model, following closely from the methods of \cite{albi2014boltzmann,albi2014kinetic,toscani2006kinetic}. The Boltzmann-type description provides a macroscopic view of the problem where we describe the evolution of the opinion density of truth-telling individuals in a large population. This description typically is computationally and analytically easier to resolve than the very large number of coupled ordinary differential equations necessary for a microscopic model.

The structure of the paper is as follows. In Section \ref{sec:leaders-followers}, we present an overview of formulations for optimal consensus control from the existing literature \cite{albi2014kinetic, albi2014boltzmann, nugent2024steering}. In Section \ref{sec:Micro-model}, we introduce the microscopic model for optimal consensus control through lying, along with a simple receding horizon strategy to calculate the corresponding instantaneous optimal control. We then present a variety of regularisation strategies, each motivated by some social convention, and compare the resulting dynamics. In Section \ref{sec:boltzmann}, we consider the large-population limit of this model, resulting in a Boltzmann-type description of the system. In this Section, we also examine the quasi-invariant limit of the Boltzmann equation and derive the corresponding Fokker-Planck equation. In Section \ref{sec:bounded-conf}, we discuss the solutions to our control model in the microscopic and kinetic form when agents interact under a bounded confidence kernel. Finally, in Section \ref{sec:multiple-liars} we perform numerical simulations to explore the behaviour of a system under the influence of multiple liars and conclude with a discussion in Section \ref{sec:conclusions}.

\section{Optimal consensus control}\label{sec:leaders-followers}

The classical Hegselmann-Krause (HK) formulation \cite{krause2000discrete} for modelling opinion dynamics is based on the idea that agents should update their opinions in a manner that promotes consensus. We assume a population of $N$ agents, each with some initial opinion $x_i(0)$ in an interval $\mathcal{I}\subset\mathbb{R}$ for $i=1, ..., N$. The opinion variables $x_i$ then evolve according to the following set of $N$ ordinary differential equations (ODEs),
\begin{equation}\label{eqn:HK-standard}
    \dot{x}_i = \frac{1}{N}\sum_{j=1}^NP(x_i, x_j)(x_j - x_i), \quad i=1, ..., N,
\end{equation}
where $P(x_i, x_j)$ is an interaction function, encoding the extent to which agent $i$ compromises with agent $j$. There are two principal types of interaction function that we will consider in this work. The first is the case where $P(x_i, \cdot) =P(x_i)$, i.e., the extent to which an individual with opinion $x_i$ interacts with an individual with opinion $x_j$ depends only on their own opinion, $x_i$. The second is the so-called \textit{bounded confidence} kernel, 
\begin{equation}\label{eq:bounded-con-introduciton}
    P(x_i, x_j) = \begin{cases}
    1\quad &\text{if } |x_j - x_i|<R,\\
    0\quad &\text{otherwise,}
\end{cases}
\end{equation}
where $R\in(0,|\mathcal{I}|]$ is a constant. In this regime, agents only interact with individuals with a similar opinion to their current opinion. This model has been studied extensively in \cite{motsch2014heterophilious, nugent2024bridging} and contains strong similarities to the Cucker-Smale model for bird flocking, discussed in \cite{cucker2007mathematics, lear2020grassmannian, shvydkoy2021dynamics}. The field of optimal consensus control is concerned with answering how best to influence agents in order to achieve some desired consensus opinion, $x_d\in\mathcal{I}$. Here, we define $x_d$ to be a consensus point if for all $i\in\{1, ..., N\}$, 
$$\lim_{t\to\infty}x_i(t) = x_d.$$

One formulation for optimal consensus control is to act on the evolution of agents' opinions directly with a control \cite{albi2014kinetic,albi2017mean}. In this case, the HK model (\ref{eqn:HK-standard}) becomes
\begin{equation*}
    \dot{x}_i = \frac{1}{N}\sum_{j=1}^NP(x_i, x_j)(x_j - x_i) + u,
\end{equation*}
for $i=1, ..., N$ where $u$ is a control that is equal for all agents but can depend on time, determined by minimising some cost functional
\begin{equation}\label{eq:directly-affect-cost}
    C(u) = \int_0^T\left(\frac{1}{N}\sum_{i=1}^N\frac{1}{2}(x_i - x_d)^2 + \frac{\nu}{2} u^2\right)\, dt.
\end{equation}
Here, $T>0$ is some arbitrary time horizon, $\nu>0$ is a regularisation parameter and $x_d\in\mathcal{I}$ is a desired consensus opinion. In this setting, we can interpret $u$ as being some external forcing on the agents' opinion, for example from advertisement or media. The possibility of targeted media, where each agent has a corresponding personal control $u_i$ is also explored in \cite{albi2017mean}.

Another formulation for optimal consensus control is known as `leaders and followers' \cite{albi2014boltzmann, during2009boltzmann}. The model assumes that we have two distinct populations of agents. The first is \textit{followers}, with opinions $x_i$ for $i=1, ..., N$, who update their opinions by communicating with each other and also by communicating with leaders. In contrast, the \textit{leaders}, with opinions $\tilde{x}_i$ for $i=1, ..., M$, update their opinions by communicating with each other but are also influenced by some external control. Mathematically, this model is written as
\begin{align*}
    \dot{x}_i &= \frac{1}{N}\sum_{j=1}^NP(x_i, x_j)(x_j - x_i) + \frac{1}{M}\sum_{k=1}^MP(x_i, \tilde{x}_k)(\tilde{x}_k - x_d),\quad &i=1, ..., N,\\
    \dot{\tilde{x}}_k &= \frac{1}{M}\sum_{j=1}^MP(\tilde{x}_k, \tilde{x}_j)(\tilde{x}_j - \tilde{x}_k) + u,\quad &k=1, ..., M.
\end{align*}
The control $u$ is chosen to minimise a cost functional,
\begin{equation}\label{eq:leaders-and-followers-cost}
    C(u) = \int_0^T\left(\frac{\psi}{M}\sum_{k=1}^M\frac{1}{2}(\tilde{x}_k - x_d)^2 + \frac{\mu}{M}\sum_{k=1}^M\frac{1}{2}(\tilde{x}_k - m_F(t))^2 + \frac{\nu}{2}u^2\right)\, dt,
\end{equation}
where $m_F(t)$ is the mean opinion of the followers, $\psi, \mu, \nu$ are non-negative constant parameters with $\psi + \mu=1$, representing the relative importance of each term in the cost functional. In staying close to the followers, the leaders strategy is to move close to the goal opinion while maintaining the amount of influence they have over the followers. In this setting, the control $u$ can be interpreted as an external forcing affecting the leaders. Therefore, the leaders and followers formulation could be relevant for modelling social media influencers and the effects of advertisement, for example.

Thus far, we have implicitly assumed that the agents lie on a fully connected network and that interactions are dependent only on the interaction function. The final control formulation that we shall discuss here is the idea of controlling the network on which the agents interact. In \cite{nugent2024steering}, the authors suppose that the agents are placed on a dynamic network, with weights $w_{ij}\in[0,1]$ encoding the strength of connection between agents $i$ and $j$. Mathematically, the model is written with $N$ ODEs describing the evolution of opinions $x_i$, and $N^2$ ODEs describing the evolution of the network weights:
\begin{align*}
    \dot{x}_i &= \frac{1}{k_i}\sum_{j=1}^Nw_{ij}\phi(x_j - x_i)(x_j - x_i),  &i=1, ..., N,\\
    \dot{w}_{ij} &=s(u_{ij})(\ell(u_{ij}) - w_{ij}),  &i, j = 1, ..., N, i\neq j.
\end{align*}
Here, $u_{ij}$ are control variables, $k_i$ is the in-degree of agent $i$ and $s, \ell$ are functions describing the speed and direction of the control respectively. The controls are chosen to minimise a cost functional that encourages agents toward the goal opinion $x_d$,
\begin{equation}\label{eq:network-control-cost}
    C(u) = \int_0^T \left(\sum_{i=1}^N\frac{1}{2}(x_i - x_d)^2 + \frac{\nu}{2} \sum_{i=1}^N\sum_{j=1}^N u_{ij}^2\right)\, ds.
\end{equation}
In the real world, social networks are often manipulated to enhance or force consensus; for example through the manipulation of social media algorithms. 

Each model we have discussed thus far relies on assuming some external figure exists outside of the model who wishes to drive the dynamics in a certain manner. The manner by which this figure can affect dynamics is a modelling choice, resulting in different levels of controllability. In what follows, we will present a formulation where the dynamics are controlled by placing assumptions on the behaviour of a particular agent in the system, a liar, rather than any external figure.

\section{Microscopic model of liars}\label{sec:Micro-model}

We now introduce a formulation for optimal consensus control where the population is driven to a desired consensus by the influence of a particular agent who is able to lie.

Suppose we have $N$ agents with opinions, $x_i\in\mathcal{I}$ for $i=1, ..., N$. Let $N-1$ of the agents' opinions evolve in the standard way, according to the HK model (\ref{eqn:HK-standard}). One agent lies by presenting an \textit{apparent opinion}, $y_i$, when in conversation with agent $i$, similarly to \cite{barrio2015dynamics}. We impose that the liar cannot present an opinion that is outside of the interval $\mathcal{I}$. For convenience and without loss of generality, we let the liar be agent 1. The aim of the liar is to bring the the system to consensus at their opinion $x_d\in\mathcal{I}$, which we assume does not change in time. We use the terms `goal opinion' and `true opinion' to describe $x_d$ interchangeably. Given our assumptions above, the dynamics are given by 
\begin{subequations}\label{eq:xidyn}
    \begin{equation}
        \dot{x}_1 = 0, \quad x_1(0) = x_d,
    \end{equation}
    and
    \begin{equation}
        \dot{x}_i = \frac{1}{N}\left(P(x_i, y_i)(y_i-x_i) + \sum_{j=2}^N P(x_i,x_j)(x_j-x_i) \right), \quad x_i(0) = x_i^0,
    \end{equation}
\end{subequations}
for $i=2, ..., N$ where $P(\cdot, \cdot)$ is an interaction function introduced in Section \ref{sec:leaders-followers}. The apparent opinion expressed by the liar, $y_i$, simply replaces the true opinion of agent 1 in the HK formulation. We write $y=(y_2, ..., y_N)^T$, defined as the solution of a minimisation problem with cost:
\begin{equation*}
    C(y) = \int_0^T\left( \frac{1}{N}\sum_{i=1}^N \frac{1}{2}(x_i - x_d)^2 + \nu\,\Psi(x,y)\, \right)\, ds, 
\end{equation*}
subject to the dynamics on $x_i$ given by (\ref{eq:xidyn}). Here, $T>0$ is some arbitrary time horizon, $\Psi$ is a convex function representing a cost or penalty associated with lying, known as a regularisation function, and $\nu>0$ is a constant parameter representing the relative weight of regularisation. We can think of $y_i$ for $i=2, ..., N$ as functions of $x(t)$ at time $t\geq 0$, where $x(t) = [x_2(t), ..., x_N(t)]^T\in\mathcal{I}^{N-1}$, mapping onto the interval $\mathcal{I}$, i.e., $y_i(x(t)): \mathcal{I}^{N-1}\to\mathcal{I}$ for all $t\geq0$ and $i=2, ..., N$.

We adopt a receding horizon strategy, often also referred to as an instantaneous control strategy, similar to that used by Albi, Herty and Pareschi \cite{albi2014boltzmann,albi2014kinetic}. The idea is to split the time interval $[0,T]$ into a finite number of $M+1$ time intervals with length $\Delta t$, defining $t^n = n\Delta t$ for $n=0, ..., M$. We assume that the control $y(t)$ is equal to a constant $y^n$ on each time interval $[t^n, t^{n+1}]$. Then each control value $y^n$ is determined as a solution to a reduced optimisation problem:
\begin{equation}\label{eq:general-reduced-problem}
    y^n = \text{arg}\min_{y\in\mathcal{I}^{N-1}}\left\{\int_{t^n}^{t^{n+1}} \left(\frac{1}{N}\sum_{i=1}^N \frac{1}{2}(x_i - x_d)^2 + \nu\,\Psi(x,y)\right)\, ds\right\}, 
\end{equation}
where $x_i$ evolves according to the dynamics (\ref{eq:xidyn}) with initial condition $\bar{x}_i = x_i(t^n)$.

The resulting control is typically suboptimal; however, the simplicity of the expressions that arise for the control will be useful when we consider a kinetic description of the model in Section \ref{sec:boltzmann}.

\subsection{No regularisation}\label{sec:Micro-model-no-reg}

We shall first consider the case where the liar acts under no regularisation, meaning that $\Psi(x,y) \equiv 0$ in (\ref{eq:general-reduced-problem}). Then, $y^n$ is defined according to the solution of the minimisation problem
\begin{equation*}
    y^n = \text{arg}\min_{y\in\mathcal{I}^{N-1}}\left\{\int_{t^n}^{t^{n+1}} \frac{1}{N}\sum_{i=1}^N \frac{1}{2}(x_i - x_d)^2 \, ds\right\},
\end{equation*}
subject to the dynamics of equation (\ref{eq:xidyn}). The Hamiltonian obtained via the maximum principle \cite{evans2005introduction} is
\begin{equation*}
    H(x, p, y) = \sum_{i=2}^Np_i\frac{1}{N}\left[ \sum_{j=2}^N P(x_i, x_j)(x_j-x_i) + P(x_i,y_i)(y_i-x_i)\right] -\frac{1}{N}\sum_{i=1}^N\frac{1}{2}(x_i-x_d)^2,
\end{equation*}
where $p=(p_2, .., p_N)$ is a vector of Lagrange multipliers. Taking partial derivatives with respect to $y_i$, $x_i$, and $p_i$ gives the following set of conditions for optimality:
$$ p_i\partial_{y_i}\{ P(x_i,y_i)(y_i-x_i)\}=0, \quad i=2, ..,, N,$$
$$ \dot{p}_i(t) = \frac{1}{N}(x_i-x_d) - \frac{1}{N}\left[ \sum_{j=2}^N R_{ij} + p_i\partial_{x_i}\{P(x_i, y_i)(y_i-x_i)\}\right], \quad p_i(t^{n+1}) = 0,\quad i=2, ..., N,$$
$$\dot{x}_i(t) = \frac{1}{N}\left(\sum_{j=2}^N P(x_i,x_j)(x_j-x_i) + P(x_i, y_i)(y_i-x_i)\right), \quad x_i(t^n) = \bar{x}_i,\quad i=2, ..., N$$
where
\begin{equation}\label{eq:R}
 R_{ij} = p_i\partial_{x_i}\{P(x_i, x_j)(x_j-x_i)\} + p_j\partial_{x_j}\{P(x_j, x_i)(x_i-x_j)\}.
 \end{equation}
We then discretise these equations. As there are no terminal costs, $p_i(t^{n+1}) = 0$ and we can discretise the equation for $\dot{p}_i(t)$ via a backward Euler scheme and obtain
\begin{equation}\label{eq:pin-basic}
    p_i^n = -\frac{\Delta t}{N}(x_i^{n+1} -x_d).
\end{equation}
The first equation in our list of optimality conditions then gives us the following
\begin{equation}\label{eq:noreg_optcontrol}
     (x_i^{n+1}-x_d)\partial_{y_i}\{P(x_i,y_i)(y_i-x_i)\}|_{(\bar{x}_i, y_i^n)} = 0,
\end{equation}
where $x_i^{n+1}$ is given by (\ref{eq:xidyn}) with initial conditions $\bar{x}=(\bar{x}_2, ..., \bar{x}_N)^T$. Therefore, if the position of $x_i$ allows it, the optimal choice of $y_i$ is to choose $y_i$ such that $x_i^{n+1} = x_d$, i.e., choose an admissible lie such that agent $i$ is in agreement with the liar in the next time-step. If there is no such admissible lie that can be told, then we choose $y_i$ such that
\begin{equation}\label{eq:no-regyi}
    \partial_{y_i}\{P(x_i, y_i)(y_i-x_i)\}|_{(\bar{x}_i, y_i^n)} = 0,
\end{equation}
i.e., choose $y_i$ to maximise the distance that $x_i$ can move in the given time-step toward $x_d$. If (\ref{eq:no-regyi}) also does not have an admissible solution, then the optimal control is given by a projection onto the interval $\mathcal{I}=[-1,1]$, say,
$$ y_i^n = \mathbb{P}_{\mathcal{I}}(\tilde{y}_i^n) = \begin{cases}
    1, \quad &\text{if } \tilde{y}^n_i>1,\\
    \tilde{y}_i^n, \quad &\text{if } \tilde{y}_i^n\in \mathcal{I}, \\
    -1, \quad &\text{if } \tilde{y}_i^n<-1,
\end{cases}$$
where $\tilde{y}_i^n$ is the solution to $x_i^{n+1} - x_d=0$ and $\mathbb{P}_{\mathcal{I}}$ is an operator projecting the solution $\tilde{y}_i^n$ onto the interval $\mathcal{I}$.
\begin{figure}[htb!]
    \centering
    \includegraphics[width=0.9\textwidth]{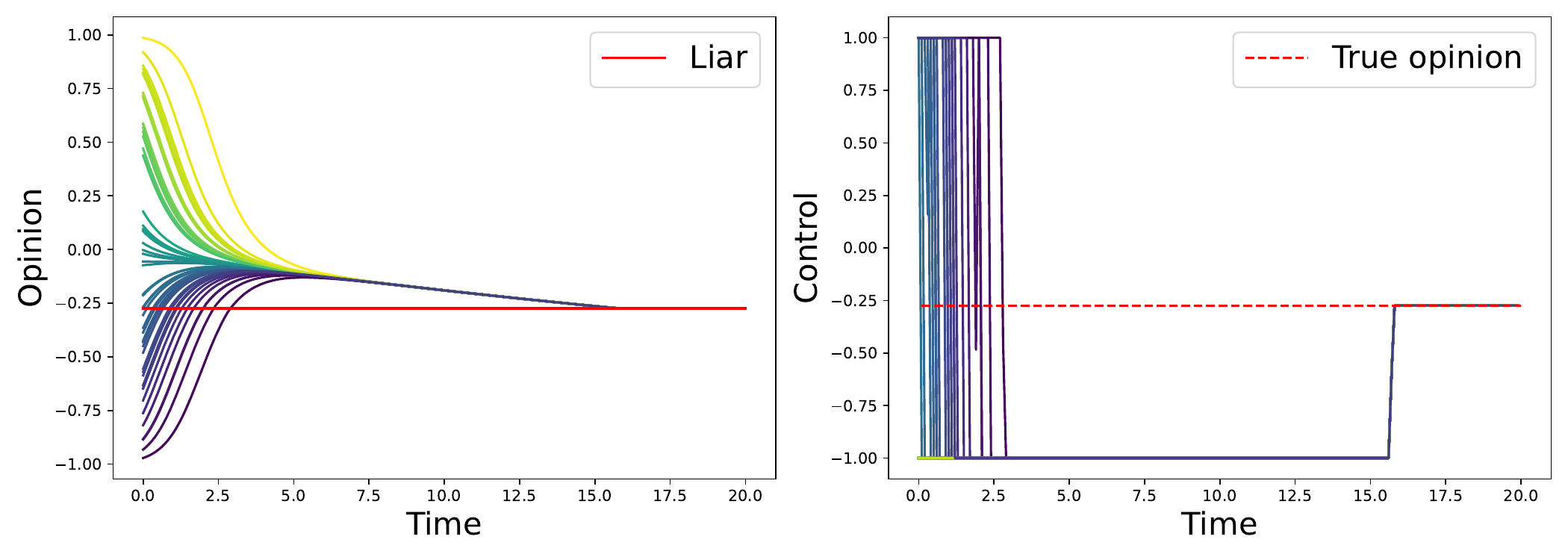}
    \caption{Trajectory of agents (left) and the corresponding controls, or lies, (right) in the case where the liar has no regularisation. Here, $N=50$ and $P(x,\cdot\,) = P(x) = 1-x^2$.}
    \label{fig:no-reg}
\end{figure}
\vspace{-0.5cm}

In Figure \ref{fig:no-reg}, we plot a simple example with an interaction function $P(x,\cdot\,) = P(x) = 1-x^2$ and $N=50$ agents. We observe that consensus at the liar's goal opinion is reached in finite time (which is not necessarily the case for all control formulations; see \cite{nugent2024steering} for examples) and quite rapidly. The corresponding control is generally bang-bang - the liar presents an opinion of $\pm 1$ initially, depending on the relative position of the agent they are interacting with. When communicating with an individual with opinion $x_i$, if $x_i>x_d$, then the liar will generally present an opinion $y_i = -1$ to cause the greatest rate of change in the opinion of agent $i$ toward $x_d$. Similarly, if $x_i<x_d$, then the liar will present an apparent opinion $y_i = 1$. The liar continues with this strategy until the opinions of the truth-telling agents are sufficiently close to $x_d$ and then suddenly changes to expressing their true/goal opinion so that consensus is achieved. It is interesting to note that qualitatively, this strategy is observed in the real world, when exaggeration leads to liars achieving their goals. An example of this is the formation of the Good Friday peace agreement in Northern Ireland, where both sides inflated early progress to encourage compromise, later pulling back to a milder narrative to consolidate consensus \cite{dixon2002political}.

However, we do note occasional moments where the liar presents an opinion $y_i^n$ that does not belong to the set of bang-bang controls discussed above, $y_i^n\notin\{\pm1, x_d\}$ (see the outlier around $t=2$ in the plot on the right in Figure \ref{fig:no-reg}). These moments correspond to when a truth-teller is very close to crossing the desired opinion $x_d$ at time-step $n$. The liar presents an opinion to ensure $x_i^{n+1}=x_d$ but the truth-telling dynamics are still active and the truth-teller moves away again from $x_d$ as they are influenced by their peers. This is a drawback of using a reduced horizon strategy where the control can only determine the optimal lie based on information obtained from looking one time-step into the future.

\remark{There are conditions under which the control corresponding to a system with interaction function $P(x_i,\cdot) = 1-x_i^2$ is not bang-bang. For example, suppose that agents are initialised such that half of the truth-telling agents have initial opinion $x_d+\varepsilon$ and the other half have opinion $x_d-\varepsilon$ with $x_d\in(-1,1)$ and $\varepsilon>0$. For sufficiently small $\varepsilon$, it can be shown that $\tilde{y}_i^n\in\mathcal{I}$ and $\tilde{y}_i^n\neq x_d$ for all $i=2, ..., N$.}\normalfont

The control profile can look very different for different interaction functions. Since the interaction function used to create Figure \ref{fig:no-reg} carries no dependence on how extreme the view of the liar is, the bang-bang control is very natural. In contrast, if we adopt a bounded confidence interaction kernel, then the liar is forced to hug the radius of interaction to be able to influence agents. This case will be discussed further in Section \ref{sec:bounded-conf}.

\subsection{Classic regularisation}\label{sec:Micro-model-std-reg}

In optimal control problems, a penalty is typically added to the cost functional for the control being `switched on'. We have seen examples of such cost functionals in Section \ref{sec:leaders-followers} in equations (\ref{eq:directly-affect-cost}), (\ref{eq:leaders-and-followers-cost}) and (\ref{eq:network-control-cost}). In our setting, the control being `switched off' corresponds to the liar telling the truth, i.e., $y_i = x_d$. We could interpret this penalty of deviating too far from the true opinion as the psychological cost of lying \cite{van2012learning}. 

In this case, we consider the evolution of truth-tellers' opinions to be given by (\ref{eq:xidyn}) but this time we write our reduced cost function as
\begin{equation}\label{eq:std-reg-cost}
    y^n = \arg\min_{y\in{\mathcal{I}^{N-1}}}\left\{\int_{t^n}^{t^{n+1}} \left(\frac{1}{N}\sum_{i=1}^N \frac{1}{2}(x_i - x_d)^2 + \frac{1}{N}\sum_{i=2}^N \frac{\nu}{2}(y_i-x_d)^2\right)\, ds\right\},
\end{equation}
where $\nu>0$ is again a regularisation parameter.

We proceed to derive the optimal instantaneous control with the same method as in Section \ref{sec:Micro-model-no-reg}, with details given in Appendix \ref{app:optimality-conditions}. From the optimality system, we conclude that the control will be given by a projection of the solution of the implicit equation
\begin{equation}\label{eq:SR-control}
    \tilde{y}_i^n = x_d - \frac{\Delta t}{\nu N}(x_i^{n+1}-x_d)\partial_{\tilde{y}_i}\{P(x_i,\tilde{y}_i)(\tilde{y}_i-x_i)\}|_{(\bar{x}_i, \tilde{y}_i^n)},
\end{equation}
onto the interval $\mathcal{I}$, written $y_i^n = \mathbb{P}_{\mathcal{I}}(\tilde{y}_i^n)$. As a sanity check, we can note that when $\nu$ is large (and therefore the regularisation term dominates in (\ref{eq:std-reg-cost})), the solution (\ref{eq:SR-control}) tends to $y_i^n = x_d$, i.e., the agent presents its true opinion. This makes intuitive sense as this regularisation penalises lying.

\begin{figure}[htb!]
    \centering
    \begin{subfigure}[b]{0.45\textwidth}
        \centering
        \includegraphics[width=\textwidth]{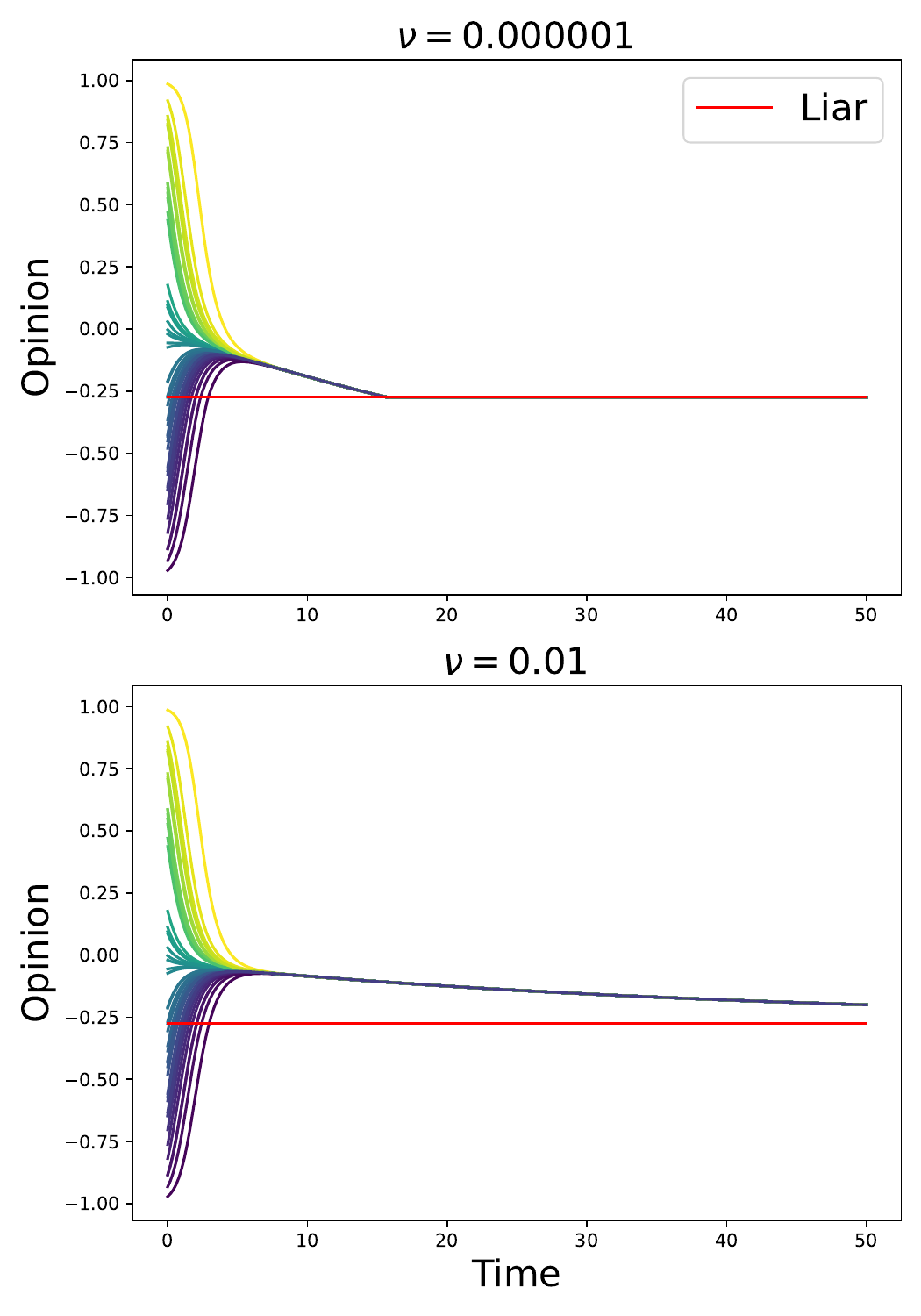}
        \caption{Trajectory}
    \end{subfigure}
    \begin{subfigure}[b]{0.45\textwidth}
        \centering
        \includegraphics[width=\textwidth]{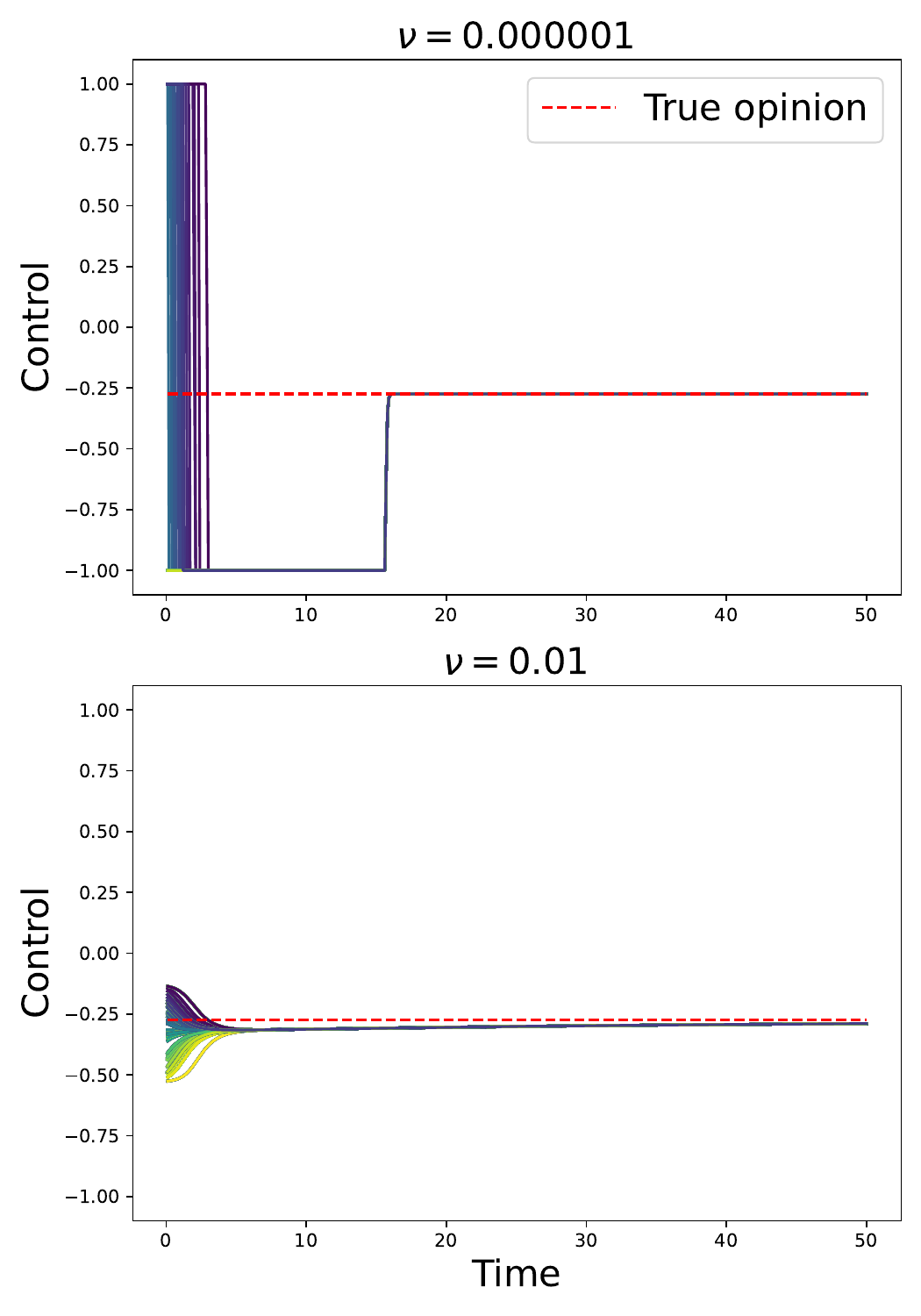}
        \caption{Control}
    \end{subfigure}
    \caption{Trajectory of agents and the corresponding controls in the case where the liar is penalised for deviating from their true opinion. Here, $N=50$, $P(x,\cdot\,) = P(x) = 1-x^2$ and we compare values of $\nu = 0.000001$ and $\nu = 0.01$. The figures on the left show the trajectory of the system under each of these regularisations and the figures on the right show the corresponding controls/lies.}
    \label{fig:std-reg}
\end{figure}

The plots in Figure \ref{fig:std-reg} display the outcome of this control strategy for increasing values of the regularisation parameter $\nu$. For very small $\nu$, the control is very similar to that which we have seen in Section \ref{sec:Micro-model-no-reg} with no regularisation. As $\nu$ increases, we see that the liar deviates less and less from their true opinion. Note that $P(x_i,\cdot\,) = P(x_i) = 1-x_i^2>0$ apart from when $x_i = \pm1$ and the probability of $x_i^0=\pm1$ for $x_i^0\sim U[-1,1]$ is zero. Therefore, the population must reach consensus \cite{motsch2014heterophilious} and thus, even when the liar can only tell the truth, the fact that their opinion remains constant means that consensus is reached at $x_d$.

A feature of the dynamics that are observed in Figure \ref{fig:std-reg} is that once the truth-tellers have reached an asymptotic relative consensus amongst themselves and are only influenced by the liar, the lie told to (or control on) every individual is the same. Hence, in this regime we can simplify the dynamics to two particles, one representing the truth-telling population and the other representing the lie told by the liar. With an Euler discretisation for the dynamics given by (\ref{eq:xidyn}) and control given by (\ref{eq:SR-control}), we have that the state of the system once the two-particle limit is reached is given by
\begin{align*}
    x^{n+1} &= x^n + \frac{\Delta t}{N}P(x^n, y^n)(y^n - x^n),\\
    y^n &= x_d -\frac{\Delta t}{\nu N}(x^{n+1}-x_d)\partial_y(P(x^n,y)(y-x^n))|_{y=y^n},
\end{align*}
where $x^n$ is the opinion of the truth-tellers in the population at time-step $n\geq0$ and $y^n$ is the lie told to the population at time-step $n\geq0$. Note that here we have removed the assumption that $y^n\in\mathcal{I}$ in order to simplify the dynamics. We consider the simple case where $P(x^n,y^n) = 1$ and recover first order difference equations where the solution $x^n$ is given by
\begin{equation}
    x^n = x_d\left(1 - \left(\frac{1-\frac{\Delta t}{N}}{1+ \frac{\Delta t^2}{\nu N^2}}\right)^n\right) + \left(\frac{1-\frac{\Delta t}{N}}{1+ \frac{\Delta t^2}{\nu N^2}}\right)^nx^0,
\end{equation}
for $n\geq 0$, where $x^0$ is the value of $x_i$ at the point where the truth-tellers meet relative consensus for $i=2, ..., N$. Since the liar is allowed to lie with any opinion $y\in\mathbb{R}$, we can take $x_d = 0$ without loss of generality in order to consider convergence rates toward $x_d$. The half life of the convergence, meaning the time taken for $x^0$ to decrease to $x^0/2$ is
$$ t_* = \frac{-\Delta t\log 2}{\log\left(1-\frac{\Delta t}{N}\right) - \log\left(1+ \frac{\Delta t^2}{\nu N^2}\right)},$$
where we define $t_* = n_*\Delta t$ for $n_*$ the smallest value of $n$ such that $x^n<x^0/2$. In the case where the liar may only tell the truth (taking $\nu\to\infty$), the half life is given by
$$ \tilde{t}_* = \frac{-\Delta t\log 2}{\log\left(1-\frac{\Delta t}{N}\right)}.$$
We see that for $\nu>0$, $t_*<\tilde{t}_*$ so lying is always a more efficient way to achieve consensus. Finally, we can consider how small we require our regularisation parameter to be in order to achieve a specified improvement on the half life. Suppose we wish to decrease the half life of the long-time dynamics by a factor of $k$, then we require
$$ \nu \leq \frac{\Delta t^2}{N^2}\left(\frac{1}{\left(1-\frac{\Delta t}{N}\right)^{1-k}-1}\right).$$
Considering the right hand side of this inequality as a function of $k$, we see that as $k$ increases our condition becomes more constrictive and we need to take $\nu$ to be very small. This demonstrates a trade-off between speed of convergence to the liar's goal opinion and the risk associated with taking $\nu$ to be small.

\subsection{Consistent control}\label{sec:Micro-model-consis-reg}

We will now consider a different regularisation on the liar, motivated by a risk or cost associated with telling different lies to different people. In psychology, this is referred to as the multiple-audience problem \cite{bond2004maintaining}. The dynamics of $x_i$ remain unchanged but we optimise over a reduced cost functional
\begin{equation}\label{eq:cost-fun-consis-control}
    y^n = \arg\min_{y\in{\mathcal{I}^{N-1}}}\left\{\int_{t^n}^{t^{n+1}} \left(\frac{1}{N}\sum_{i=1}^N \frac{1}{2}(x_i - x_d)^2 + \frac{\nu}{2N^2}\sum_{i=2}^N\sum_{j=2}^N\frac{1}{2}(y_i-y_j)^2\right)\, ds\right\},
\end{equation}
where again, $\nu>0$ is a regularisation parameter.

Following similar Hamiltonian methods as in Section \ref{sec:Micro-model-no-reg} (see Appendix \ref{app:optimality-conditions} for more details), the optimal control at time-step $n$, $y_i^n$, is given as the projection of the solution to
\begin{equation}\label{eq:CP-control}
    \tilde{y}_i^n = \frac{1}{N-1}\sum_{j=2}^N \tilde{y}_j^n -\frac{\Delta t}{\nu(N-1)}(x_i^{n+1}-x_d) \partial_{\tilde{y}_i}\{P(x_i, \tilde{y}_i)(\tilde{y}_i-x_i)\}|_{(\bar{x}_i, \tilde{y}_i^n)},
\end{equation}
onto the interval $\mathcal{I}$, i.e., $y_i^n = \mathbb{P}_{\mathcal{I}}(\tilde{y}_i^n)$. We can perform a sanity check by considering equation (\ref{eq:CP-control}) in the limit where $\nu\to\infty$ and the liar's strategy is dominated by the need to show a consistent projection. In this case, $y_i^n\to \frac{1}{N-1}\sum_{j=2}^N y_j^n = \bar{y}^n$ - the average of the projected opinions at time-step $n$. This demonstrates that there will be complete consistency in the lies that the liar tells the general populace. 
\begin{figure}[htb!]
    \centering
    \begin{subfigure}[b]{0.45\textwidth}
        \centering
        \includegraphics[width=\textwidth]{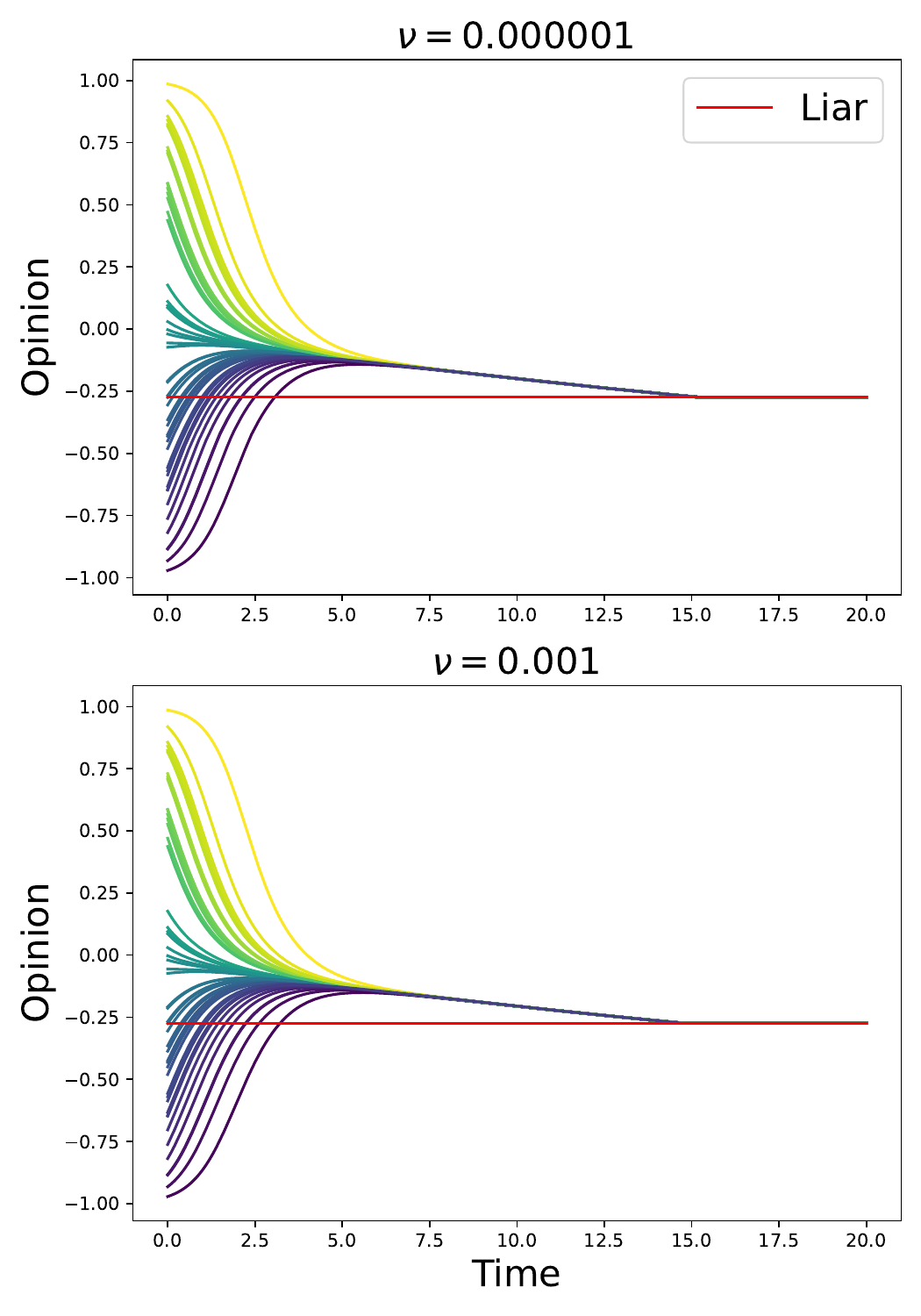}
        \caption{Trajectory}
    \end{subfigure}
    \begin{subfigure}[b]{0.45\textwidth}
        \centering
        \includegraphics[width=\textwidth]{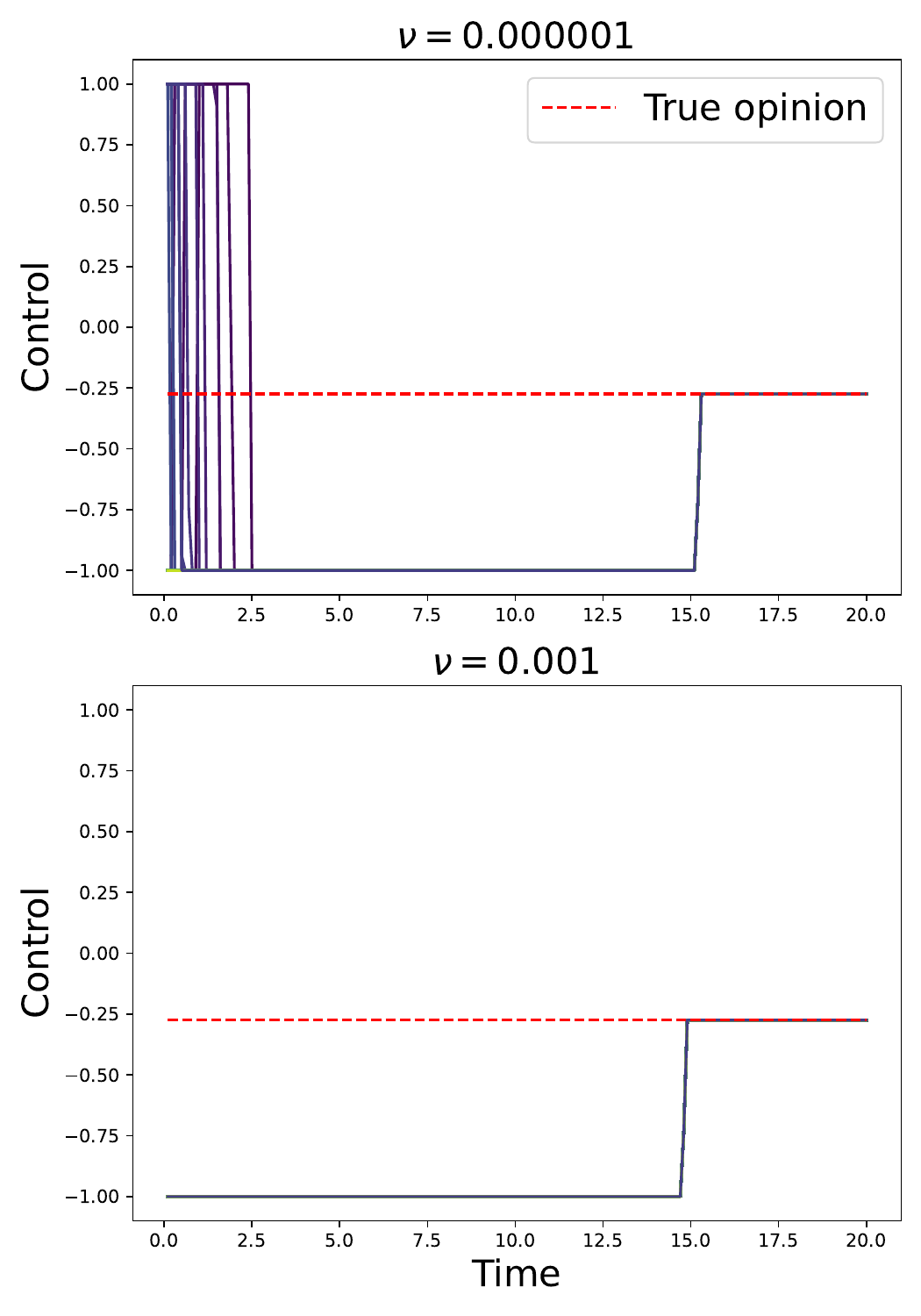}
        \caption{Control}
    \end{subfigure}
    \caption{Trajectory of agents and the corresponding controls in the case where the liar is penalised for expressing opinions inconsistently. Here, $N=50$, $P(x,\cdot\,) = P(x) = 1-x^2$ and we compare values of $\nu = 0.000001$ and $\nu = 0.001$.}
    \label{fig:consis-reg}
\end{figure}

An example of the results of this regularisation is shown in Figure \ref{fig:consis-reg}. For very small $\nu$, we see a bang-bang type control similar to that in Section \ref{sec:Micro-model-no-reg} but a small increase in $\nu$ leads to the liar being totally consistent with their lies. In the case $\nu = 0.001$, the fact that initially the optimal consistent lie is to take $y_i=-1$ for all $i=2, ..., N$ reflects the fact that the liar is placed asymmetrically in the opinion space, with more agents having $x_i>x_d$ than $x_i<x_d$. However, the sudden jump in projected opinion from $y_i = -1$ to $y_i = x_d$ is perhaps unrealistic, since it is generally frowned upon to drastically change opinion in such a short period of time. This `flip-flopping' can lower trust and leadership credibility \cite{allgeier1979waffle}, signal unreliability or be seen as opportunism \cite{nasr2024times}. This motivates our next regularisation.

\subsection{Time-consistent control}\label{sec:Micro-model-timecon}

In order to penalise the rate with which the liar changes their opinion in time, we need some approximation for the change of the lie in time, represented by $\dot{y}_i$. For consistency with our receding horizon strategy, we consider the Euler derivative approximation and at each time-step $n\geq 1$, given $y_i^{n-1}$, we penalise the quantity
$$ (\dot{y}_i^n)^2 \approx\left(\frac{y_i^n - y_{i}^{n-1}}{\Delta t}\right)^2.$$
The reduced problem at time-step $n$ is then
\begin{equation}\label{eq:TC-euler}
    y^n = \arg\min_{y\in{\mathbb{R}^{N-1}}}\left\{\int_{t^n}^{t^{n+1}} \left(\frac{1}{N}\sum_{i=1}^N \frac{1}{2}(x_i - x_d)^2 + \frac{1}{N}\sum_{i=2}^N \frac{\nu}{2}\left(\frac{y_i - y_i^{n-1}}{\Delta t}\right)^2\right)\, ds\right\}, \quad n\geq 1.
\end{equation}
Of course we also need to set initial conditions for $y_i^0$. A simple option is to take $y_i^0$ to be the solution to any of the previous optimisation problems, given by (\ref{eq:noreg_optcontrol}), (\ref{eq:SR-control}) or (\ref{eq:CP-control}) at time $t=0$. 

The optimal $y_i^n$ for the reduced problem, with optimality system detailed in Appendix \ref{app:optimality-conditions}, is given by a projection of $\tilde{y}_i^n$ onto $\mathcal{I}$ where $\tilde{y}_i^n$ satisfies
\begin{equation*}
    \tilde{y}_i^n = y_i^{n-1} - \frac{\Delta t^3}{\nu N}(x_i^{n+1}- x_d) \partial_{\tilde{y}_i}\{P(x_i, \tilde{y}_i)(\tilde{y}_i-x_i)\}, \quad n\geq 1.
\end{equation*}
Hence, $y_i^n = \mathbb{P}_{\mathcal{I}}(\tilde{y}_i^n).$ Our sanity check again shows that in a regularisation dominated setting where $\nu\to\infty$, we have $y_i^n\to y_i^{n-1}$ so the liar's projected opinion does not change in time. 

\begin{figure}[htb!]
    \centering
    \begin{subfigure}[b]{0.45\textwidth}
        \centering
        \includegraphics[width=\textwidth]{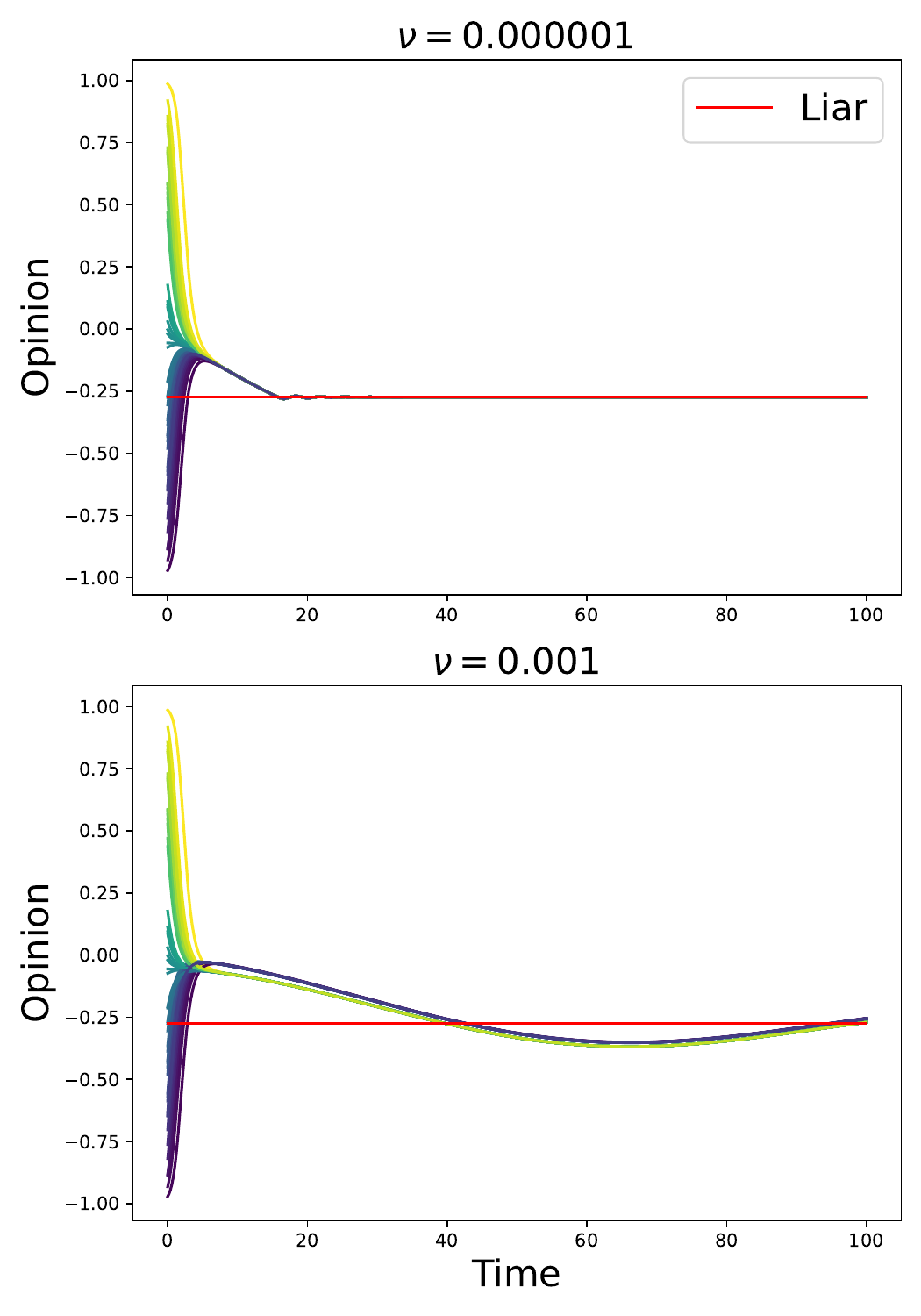}
        \caption{Trajectory}
    \end{subfigure}
    \begin{subfigure}[b]{0.45\textwidth}
        \centering
        \includegraphics[width=\textwidth]{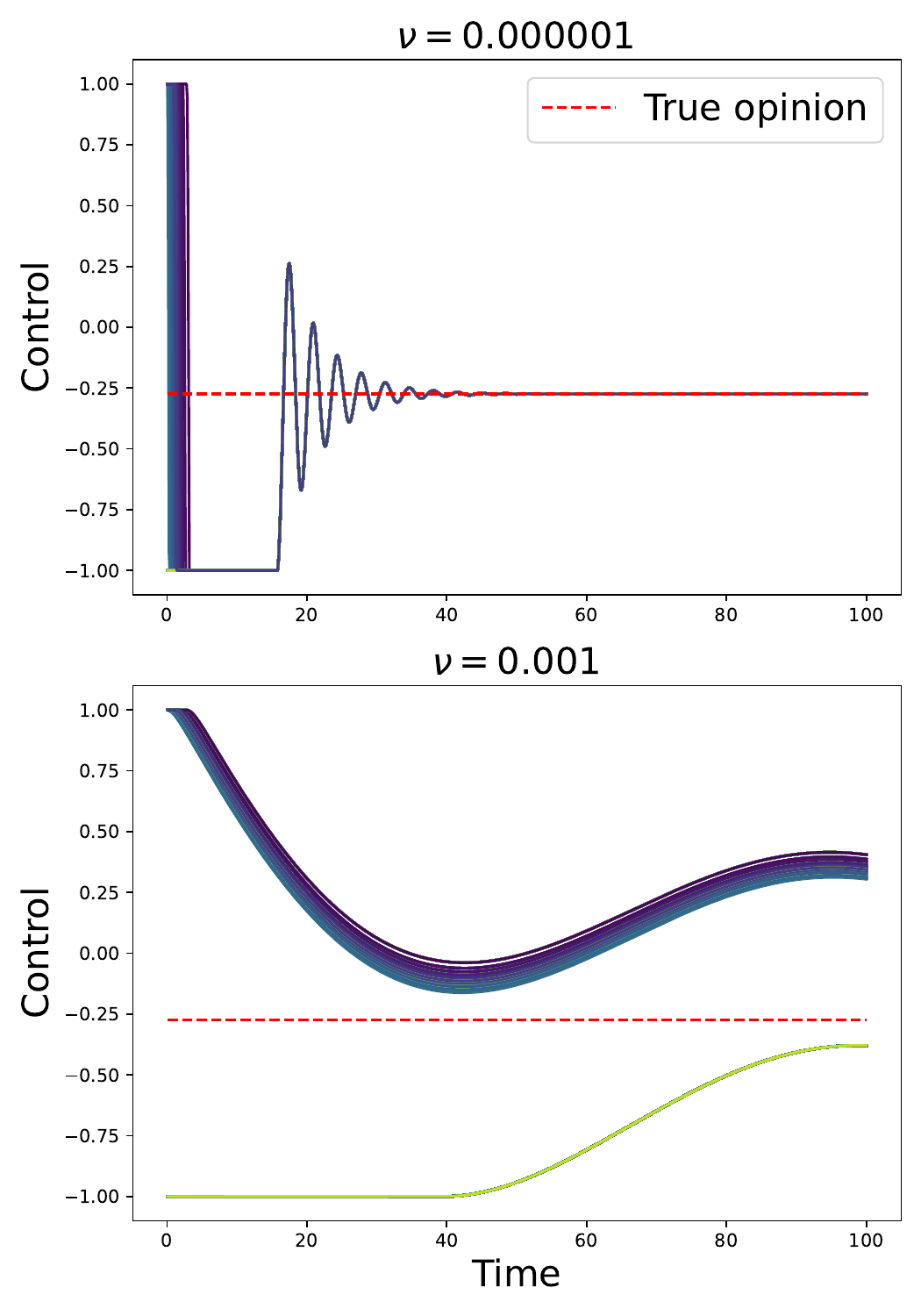}
        \caption{Control}
    \end{subfigure}
    \caption{Trajectory of agents and the corresponding controls in the case where the liar is penalised for expressing opinions inconsistently in time. The initial condition $y_i^0$ is given by the no-regularisation case, equation (\ref{eq:noreg_optcontrol}). Here, $N=50$, $P(x,\cdot\,) = P(x) = 1-x^2$ and we compare values of $\nu = 0.000001$ and $\nu = 0.001$.}
    \label{fig:timeconsis-noreg}
\end{figure}

Figure \ref{fig:timeconsis-noreg} shows some examples of trajectories and corresponding controls for this problem when we take the initial conditions $y_i^0$ to be given by the liar's optimal apparent opinions when there is no regularisation, i.e., $y_i^0$ satisfies equation (\ref{eq:noreg_optcontrol}). We observe that for small $\nu$, the control is almost the bang-bang control we have seen before but this time with some oscillation before settling at $y_i = x_d$. These oscillations persist as $\nu$ increases. This is due to the fact that it is too sharp of a change in opinion to suddenly jump to $y_i = x_d$. Instead, the control oscillates about $x_d$ to keep the rate of change in $y_i$ small before converging to $y_i = x_d$. We can analyse the appearance of oscillations in our simulations by again taking a two-particle simplification of the model in a similar way as in Section \ref{sec:Micro-model-std-reg}. We make the same assumptions: the truth tellers have reached an asymptotic consensus so are only influenced by the liar, $P(x_i,x_j) \equiv1$ and the liar is unconstrained so $y_i^n= \tilde{y}_i^n\in\mathbb{R}$. The problem is then reduced to the following discrete dynamical system,
\begin{align}\label{eqn:discrete-dyn-timecon}
\begin{split}
    x^{n+1} &= x^n + \frac{\Delta t}{N}(y^n - x^n), \quad n\geq0\\
    y^n &= y^{n-1} - \frac{\Delta t^3}{\nu N}(x^{n+1} - x_d) ,\quad n\geq 1,
\end{split}
\end{align}
with initial conditions $x^0 = X^0$ and $y^0 = Y^0$ given. If we again assume without loss of generality that $x_d=0$, we can solve the above system and obtain the following condition for the appearance of oscillatory solutions,
\begin{equation*}
    \nu<4\Delta t^2\left(1-\frac{\Delta t}{N}\right).
\end{equation*}
An example of this regime is shown in Figure \ref{fig:timeconsis-noreg} when $\nu = 0.000001$, $N=50$ and $\Delta t = 0.1$. This analysis demonstrates that oscillatory behaviour can be independent of the value at which the truth-tellers reach asymptotic relative consensus in simple cases of the interaction function $P$.

\begin{figure}[htb!]
    \centering
    \begin{subfigure}[b]{0.9\textwidth}
        \centering
        \includegraphics[width=\textwidth]{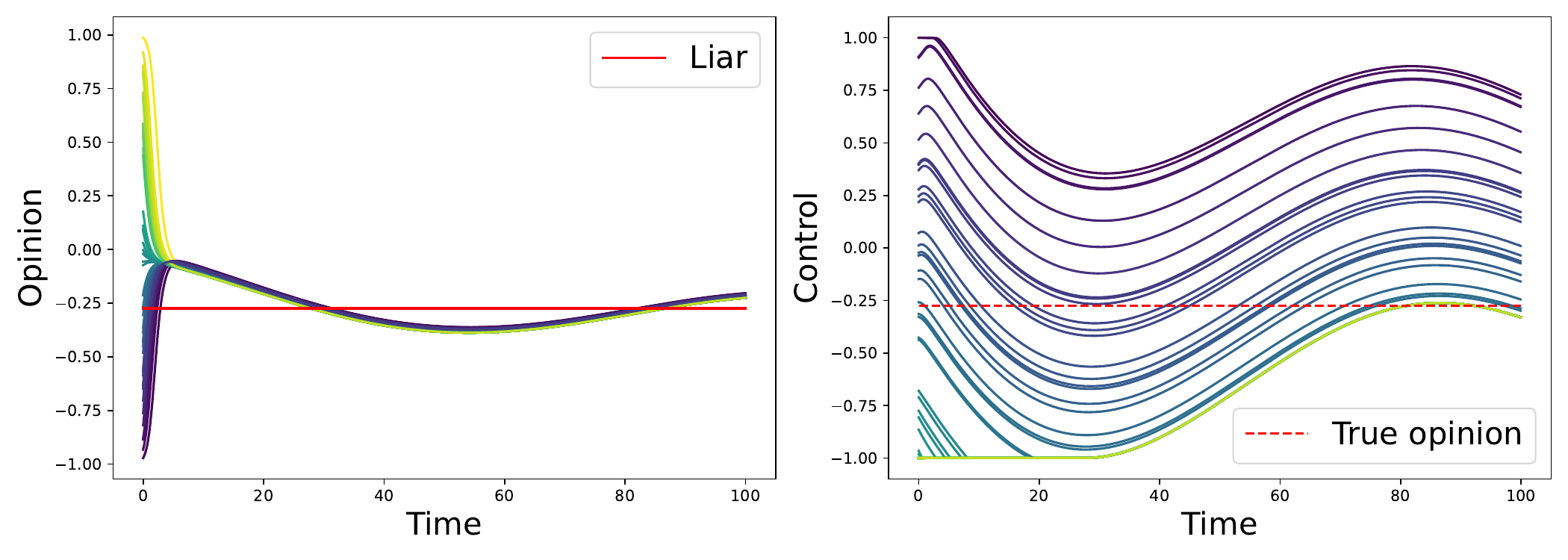}
        \caption{Classic regularisation control}
        \label{fig:timecon-reg-std-init}
    \end{subfigure}
    \begin{subfigure}[b]{0.9\textwidth}
        \centering
        \includegraphics[width=\textwidth]{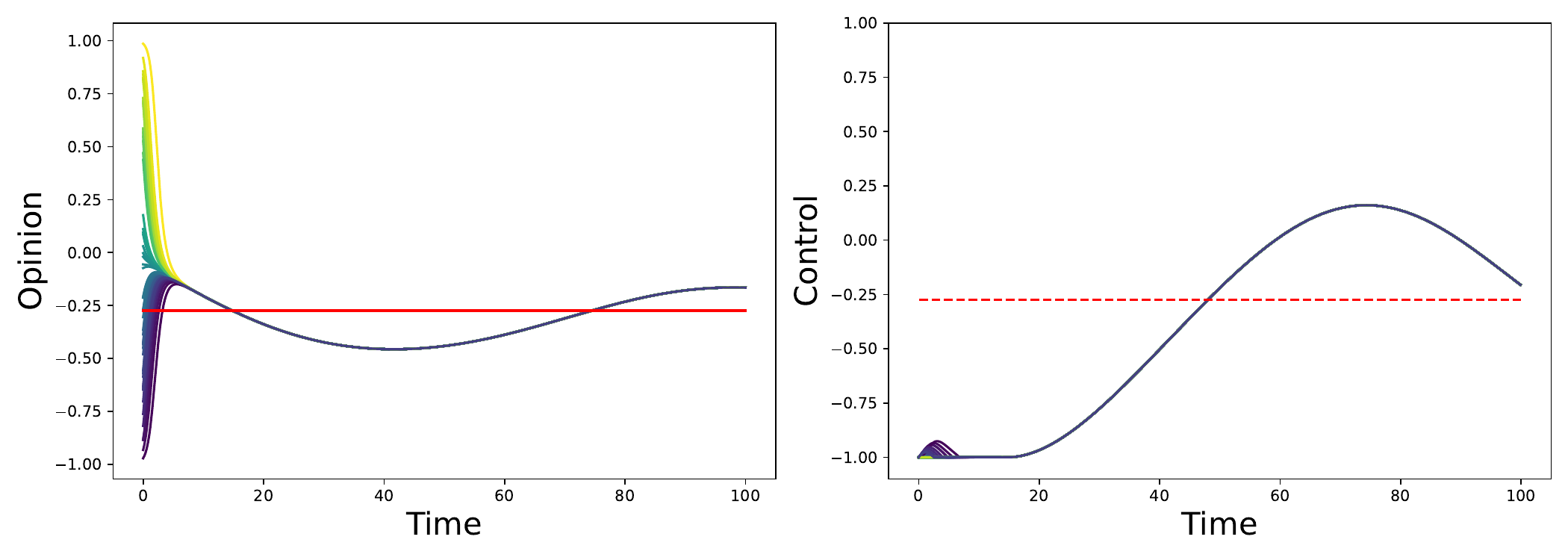}
        \caption{Consistent control}
        \label{fig:timecon-reg-consis-init}
    \end{subfigure}
    \caption{Trajectory of agents and the corresponding controls in the case where the liar is penalised for expressing opinions inconsistently in time. Here, $N=50$, $P(x,\cdot\,) = P(x) = 1-x^2$, $\nu = 0.001$ and we compare two different regularisations as the initial condition $y_i^0$. Plot \textbf{(a)} shows the case where the liar is initially penalised for straying from their true opinion, as in Section \ref{sec:Micro-model-std-reg}. Plot \textbf{(b)} shows the case where the liar is initially penalised for presenting different opinions to different agents, as in Section \ref{sec:Micro-model-consis-reg}.}
    \label{fig:timecon-reg-different-init}
\end{figure}

Figure \ref{fig:timecon-reg-different-init} demonstrates the effects of choosing different regularisations for the initial control $y_i^0$. Figure \ref{fig:timecon-reg-std-init} shows an example where the initial control is chosen according to our classic regularisation (\ref{eq:SR-control}) and Figure \ref{fig:timecon-reg-consis-init} shows the same example but where the initial control is chosen according to the consistent regularisation (\ref{eq:CP-control}). We observe familiar oscillations in both plots but the initial value $y_i^0$ makes a large difference to the overall dynamics and success of the control. It would be possible to choose a more optimal trajectory for the liar under this regularisation without setting $y_i^0$ arbitrarily, for example through dynamic programming, however we defer this for future work.

\subsection{A variant on time-consistent control}\label{sec:Micro-model-var-timecon}

One could argue that rather than there being a cost associated with changing opinion, as in Section \ref{sec:Micro-model-timecon}, there should be a cost associated with changing opinion in a way that is unlike the changing of opinion of truth-telling agents. Here, we assume that the truth-telling agents have some knowledge of the opinion formation process and define the `expected' change in opinion of the liar, according to an observing agent $i$, by
$$ \dot{\tilde{x}}_{1,i} = \frac{1}{N}\sum_{j=2}^N P(\tilde{x}_{1,i}, x_j)(x_j - \tilde{x}_{1,i}), \quad \tilde{x}_{1,i}(0) = y_i^0,$$
where $y_i^0$ is the initial lie told to agent $i$. Then $y_i$ is chosen to minimise a cost functional 
\begin{equation}\label{eq:TC-var}
    C(y) = \int_0^T \left(\frac{1}{N}\sum_{i=1}^N \frac{1}{2}(x_i - x_d)^2 + \frac{1}{N}\sum_{i=2}^N \frac{\nu}{2}\left(\frac{dy_i}{dt}- \dot{\tilde{x}}_{1,i}\right)^2\right)\, ds.
\end{equation}
As in Section \ref{sec:Micro-model-timecon}, at each time-step $n\geq1$, $\frac{dy_i}{dt}$ is set to be the forward Euler approximation to the time derivative of $y_i$. Furthermore, at each time-step $n\geq1$, we approximate the `expected' derivative $\dot{\tilde{x}}_{1,i}^n$ by
\begin{equation}\label{eq:xtilde-person-i}
    \dot{\tilde{x}}_{1,i}^n =\frac{1}{N}\sum_{j=2}^NP(y_i^{n-1}, x^{n-1}_j)(x^{n-1}_j - y_i^{n-1}).
\end{equation}
Hence the reduced optimisation problem for $y^n = (y_2^n, ..., y_N^n)^T$ is
\begin{equation}\label{eq:var-time-con-prob}
    y^n = \arg\min_{u\in{\mathcal{I}^{N-1}}}\left\{\int_{t_n}^{t_{n+1}} \left(\frac{1}{N}\sum_{i=1}^N \frac{1}{2}(x_i - x_d)^2 + \frac{1}{N}\sum_{i=2}^N\frac{\nu}{2}\left(\frac{y_i - y_i^{n-1}}{\Delta t} - \tilde{x}_{1,i}^n\right)^2\right)\right\}, \quad n\geq 1,
\end{equation}
where $\tilde{x}_{1,i}^n$ is given by equation (\ref{eq:xtilde-person-i}). As an initial condition, we can take $y_i^0$ to be given by the solution to equation (\ref{eq:noreg_optcontrol}), (\ref{eq:SR-control}) or (\ref{eq:CP-control}). Our optimal control is derived through Hamiltonian techniques discussed in Section \ref{sec:Micro-model-no-reg} and detailed in Appendix \ref{app:optimality-conditions}. The resulting control for $n\geq1$ is given by a projection onto $\mathcal{I}$ of $\tilde{y}_i^n$ where $\tilde{y}_i^n$ satisfies
\begin{equation}\label{eq:control_likexi}
    \frac{\tilde{y}_i^n - y_i^{n-1}}{\Delta t} - \frac{1}{N}\sum_{j=2}^N P(y_i^{n-1}, x_j^{n-1})(x_j^{n-1} - y_i^{n-1}) = -\frac{\Delta t^2}{\nu N}(x_i^{n+1}- x_d) \partial_{\tilde{y}_i}\{P(x_i, \tilde{y}_i)(\tilde{y}_i-x_i)\}|_{(\bar{x}_i, \tilde{y}_i^n)}.
\end{equation} 
We can see that taking $\nu\to\infty$ means that the liar's apparent opinions update in the same way as they would if they were the true opinions of a truth-telling agent. 

\begin{figure}[htb!]
    \centering
    \begin{subfigure}[b]{0.45\textwidth}
        \centering
        \includegraphics[width=\textwidth]{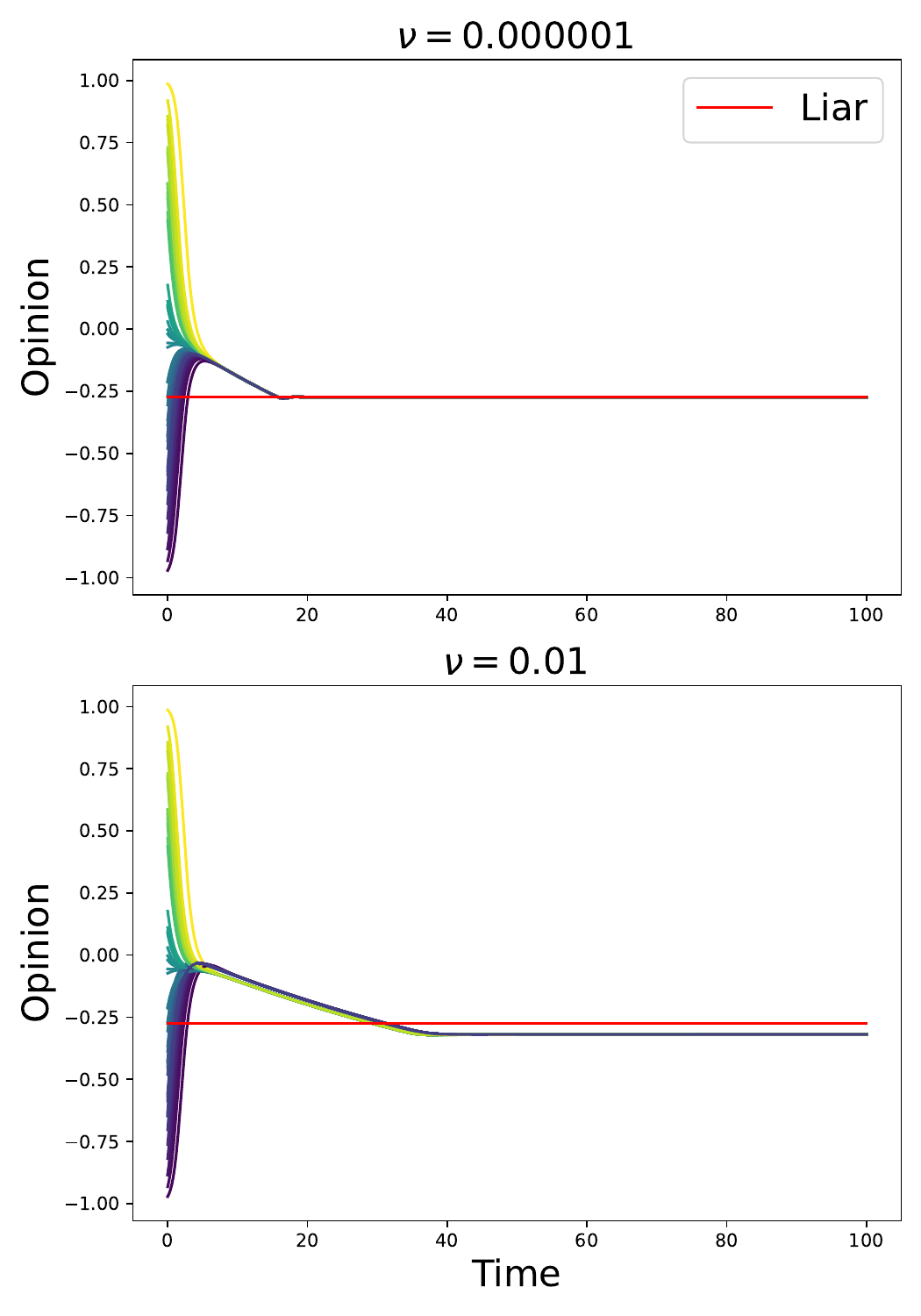}
        \caption{Trajectory}
    \end{subfigure}
    \begin{subfigure}[b]{0.45\textwidth}
        \centering
        \includegraphics[width=\textwidth]{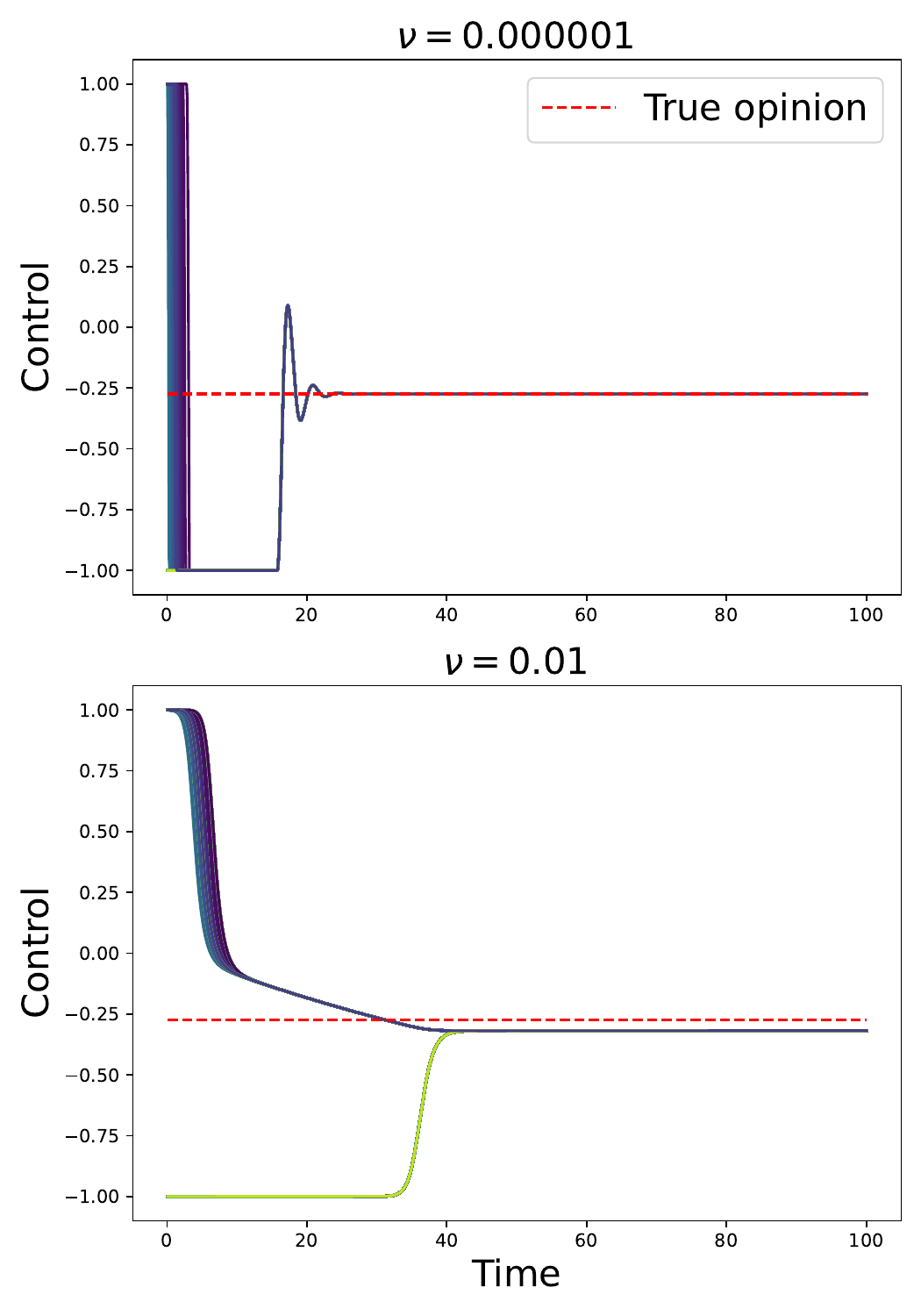}
        \caption{Control}
    \end{subfigure}
    \caption{Trajectory of agents and the corresponding controls in the case where the liar is penalised for changing their opinion in a way which is different to that of a standard agent. Here, $N=50$, $P(x,\cdot\,) = P(x) = 1-x^2$, $y_i^0$ is chosen according to (\ref{eq:noreg_optcontrol}) and we compare values of $\nu = 0.000001$ and $\nu = 0.001$.}
    \label{fig:timeconsis-reg-likeagent}
\end{figure}

Figure \ref{fig:timeconsis-reg-likeagent} shows a simple example of the trajectory and corresponding controls of a system with a liar acting according to this regularisation with the initial condition of the control given by (\ref{eq:noreg_optcontrol}). In contrast to the previous sections, we see that for sufficiently large values of $\nu$, the regularisation prohibits consensus at the liar's goal opinion. This is due to the fact that consensus is a fixed point of the truth-telling dynamics \cite{motsch2014heterophilious}. When the apparent opinion of the liar is in agreement with the consensus of the truth-telling agents, the preference of the liar to act like a truth-telling agent and present an opinion in agreement with the consensus overrides the liar's goal to bring consensus to $x_d$. Again, it is important to point out that the success of the liar's strategy is highly dependent on the initial condition on the control, $y_i^0$.

\subsection{Sparse control}

Finally, we wish to minimise the number of people to which the liar has to lie. If the liar has true opinion $x_d$ then the amount by which they lie to agent $i$ is $|y_i - x_d|$. Hence, while the dynamics evolve according to (\ref{eq:xidyn}), we seek to minimise the cost functional
\begin{equation}\label{eq:SC-prob}
    C(y) = \int_0^T \left(\frac{1}{N}\sum_{i=1}^N \frac{1}{2}(x_i - x_d)^2 + \frac{1}{N}\sum_{i=2}^N \nu|y_i-x_d|\right)\, ds.
\end{equation}

\subsubsection{Nonlinear model predictive control}

The optimisation problem (\ref{eq:SC-prob}) is analytically intractable due to the non-differentiability of the $l_1$ norm, therefore we cannot apply the Hamiltonian method described in Section \ref{sec:Micro-model-no-reg}. We will instead use a numerical method known as \textit{nonlinear model predictive control} (NMPC) \cite{mayne2000constrained, borghi2024kinetic, bailo2018optimal}. We consider a prediction horizon of $H$ steps for $H\in\{1, ..., M - 1\}$, a discrete time version of the dynamics (\ref{eq:xidyn}) and minimise the following cost function
\begin{equation}\label{eq:sparse-cost-NMPC}
    C_n(y) = \sum_{h=0}^H\frac{\Delta t}{N}\left(\sum_{i=1}^N\frac{1}{2}(x_i^{n+h+1}-x_d)^2 + \sum_{i=2}^N\nu|y_i^{n+h} - x_d|\right).
\end{equation}
The method produces a sequence of $H+1$ controls $(y^n, ..., y^{n+H})$ where only the first term $y^n$ is used to advance the dynamics from time-step $n$ to $n+1$. The full nonlinear model predictive control method, with particle swarm optimisation (PSO) used to find the minimiser can be found in Algorithm \ref{alg:sparse-nummethod} in Section \ref{app:nmpc-algorithm} of the Appendix. The optimisation algorithm is based on Algorithm 3 in \cite{borghi2024kinetic}.

\begin{figure}[htb!]
    \centering
    \begin{subfigure}[b]{0.495\textwidth}
        \centering
        \includegraphics[width=\textwidth]{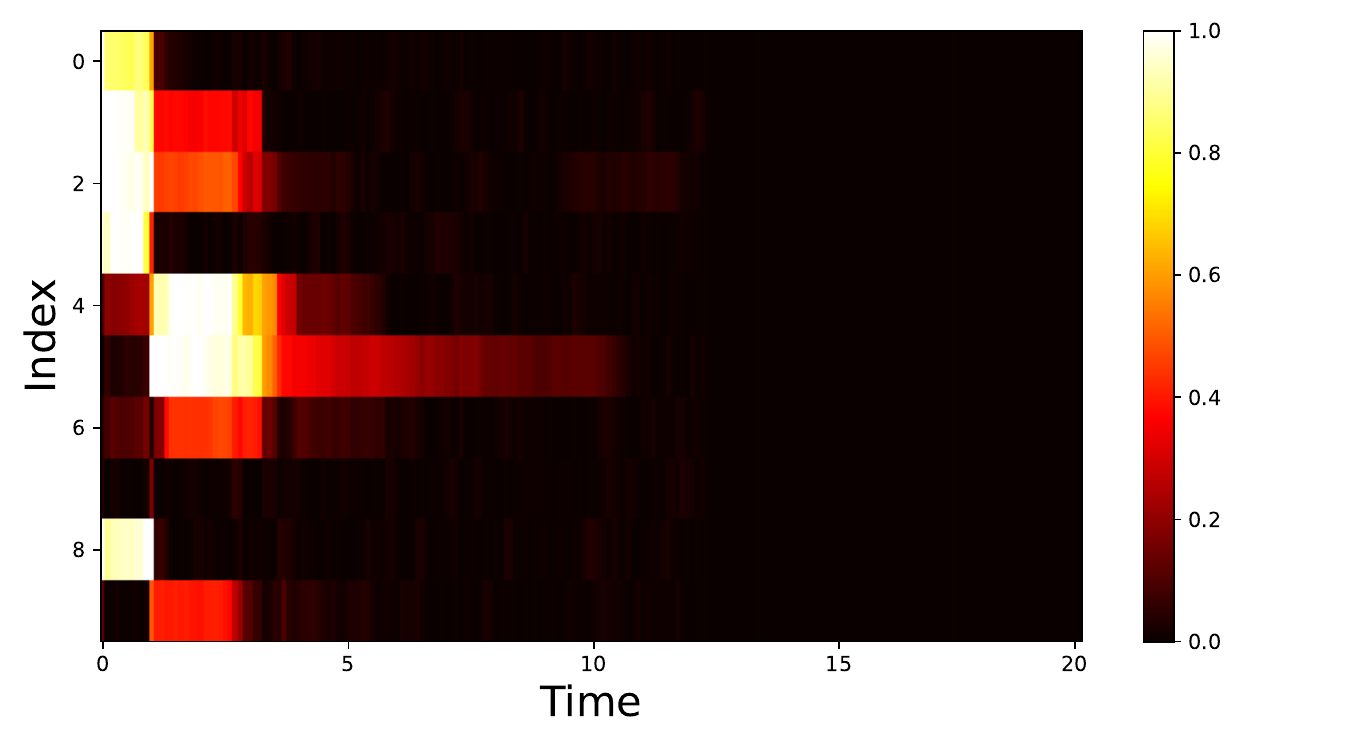}
        \caption{$l_1$ norm regularisation with $P(x,\cdot) = P(x) = 1-x^2$.}
        \label{subfig:l1norm-NMPC}
    \end{subfigure}
    \begin{subfigure}[b]{0.495\textwidth}
        \centering
        \includegraphics[width=\textwidth]{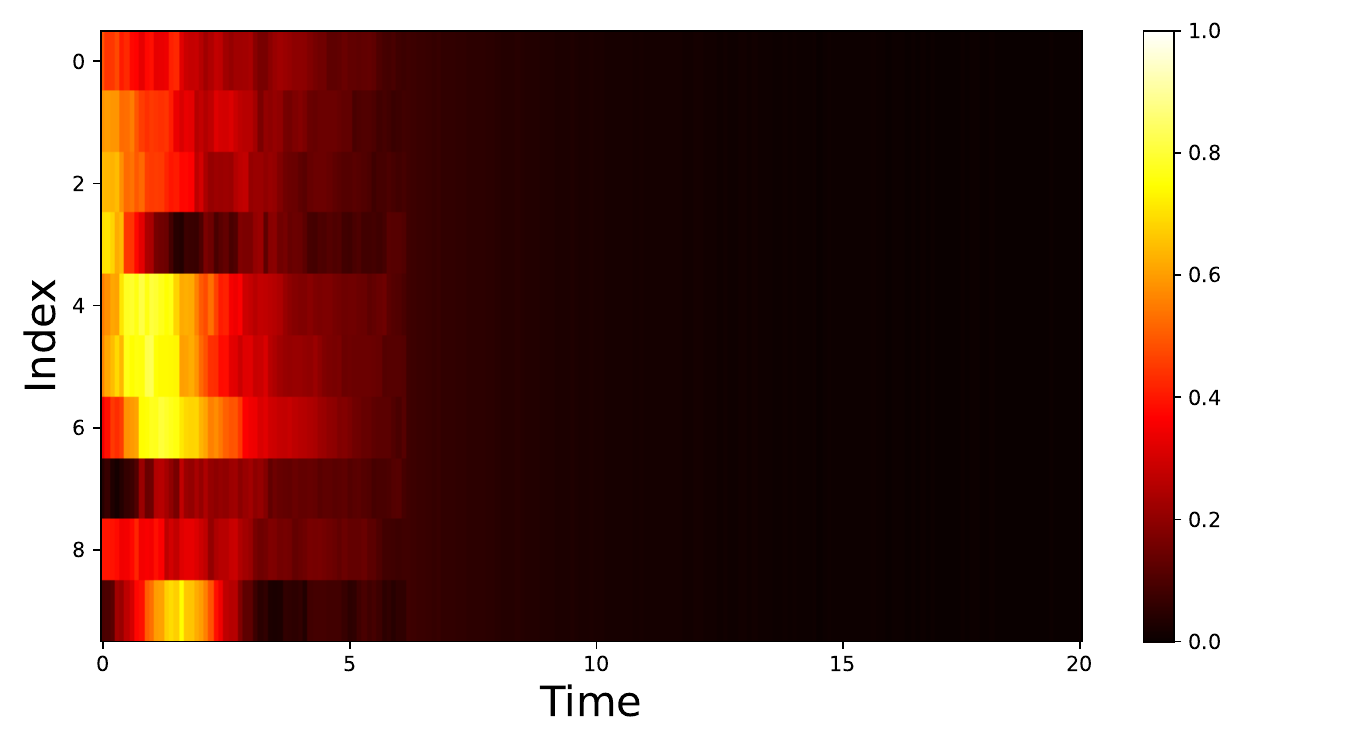}
        \caption{$l_2$ norm regularisation with $P(x,\cdot) = P(x) = 1-x^2$.}
        \label{subfig:l2norm-NMPC}
    \end{subfigure}
    \begin{subfigure}[b]{0.495\textwidth}
        \centering
        \includegraphics[width=\textwidth]{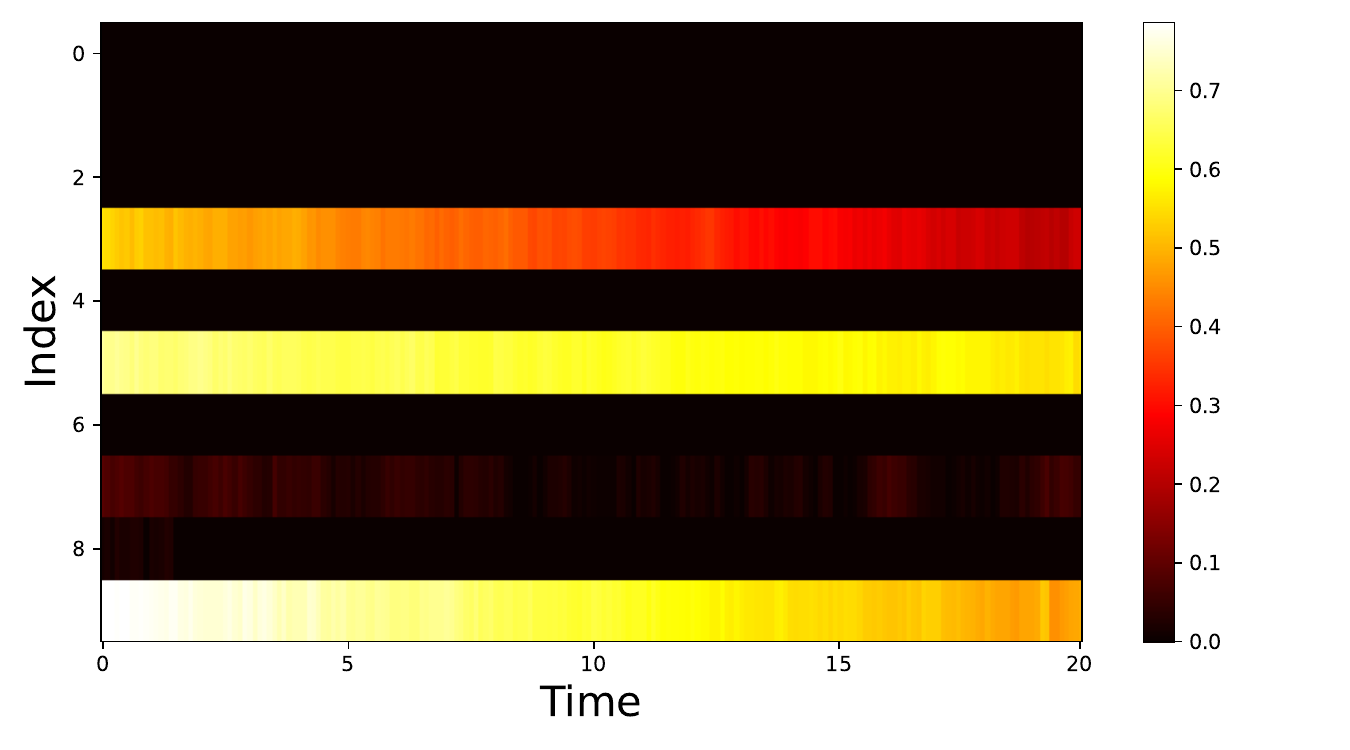}
        \caption{\centering $l_1$ norm regularisation with $P(x,y) = \chi(|y-x|<R)$ for $R=0.2$.}
        \label{subfig:l1norm-NMPC-boundedcon}
    \end{subfigure}
    \begin{subfigure}[b]{0.495\textwidth}
        \centering
        \includegraphics[width=\textwidth]{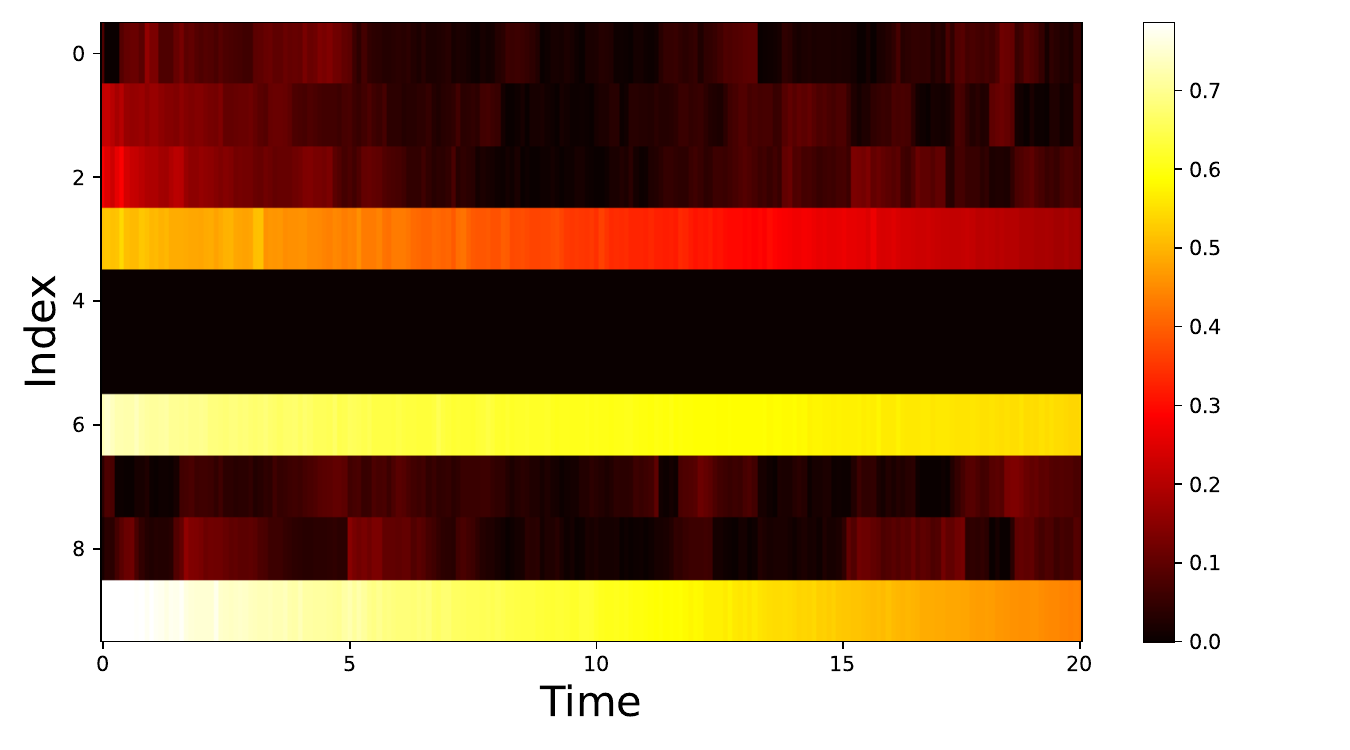}
        \caption{\centering$l_2$ norm regularisation with $P(x,y) = \chi(|y-x|<R)$ for $R=0.2$.}
        \label{subfig:l2norm-NMPC-boundedcon}
    \end{subfigure}
    \caption{Comparison of the size of the $l_1$ norm of the variation of the control from the liar's true opinion, $x_d$, in the cases where the optimisation problem is solved with penalty incurred for large $l_1$ norm (sparse control) versus large $l_2$ norm. In all cases, the controls are calculated via NMPC-PSO with $H=3$ and we consider $N=10$ truth-telling agents. In sub-figures (a) and (b), the interaction function $P(x,\cdot) = 1-x^2$ whereas in figures (c) and (d) we consider a bounded confidence interaction function $P(x,y) = \chi(|y-x|<R)$ where $\chi$ is an indicator function and $R=0.2$. In figures (a) and (b), we take regularisation parameter $\nu = 0.0155$ whereas in figures (c) and (d) we take $\nu = 0.0225$.}
    \label{fig:sparse-control-nmpc}
\end{figure}

To present an example of sparse optimal control in practice, we will consider a system with only 10 truth-telling agents. Figure \ref{fig:sparse-control-nmpc} displays heat maps of the matrix $\mathcal{H}$, where $\mathcal{H}_{i,n} = |y_i^n - x_d|$ shows the magnitude of the lie told to agent $i$ at time-step $n$. The $y$-axis shows the index of each agent, the $x$-axis shows time and a bright square corresponds to a lie of large magnitude whereas a dark square corresponds to little/no lie told. In this Figure, we consider the resulting controls for the system given two interaction functions. The first is $P(x,\cdot) = P(x) = 1-x^2$, which we have discussed throughout Sections \ref{sec:Micro-model-no-reg} to \ref{sec:Micro-model-var-timecon}, shown in Figures \ref{subfig:l1norm-NMPC} and \ref{subfig:l2norm-NMPC}. The second is the bounded confidence interaction function $P(x,y) = \chi(|y-x|<R)$ for $\chi$ an indicator function and $R=0.2$ in this case, first introduced in equation (\ref{eq:bounded-con-introduciton}) and discussed further in Section \ref{sec:bounded-conf}. The corresponding controls for this interaction function are found in Figures \ref{subfig:l1norm-NMPC-boundedcon} and \ref{subfig:l2norm-NMPC-boundedcon}.

Figures \ref{subfig:l1norm-NMPC} and \ref{subfig:l1norm-NMPC-boundedcon} show the results of our NMPC-PSO algorithm acting on the cost function (\ref{eq:sparse-cost-NMPC}) where we refer to the regularisation $|y_i^{n+h}-x_d|$ as a regularisation on the $l_1$ norm. Figures \ref{subfig:l2norm-NMPC} and \ref{subfig:l2norm-NMPC-boundedcon} show the result of the same algorithm but with the $l_1$ norm constraint on the control replaced with an $l_2$ constraint, $|y_i^{n+h}-x_d|^2$ in (\ref{eq:sparse-cost-NMPC}). The 10 truth-telling agents have random initial opinions in the interval $[-1,1]$, given in Table \ref{tab:table1}, and the liar has goal opinion $x_d =0$. Hence, in this case, the brightest spot in the heat map (the maximum absolute value of lie that the liar can tell) is 1. Comparing Figures \ref{subfig:l1norm-NMPC} and \ref{subfig:l2norm-NMPC}, we observe sparser results in the $l_1$ case, Figure \ref{subfig:l1norm-NMPC}, where the liar tells large lies to 2-4 of the agents at each time-step. In contrast, in Figure \ref{subfig:l2norm-NMPC}, where the liar acts under the $l_2$ norm regularisation, the liar tells lots of smaller lies to every agent at each time-step. Similar results are obtained in Figures \ref{subfig:l1norm-NMPC-boundedcon} and \ref{subfig:l2norm-NMPC-boundedcon} for the bounded confidence case.

\begin{table}[htb!]
  \begin{center}
    \begin{tabular}{l|c|c|c|c|c|c|c|c|c|c}
      \textbf{index} & 0 & 1 & 2 & 3 & 4 & 5 & 6 & 7 & 8 & 9\\ 
      \hline
      $x_i^0$ & -0.275 & -0.387 & -0.431 & 0.716 & -0.883 & -0.885 & -0.932 & 0.177 & -0.208 & 0.986\\ 
    \end{tabular}
    \caption{Initial opinions of truth-telling agents being controlled in Figure \ref{fig:sparse-control-nmpc} to 3 decimal places}
    \label{tab:table1}
  \end{center}
\end{table}

When $P(x, \cdot)=1-x^2$, it is natural for the liar to lie the most at the start of the simulation. This is because the model is contracting, so the truth-telling agents are the furthest away from the liar at initial times. This fact is reflected in both Figures \ref{subfig:l1norm-NMPC} and \ref{subfig:l2norm-NMPC} where at large times we see that the control is switched off completely whereas at early times we use quite a lot of control for both regularisations. 

A feature of the model that we exploit in our numerical scheme is that when the truth-tellers are in agreement with each other, the liar has no incentive to lie more to any one truth-teller over another. Hence in this limit the liar presents the same lie to every truth-telling agent and the size of the problem can be reduced to 2 agents, the truth-teller (with a weighting of $N-1$) and the liar. This reduces the computational complexity of the algorithm at later times and gives more accurate results.

The random initialisation of truth-telling agents is such that the agent with the most extreme initial opinion (initial opinion the furthest away from $x_d$) has index 9 (see Table \ref{tab:table1}). The initial mean opinion of the truth-telling agents is $\approx-0.212$. Interestingly, in the case where $P(x,\cdot) = 1-x^2$, we see that in Figure \ref{subfig:l1norm-NMPC} the liar spends relatively little time lying to agent 9, despite this agent having the opinion furthest from the goal opinion $x_d$. The liar's strategy seems to be to let the other truth-telling agents do the hard work of bringing agent 9 toward the truth-telling consensus, which happens to be further away from agent 9's initial opinion than $x_d$ is.

When $P(x,y) = \chi(|y-x|<R)$, we observe very different strategies from the liar compared to the case where $P(x,\cdot) = 1-x^2$. In the sparse case, Figure \ref{subfig:l1norm-NMPC-boundedcon}, the liar's strategy is to target the truth-telling agents that are the furthest away (or most extreme) by choosing to lie the most to the agent(s) closest to them who is also connected to an extreme agent. In this way, the liar reduces the cost from the $l_1$ norm constraint while still manipulating the most extreme agent through the less extreme agents that they are connected to. Indeed, we see that in Figure \ref{subfig:l1norm-NMPC-boundedcon} the liar focuses their efforts on influencing agents with indices 3, 5 and 9. From Table \ref{tab:table1}, we see that these agents correspond well to the strategy outlined above. In contrast, we see that in Figure \ref{subfig:l2norm-NMPC-boundedcon} when we penalise the liar according to the $l_2$ norm, the liar targets a similar set of agents but chooses to lie to agent 6 who is further away than agent 5. There are also a lot of lies of small magnitude told to agents that are close to the liar that are absent in the case of sparse control. For both cases, we see that the liar tells lies for a much longer time than it does when $P(x,\cdot) = P(x) = 1-x^2$. This is due to the fact that agents are very poorly connected in the bounded confidence case so it takes much longer to reach consensus.

A difficulty with using an NMPC-PSO method to calculate $y_i^n$ in the case where $P(x,y)$ is a bounded confidence kernel is that under this kernel, the optimal control can be non-unique. Indeed, in the case where $|x_i-x_d|>R$, there is a trade-off between telling a large lie to agent $i$ which would be costly or not lying at all and therefore not effecting the dynamics of $x_i$ which also invokes a cost. Hence, due to the inherent stochasticity of the algorithm, we can recover quite different solutions for the same parameter values $\nu$ and $H$.

In terms of choosing an appropriate value of $H$, there is a trade-off between the information given to the algorithm and the computational cost. Indeed, for small $H$, not many iterations of PSO are necessary to find a good minimiser for (\ref{eq:sparse-cost-NMPC}) but the liar has very little information about the trajectory of agents so generally will not produce a very sparse result. In contrast, when $H$ is high, the liar has a lot of information so will generally produce sensible, sparse results. The drawback, however, is that the PSO takes many more iterations to converge since $y=(y^n, ..., y^{n+H})$ is of high dimension, increasing the computational cost dramatically. Another impact of choosing a large value of $H$, observed for the case when $P(x,\cdot) = 1-x^2$, is that the algorithm learns that the model is contracting and so optimal strategy is to lie a lot at the beginning of the simulation and then not at all later in the simulation. While holding $\nu$ constant, this can produce results that are sparse in time but not in space since the liar learns to lie a lot to everyone at initial time.

A slightly different sparse optimal control problem that is also worth considering is to solve for the optimal control when the liar can only lie to a maximum of $N_L<N$ individuals. This approach is known as \textit{cardinality constrained optimisation} \cite{tillmann2024cardinality, kanzow2021sequential}. Another strategy to promote sparse results would be to impose that the liar can only tell lies for a certain number of time-steps. These formulation are left for future work.

\section{Boltzmann-type control}\label{sec:boltzmann}

In a regime where our population size is very large, it can be extremely computationally costly to solve the large number of ODEs and calculate the many individual controls necessary for the microscopic problem outlined in Section \ref{sec:Micro-model}. We will now consider a Boltzmann-type dynamic corresponding to the instantaneous control formulation described in Section \ref{sec:Micro-model}. This allows us to study the behaviour of the system in the limit as the number of agents $N\to \infty$. The Boltzmann equation is a classical model for describing the macroscopic behaviour of interacting gas particles \cite{livi2017nonequilibrium}. The idea behind using the Boltzmann equation for our setting is that we can liken the collisions between gas particles to the interactions taking place between agents of different opinions in the limit $N\to \infty$ \cite{pareschi2013interacting, toscani2006kinetic}.

In this Section, we will only consider the kinetic description of the system when the liar is acting in accordance with the standard regularisation control case discussed in Section \ref{sec:Micro-model-std-reg}. Recall the cost function
$$ C(y) = \int_0^T\left(\frac{1}{N}\sum_{i=1}^N\frac{1}{2}(x_i - x_d)^2 + \frac{1}{N}\sum_{i=2}^N\frac{\nu}{2}(y_i - x_d)^2\right)\, ds.$$
We will make the additional assumption that $P(x_i,\cdot\,) = P(x_i)$ i.e., that the interaction between agents $i$ and $j$ only depends on the opinion of agent $i$, similarly to the examples discussed throughout Section \ref{sec:Micro-model}. This assumption means that we can write an explicit equation for $\tilde{y}_i^n$. The more complex case where $P(x_i, x_j)$ is the bounded confidence kernel (\ref{eq:bounded-con-introduciton}) will be discussed in Section \ref{sec:bounded-conf}.

Given our assumptions and recalling equation (\ref{eq:SR-control}), at each time-step, $y_i^n$ is given by a projection of the solution to the equation
\begin{equation*}\label{eqn:simplified-binary-control}
    \tilde{y}_i^n = x_d - \frac{\Delta t}{\nu N}P(x_i^n)(x_i^{n+1}-x_d)
\end{equation*}
onto the interval $\mathcal{I}$. We write $x_i^{n+1}$ as one Euler step so
\begin{equation}\label{eq:eulerstep}
    x_i^{n+1} = x_i^n + \frac{\Delta t}{N}\left(P(x_i^n, \tilde{y}_i^n)(\tilde{y}_i^n - x_i^n)+\sum_{j=2}^N P(x_i^n, x_j^n)(x_j^n - x_i^n)\right),
\end{equation}
for $i = 2, ..., N$. We now consider the model predictive control system in the simplified case of binary interactions. We have that when two truth-tellers interact their opinions update as
\begin{equation}\label{eq:binary-interactions-boltzmann}
    \begin{split}
        x_i^{n+1} &= x_i^n + \frac{\Delta t}{2}P(x_i^n)(x_j^n - x_i^n),\\
    x_j^{n+1} &= x_j^n + \frac{\Delta t}{2}P(x_j^n)(x_i^n - x_j^n).
    \end{split}
\end{equation}
When an arbitrary truth-teller $k$ interacts with the liar (agent 1), the opinions of agent $k$ and the liar update as follows:
\begin{align}\label{eqn:x_i-bin-interact}
    x_k^{n+1} &= x_k^n + \frac{\Delta t}{2}P(x_k^n)(y_k^n - x_k^n),\\
    x_1^{n+1} &= x_1^n = x_d ,\nonumber
\end{align}
where $y_k^n$ is the lie told by the liar to person $k$ under the binary interaction dynamics,
\begin{equation}\label{eq:bin-control}
    \tilde{y}_k^n = x_d - \frac{\Delta t}{2\nu}P(x_k^n)(x_k^{n+1} - x_d).
\end{equation}
It is worth comparing the equations for the control in the binary interaction dynamic (\ref{eq:bin-control}) and the microscopic model (\ref{eq:SR-control}). We see that in considering only binary interactions, the liar loses knowledge of how the truth-tellers in the population will effect the opinion $x_i^{n+1}$. This generally means that the liar is less effective in their goal to quickly achieve a consensus at $x_d$. 

\begin{figure}[htb!]
    \centering
    \begin{subfigure}[b]{0.45\textwidth}
        \centering
        \includegraphics[width=\textwidth]{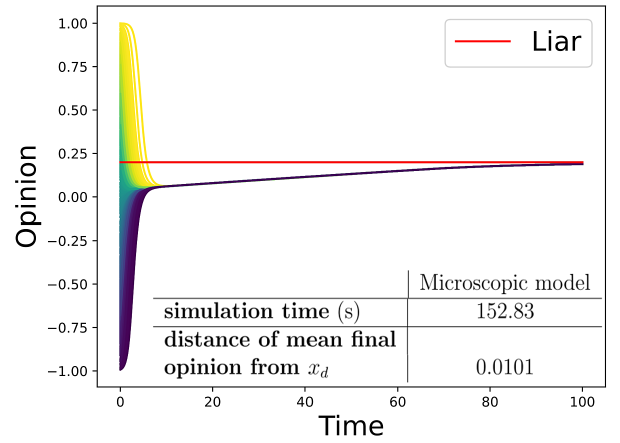}
        \caption{Microscopic model}
    \end{subfigure}
    \begin{subfigure}[b]{0.45\textwidth}
        \centering
        \includegraphics[width=\textwidth]{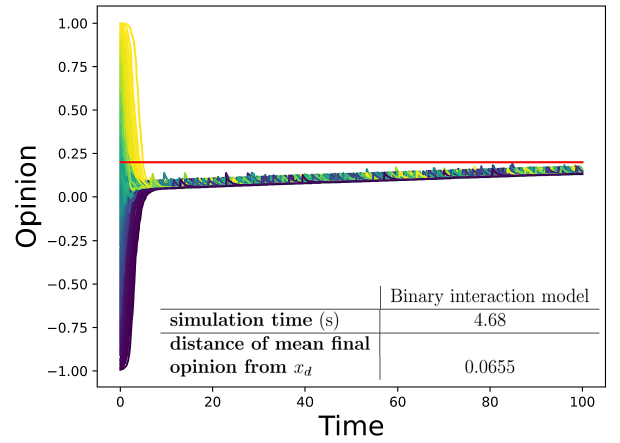}
        \caption{Binary interaction model}
    \end{subfigure}
    \caption{Comparison of trajectories of truth-telling agents moving under microscopic dynamics (a) and binary interaction dynamics (b). The microscopic dynamics are given by equation (\ref{eq:xidyn}) and control (\ref{eq:SR-control}) and the binary interaction dynamics are given by equations (\ref{eq:binary-interactions-boltzmann}), (\ref{eqn:x_i-bin-interact}) and control (\ref{eq:bin-control}). We display the time taken to solve the for the dynamics up to $T=100$ as well as the mean distance from the desired opinion at time $T$. Here, we take $N=500$, $P(x,\cdot) = 1-x^2$ and $\nu = 0.00001$. The goal opinion of the liar is shown in red.}
    \label{fig:comparison-microscopic-binary}
\end{figure}

To explore the differences between the microscopic and binary interaction models further, we consider a population of size $N=500$. Figure \ref{fig:comparison-microscopic-binary} shows a simple example of a simulation when $P(x,\cdot) = 1-x^2$ for both the microscopic model and the binary interaction model. We simulate binary interactions in a probabilistic manner. At each time-step, truth-tellers will interact with another truth-teller (chosen at random) with probability $\frac{N-1}{N}$ and will update their opinions according to equation (\ref{eq:binary-interactions-boltzmann}) or interact with the liar with probability $\frac{1}{N}$ and update their opinion according to equation (\ref{eqn:x_i-bin-interact}). The first difference we notice in the two models is that the computational cost of the microscopic model is orders of magnitude higher than the binary interaction model. We see that the dynamics of the microscopic model are largely well approximated by the dynamics of the binary interaction model but that the control is indeed less effective, evidenced by the distance between $x_d$ and the mean opinion of truth-telling particles at time $T$, $\left|\sum_{j=2}^Nx_j(T) - x_d\right|$, being greater for the binary interaction model. Overall however, we conclude that the loss in accuracy is small in comparison to the gain we achieve in computational cost when we consider only binary interactions in the large population limit $N\gg 1$.

In order to simplify our dynamics for the binary interaction model, we will determine a condition on $\nu$ such that $\tilde{y}_i^n\in\mathcal{I}=[-1,1]$ and hence $y_i^n = \tilde{y}_i^n$.
\proposition{If $x_d\in(-1,1)$, the projected opinion $y_k^n$ of the liar to person $k$ at time-step $n$ is given by equation (\ref{eq:bin-control}) if
$$ \nu\geq \max_{x\in\mathcal{I}} \left\{\frac{\alpha P(x)}{1-x_d}\left(-(x - x_d) - \alpha P(x)(1-x)\right)\right\}, \quad \nu\geq \max_{x\in\mathcal{I}}\left\{\frac{\alpha P(x)}{1+x_d}((x-x_d) -\alpha P(x)(1+x))\right\}$$
where $\alpha = \frac{\Delta t}{2}$. These conditions can be combined, assuming $\mathcal{I} = [-1,1]$ and $0\leq P(x)\leq 1$, into a simpler form:
\begin{equation}\label{eq:ineq-nu-bound}
    \nu\geq\frac{\alpha(1+|x_d|) }{1-|x_d|}, \quad \text{for }|x_d|<1.
\end{equation}
When $|x_d| = 1$, we must have $\nu = \infty$.}\label{prop:bound-on-reg} \normalfont

\textbf{Proof}: First, we will determine conditions for $y_k^n\leq 1$. Combining equations (\ref{eq:bin-control}) and (\ref{eqn:x_i-bin-interact}), we have that
$$ y_k^n = \frac{\nu+ \alpha P(x_k^n)}{\nu + \alpha^2 P(x_k^n)^2}x_d - \frac{\alpha P(x_k^n)(1-\alpha P(x_k^n))}{\nu + \alpha^2P(x_k^n)^2}x_k^n, $$
where $\alpha = \frac{\Delta t}{2}$. Setting the inequality $y_k^n\leq 1$ and rearranging, we have
$$ \alpha P(x_k^n)x_d - \alpha^2P(x_k^n)^2 - \alpha P(x_k^n)(1-\alpha P(x_k^n))x_k^n\leq \nu(1-x_d).$$
Assuming that $x_d\in(-1,1)$, we can divide through by $(1-x_d)$ and maintain the inequality, giving
\begin{equation}\label{eq:nubound}
    \nu \geq \frac{\alpha P(x_k^n)}{1-x_d}(-(x_k^n-x_d) - \alpha P(x_k^n)(1-x_k^n)).
\end{equation}
Taking the maximum over the RHS of (\ref{eq:nubound}) gives the first inequality in Proposition \ref{prop:bound-on-reg}. The second is given by the same argument but requiring $y_k^n\geq -1$.

The simplified version of the two inequalities (\ref{eq:ineq-nu-bound}), is obtained by noting that the maximum and minimum values that $x$ can take are 1 and $-1$ respectively and that the maximum value that $P(x)$ can take is 1.

When $|x_d| = 1$, we must take $\nu = \infty$ (i.e., the liar must tell the truth) since in this regime where $P(x_k, y_k) = P(x_k)$, the only incentive to lie is to project an opinion more extreme than $\pm 1$. Indeed, if the liar projects an opinion closer to $x_k$ than $x_d$, then the rate of change of $x_k$ in the direction of $x_d$ will be decreased. $\blacksquare$

\remark{We note that the bound $y_k^n\leq 1$ is most likely to be violated when $x_k^n <x_d$, in which case the liar will try to compensate their difference in opinion by projecting an opinion more extreme than 1 in order to cause the greatest rate of change. The inequality (\ref{eq:nubound}) reflects the intuition explained above by the fact that the RHS is only positive, hence placing a limit on what control can be used, when $(x_k^n - x_d)$ is sufficiently negative.}\normalfont

Given Proposition \ref{prop:bound-on-reg}, we can now safely take equation (\ref{eq:bin-control}) to be our control by choosing $\nu$ to be sufficiently large for each $x_d$. Solving the simultaneous equations given by (\ref{eq:bin-control}) and (\ref{eqn:x_i-bin-interact}), we have that the post-interaction opinion of a truth-teller who has interacted with a liar is
\begin{equation}\label{eq:interliarup}
    x_{k}^{n+1} = x_k^n + \frac{\alpha P(x_k^n)(\nu + \alpha P(x_k^n))}{\nu + \alpha^2P(x_k^n)^2}(x_d - x_k^n).
\end{equation}

\subsection{Binary interaction dynamic}\label{sec:bin-interact}

We introduce a density for truth-telling agents $f_T(x,t)$, representing the density of agents with opinion $x\in\mathcal{I}$ at time $t\geq0$. Given a subdomain $\Omega\subset \mathcal{I}$, the integral 
$$ \int_{\Omega} f_T(x,t) \, dx$$
represents the number of truth-tellers with opinion included in $\Omega$ at time $t$. We normalise the density of the truth-tellers' opinions to 1, so
$$ \int_{\mathcal{I}}f_T(x,t )\, dx = 1, \quad \forall t\geq 0.$$
Similarly, we introduce a density for the liar, $f_L(x,t)$. We will impose that
$$\int_{\mathcal{I}}f_L(x,t)\, dx = \rho\leq1, \quad \forall t\geq 0,$$
where the constant $\rho>0$ represents the relative influence of the liar compared with the influence the truth-tellers have on each other. Since the true opinion of the liar is fixed, we have an exact form for $f_L(x,t)$,
$$ f_L(x,t) = \rho\delta(x-x_d), \quad \forall t\geq 0.$$
The rate of change in time of $f_T(x,t)$ depends on the interactions between truth-tellers and the interactions between truth-tellers and the liar. Hence, we have two different types of binary interaction. The first describes the interactions between the truth tellers. Given two truth-tellers with opinions $x, x'\in\mathcal{I}$, their post-interaction opinions $x_*, x_*'$ are given by
\begin{subequations}\label{eq:set-binary-interactions}
\begin{eqnarray}
    x_{*} &=& x + \alpha P(x)(x'-x) +  \Theta_1 D(x),\label{eq:binarystar1}\\
    x'_{*} &=& x' +\alpha P(x')(x-x') + \Theta_2D(x')\label{eq:binarystar2},
\end{eqnarray}
where the additional noise terms $\Theta_1D(x)$ and $\Theta_2D(x)$ take into account effects falling outside of our description of the model. The terms $\Theta_1$ and $\Theta_2$ are random variables that take values on a set $\mathcal{B}\subset\mathbb{R}$, with identical distribution of mean 0 and variance $\sigma^2$. The function $D(x)$ represents the local diffusion coefficient for a given opinion such that $0\leq D(x)\leq 1$. 

The second type of interaction is when a truth teller interacts with the liar. Given a truth-teller with opinion $x''\in\mathcal{I}$, an interaction with the liar gives updated opinion
\begin{equation*}
x''_* = x'' + \alpha P(x'') (y - x'') + \Theta'D(x''),
\end{equation*}
where $y$ is the lie told by the liar and again $\Theta'$ is a random variable identically distributed to $\Theta_1, \Theta_2$. By equation (\ref{eq:interliarup}), we have
\begin{equation}\label{eq:bininterliar}
 x''_* = x'' + \alpha P(x'')\left(\frac{\nu + \alpha P(x'')}{\nu + \alpha^2P(x'')^2}\right)(x_d - x'') + \Theta'D(x'').
 \end{equation}
\end{subequations}

In the absence of diffusion, from equations (\ref{eq:set-binary-interactions}) we see that the mean opinion is not conserved. 

Assuming $0\leq P(x)\leq 1$, in the case of truth-teller to truth-teller interaction, we have
$$ |x_* - x'_*| = |(x-x') - \alpha(x-x')(P(x) + P(x')|\leq |1-2\alpha||x- x'|.$$
Hence, for $\alpha\leq1/2$ the relative distance in opinion between two truth-tellers cannot increase after interaction. Now considering interaction between liars and truth-tellers, we have
$$ |x''_* - x_d| = \left|(x''-x_d) - \alpha P(x)\left(\frac{\nu + \alpha P(x'')}{\nu + \alpha^2P(x'')^2}\right)(x'' - x_d)\right| = \left|\frac{\nu(1 - \alpha P(x))}{\nu + \alpha^2P(x)^2}\right||x'' - x_d|\leq |x'' - x_d|,$$
for $\alpha<1$. This calculation shows us that the relative difference between the truth-teller's opinion and the liar's true opinion cannot increase after interaction.

It is essential to only consider interactions that do not produce values outside of the interval $\mathcal{I}$. In this case, we will take $\mathcal{I}=[-1,1]$ and seek necessary conditions such that the results of binary interactions remain inside $\mathcal{I}$.

\proposition{Assume that $0\leq P(x)\leq 1$, $x_d\in\mathcal{I}$ and
    \begin{equation}
        \quad (1-\alpha)d_-\leq \Theta_i\leq (1-\alpha)d_+, \quad\text{for } i = 1,2, \quad\text{and } \left(1-\frac{\alpha(\nu + \alpha)}{\nu + \alpha^2}\right)d_-\leq \Theta'\leq \left(1-\frac{\alpha(\nu+\alpha)}{\nu + \alpha^2}\right)d_+
    \end{equation}
    where $d_{\pm} = \pm\min_{x\in{\mathcal{I}}}\{(1\mp x)/D(x) , D(x)\neq 0\}$.
    Then the binary interactions (\ref{eq:set-binary-interactions}) preserve their bounds, i.e., $x_*, x_*', x_*''\in\mathcal{I}=[-1,1]$.}\label{prop:bin-inter-conditions}\normalfont

\textbf{Proof}: First, consider the case where we have no noise. The case of truth-teller interactions in equations (\ref{eq:binarystar1}, \ref{eq:binarystar2}) preserve the bounds of $x_*, x'_*\in \mathcal{I}$ since
$$ \max\{|x_*|, |x_*'|\}\leq \max\{|x|, |x'|\}.$$
In the case of a truth-teller's interaction with a liar, it is simple to verify that if $x_d, x''\in\mathcal{I}$ then $x_*''\in\mathcal{I}$.

Now we add noise back into our model. We have that
$$ x_* = x + \alpha P(x)(x'-x) + \Theta_1D(x),$$
so using the fact that $x'\leq 1$, we obtain
$$ x_*\leq x + \alpha P(x)(1-x) + \Theta_1D(x).$$
Setting the above expression to be less than or equal to one, we have
$$ \Theta_1\leq (1-\alpha P(x))\frac{1-x}{D(x)}, \quad D(x) \neq 0.$$
Taking the minimum of the right hand side and noting that $(1-\alpha P(x))$ is minimised by taking $P(x)=1$, we have the first half of the first inequality of Proposition \ref{prop:bound-on-reg}. We achieve the second inequality the same way but by stipulating that $x'\geq-1$.

In the case of the truth-teller to liar interaction, first let
\begin{equation*}\label{eq:gammaeq}
    \gamma(x) = \frac{\nu +\alpha P(x)}{\nu + \alpha^2 P(x)^2}.
\end{equation*}
Then equation (\ref{eq:bininterliar}) can be re-written
\begin{equation}\label{eq:xstargamma}
    x''_* = x'' + \alpha P(x) \gamma(x)(x_d - x'') + \Theta'D(x'').
\end{equation}
We proceed identically to the truth-teller to truth-teller interaction case and achieve the bounds found in the proposition for $\Theta'$. $\blacksquare$

Proposition \ref{prop:bin-inter-conditions} provides conditions for the results of binary interactions to remain within the opinion interval $\mathcal{I} = [-1,1]$. Moving forward, we will only consider binary interactions of this type.

\subsection{Boltzmann equation}\label{sec:boltzmann-eqnanal}

The time evolution of the density $f_T(x,t)$ for a suitable test function $\varphi(x)$ is given by an integro-differential Boltzmann-type equation in weak form. Since we have two distinguishable populations influencing the opinion density, our equation is very similar to the Boltzmann-type equation explored in \cite{albi2014boltzmann} where we have distinguishable populations of leaders and followers, as outlined in Section \ref{sec:leaders-followers}. Indeed,
\begin{equation}\label{eq:Boltzmann}
    \frac{d}{dt}\int_{\mathcal{I}}\varphi(x) f_T(x,t)\, dx = (Q_{T}(f_T,f_T), \varphi) + (Q_{L}(f_T,f_L),\varphi) ,
\end{equation}
where
\begin{subequations}\label{eq:both-boltz-operators}
\begin{equation}\label{eq:bltzoperator}
    (Q_{T}(f_T,f_T),\varphi) = \left<\int_{\mathcal{I}^2}B_{\text{int}}^{T}\,(\varphi(x_*) - \varphi(x))f_T(x,t)f_T(x',t)\, dx\, dx'\right>
\end{equation}
and
\begin{equation}
    (Q_{L}(f_T,f_L),\varphi) = \left<\int_{\mathcal{I}^2}B_{\text{int}}^{L}\,(\varphi(x_*'') - \varphi(x''))f_T(x'',t)f_L(x,t)\, dx''\, dx\right>.
\end{equation}
\end{subequations}
Here, $x_*$ and $x_*''$ denote the results of binary interactions according to equations (\ref{eq:binarystar1}) and (\ref{eq:bininterliar}) respectively, $\left<\cdot\right>$ denotes expectation relative to random variables $\Theta_1, \Theta_2, \Theta'$ and $B_{\text{int}}$ is related to interaction frequency. Guaranteeing that post-interaction opinions preserve their bounds, $B_{\text{int}}^T$ and $B_{\text{int}}^L$ are given by
\begin{eqnarray*}
B_{\text{int}}^{T} &=& B_{\text{int}}^{T}(x,x',\Theta_1,\Theta_2) = \eta_{T}\chi(|x_*|\leq1)\chi(|x'_*|\leq1),\\
B_{\text{int}}^{L} &=& B_{\text{int}}^{L}(x'',x,\Theta') = \eta_{L}\chi(|x_*''|\leq1),
\end{eqnarray*}
where $\eta_{T}, \eta_{L}>0$ are constant rates with $\eta_{L}\ll\eta_{T}$. If the conditions of Proposition \ref{prop:bin-inter-conditions} are satisfied, then we can assume that $|x_*|\leq1$, $|x_*'|\leq1$ and $|x_*''|\leq 1$ so the Boltzmann operators (\ref{eq:both-boltz-operators}) are simplified to
\begin{subequations}\label{eq:bolts-set-rhs}
\begin{equation}\label{eq:boltz-rhs}
    (Q_{T}(f_T,f_T),\varphi) = \eta_{T}\left<\int_{\mathcal{I}^2}(\varphi(x_*) - \varphi(x))f_T(x,t)f_T(x',t)\, dx\, dx'\right>,
\end{equation}
and
\begin{equation}\label{eq:boltz-liar-rhs}
    (Q_{L}(f_T,f_L),\varphi) = \eta_{L}\left<\int_{\mathcal{I}^2}(\varphi(x_*'') - \varphi(x''))f_T(x'',t)f_L(x,t)\, dx''\, dx\right>.
\end{equation}
\end{subequations}

Taking $\varphi(x) = 1$ returns to the fact that the total number of agents is conserved. When $\varphi(x) =x$, we obtain the evolution of the average opinion which can be written as
\begin{multline}\label{eq:mean-derivative-boltz}
    \frac{d}{dt}\int_{\mathcal{I}}x f_T(x,t)\, dx = \frac{\eta_{T}}{2}\left<\int_{\mathcal{I}^2}(x_* + x'_*- x - x')f_T(x,t)f_T(x',t)\, dx\, dx'\right> \\+ \eta_{L}\left<\int_{\mathcal{I}^2} (x''_* - x'')f_{T}(x'',t)f_L(x,t)\, dx''\, dx\right>.
\end{multline}
Let $m_T(t) = \int_{\mathcal{I}}xf_T(x,t)\, dx$ denote the mean opinion of the truth-telling population at time $t\geq0$. We have that $f_L(x,t) = \rho\delta(x - x_d)$ and since $\Theta_1, \Theta_2$ have zero mean, we can substitute equations (\ref{eq:binarystar1}, \ref{eq:bininterliar}) into (\ref{eq:mean-derivative-boltz}) and obtain
\begin{multline}
    \frac{d}{dt}m_T(t) = \frac{\eta_{T}}{2}\int_{\mathcal{I}^2}\left(\alpha P(x)(x'-x) + \alpha P(x')(x-x')\right)f(x,t)f(x',t)\, dx\, dx' \\+ \rho\eta_{L}\int_{\mathcal{I}}\alpha P(x) \left(\frac{\nu + \alpha P(x)}{\nu + \alpha^2P(x)^2}\right)(x_d - x)f_T(x,t)\, dx.
\end{multline} 

In the case where $P(x) = \mathcal{P}$ for $\mathcal{P}>0$ a constant, we can find $m_T(t)$ analytically. The first interaction integral is zero in this case. Simplifying the second term gives
\begin{equation}\label{eq:mean-oneliar}
    \frac{d}{dt}m(t) = \rho\eta_{L} \int_{\mathcal{I}}\alpha \mathcal{P} \left(\frac{\nu + \alpha \mathcal{P}}{\nu + \alpha^2\mathcal{P}^2}\right)(x_d - x)\, f_T(x,t)\, dx.
\end{equation}
Using the definition of $m_T(t)$ and the fact that $f_T$ is a normalised density, we have
\begin{equation*}
    \frac{d}{dt}m_T(t) = \rho\eta_{L}\left(\frac{\nu\alpha \mathcal{P} + \alpha^2\mathcal{P}^2}{\nu + \alpha^2\mathcal{P}^2}\right)(-m_T(t) + x_d).
\end{equation*}
Solving this ODE, we have
\begin{equation*}
    m_T(t) = (m_T(0)-x_d)\exp\left(-\rho\eta_{L}\left(\frac{\nu\alpha \mathcal{P} + \alpha^2\mathcal{P}^2}{\nu + \alpha^2 \mathcal{P}^2}\right)\, t\right) + x_d.
\end{equation*}
Hence, $m_T(t)\to x_d$ as $t\to \infty$ for any $\eta_{L}, \rho, \mathcal{P}>0$, regardless of the initial condition $m_T(0)$. Continuing with $P(x) = \mathcal{P}$, we consider the second moment 
\begin{equation*}
    V_T(t) = \int_{\mathcal{I}}x^2 f(x,t)\, dx.
\end{equation*}
Then
\begin{equation*}
    \frac{d}{dt}V_T(t) = \eta_{T}\left<\int_{\mathcal{I}^2}(x_*^2 - x^2)f_T(x,t)f_T(x',t)\, dx\, dx'\right> + \eta_{L}\left<\int_{\mathcal{I}^2}((x_*'')^2 - x''^2)f_T(x'',t)f_L(x,t)\, dx''\, dx\right>.
\end{equation*}
Substituting equations (\ref{eq:binarystar1}) and (\ref{eq:xstargamma}) for $x_*$ and $x''_*$ respectively, we have
\begin{multline*}
    \frac{d}{dt}V_F(t) = \eta_{T}\int_{\mathcal{I}^2}\left(\alpha^2\mathcal{P}^2(x' - x)^2 + 2\alpha \mathcal{P}x(x'-x)\right)f_T(x,t)f_T(x',t)\, dx\, dx' \\+ \rho\eta_{L}\int_{\mathcal{I}} \left(\alpha^2\mathcal{P}^2\gamma^2(x_d - x)^2 + 2\alpha \mathcal{P}\gamma x(x_d - x)\right)f_T(x,t) + (\eta_{T} + \rho\eta_{L})\sigma^2\int_{\mathcal{I}}D(x)^2f_T(x,t)\, dx,
\end{multline*}
where
$$ \gamma = \frac{\nu + \alpha \mathcal{P}}{\nu + \alpha^2 \mathcal{P}^2},$$
and we use the fact that $\Theta_i$ have zero mean and variance $\sigma^2$. In the absence of diffusion ($D(x)\equiv0$), we simplify the above expression with the following definitions
\begin{equation*}
    \int_{\mathcal{I}^2}x^2f_T(x,t)f_T(x',t)\, dx\, dx' = V_T(t), \quad \int_{\mathcal{I}^2}xx'f_T(x,t)f_T(x',t)\, dx\, dx' = m_T(t)^2, \quad
\end{equation*}
$$ \int_{\mathcal{I}^2}xx_df_T(x,t)f_T(x',t)\, dx\, dx' = x_dm_T(t),$$
to obtain
\begin{multline*}
    \frac{d}{dt}V_T(t) =(2\eta_{T}(\alpha^2\mathcal{P}^2-\alpha \mathcal{P}) + \rho\eta_{L}(\alpha^2\mathcal{P}^2\gamma^2 -2\alpha \mathcal{P}\gamma))V_T(t) + 2\eta_{T}(-\alpha^2\mathcal{P}^2 + \alpha \mathcal{P})m_T(t)^2 \\+ \rho\eta_{L}\alpha^2\mathcal{P}^2\gamma^2x_d^2 + 2\rho\eta_{L}(-\alpha^2\mathcal{P}^2\gamma^2 + \alpha \mathcal{P}\gamma)x_dm_T(t). 
\end{multline*}
Since $m_T(t)\to x_d$ as $t\to\infty$, at large times we have
\begin{equation*}
    \frac{d}{dt}V_T(t) = (-2\eta_{T}\alpha \mathcal{P}(1-\alpha \mathcal{P}) - \rho\eta_{L}\alpha \mathcal{P} \gamma(2 - \alpha \mathcal{P}\gamma))(V_T(t) - x_d^2).
\end{equation*}
A sufficient condition for the prefactor $(-2\eta_{T}\alpha \mathcal{P}(1-\alpha \mathcal{P}) - \rho\eta_{L}\alpha \mathcal{P} \gamma(2 - \alpha \mathcal{P}\gamma))$ to be non-positive is that $\alpha \mathcal{P}\leq1$ and $\alpha \mathcal{P}\gamma\leq 2$. 

Assuming this condition, we have that $V_T(t)\to x_d^2$ as $t\to\infty$. Therefore, the quantity
\begin{equation*}
    \int_{\mathcal{I}}(x-x_d)^2f_T(x,t)\, dx = V_T(t) - 2x_dm_T(t) + x_d^2 
\end{equation*}
tends to zero for large $t$. Hence, the long time density $f_{T, \infty}(x)$ tends to a Dirac delta density around the liar's true opinion
$$ f_{T, \infty}(x) = \delta(x-x_d).$$

\subsection{Quasi-invariant limit}\label{sec:quasi-invariant-limit}

In order to maintain memory of the microscopic interactions at an asymptotic level, we perform a rescaling of the variables $\eta_T,\eta_L,\alpha,\sigma^2$ and $\nu$. This rescaling is important as we want to consider the impact of small interactions in order to apply Taylor's theorem to $\varphi(x_*) - \varphi(x)$ in the Boltzmann equation. We make the following assumptions:
\begin{equation*}
    \alpha = \varepsilon,\quad \eta_T = \frac{1}{c_T\varepsilon}, \quad \eta_L = \frac{1}{c_L\varepsilon},\quad \sigma^2 = \varepsilon\zeta, \quad \nu = \varepsilon\kappa,
\end{equation*}
where $\varepsilon, \zeta,\kappa, c_T, c_L>0$ and $\varepsilon\ll1$. We require that in the limit $\varepsilon\to0$, the main macroscopic properties of the kinetic system are preserved. To this aim, consider the evolution of the scaled first two moments when $P(x) = \mathcal{P}$,
\begin{equation*}
    \frac{d}{dt}m_T(t) = \frac{\rho}{c_L}\frac{\kappa P + P^2}{\kappa + \varepsilon P}(-m_T(t) + x_d),
\end{equation*}
\begin{multline*}
    \frac{d}{dt}V_T(t) = \left(-\frac{2}{c_T}\mathcal{P}(1-\varepsilon \mathcal{P}) - \frac{\rho}{c_L}\frac{\mathcal{P}(\kappa + \mathcal{P})}{\kappa + \varepsilon \mathcal{P}^2} \left(2 - \varepsilon P\frac{\kappa + \mathcal{P}}{\kappa + \varepsilon \mathcal{P}^2}\right)\right)V_T(t) + \frac{2}{c_T}\mathcal{P}(1-\varepsilon \mathcal{P})m_T(t)^2\\ + \frac{2\rho}{c_L}\frac{\mathcal{P}(\kappa + \mathcal{P})}{\kappa + \varepsilon \mathcal{P}^2}\left(1-\varepsilon \mathcal{P} \frac{\kappa + \mathcal{P}}{\kappa + \varepsilon \mathcal{P}^2}\right)x_dm_T(t) + \frac{\rho}{c_L} \varepsilon \mathcal{P}^2\frac{(\kappa + \mathcal{P})^2}{(\kappa + \varepsilon \mathcal{P})^2}x_d^2 + \left(\frac{1}{c_T} + \frac{\rho}{c_L}\right)\zeta\int_{\mathcal{I}} D(x)^2f_T(x,t)\, dx.
\end{multline*}
In the limit $\varepsilon\to 0$, the above equations give
\begin{equation}\label{eq:boltzmann-mean-quasi-invariant}
    \frac{d}{dt}m_T(t) = \frac{\rho}{c_L}\frac{\kappa \mathcal{P} + \mathcal{P}^2}{\kappa}(-m_T(t) + x_d),
\end{equation}
and
\begin{multline*}
    \frac{d}{dt}V_T(t) = \left(-\frac{2}{c_T}\mathcal{P} - \frac{\rho}{c_L}\frac{2\mathcal{P}(\kappa + \mathcal{P})}{\kappa}\right)V_T(t) + \frac{2}{c_T}\mathcal{P}m_T(t)^2 + \frac{\rho}{c_L}\frac{2\mathcal{P}(\kappa + \mathcal{P})}{\kappa}x_dm_T(t) \\+ \zeta\left(\frac{1}{c_T} + \frac{\rho}{c_L}\right)\int_{\mathcal{I}} D(x)^2f_T(x,t)\, dx.
\end{multline*}
In this limit, the decay of $V_T(t)$ is far more obvious and we get exactly the long time behaviour that we would expect in the diffusion-free regime, $D(x)\equiv0$. 

The Boltzmann equation (\ref{eq:Boltzmann}) with our new scaling gives
\begin{multline*}
    \frac{d}{dt}\int_{\mathcal{I}}\varphi(x)f(x,t) \, dx = \frac{1}{c_T\varepsilon}\left<\int_{\mathcal{I}^2}(\varphi(x_*) - \varphi(x))f_T(x,t)f_T(x',t)\, dx\, dx'\right> \\+ \frac{1}{c_L\varepsilon}\left<\int_{\mathcal{I}^2}(\varphi(x_*'') - \varphi(x''))f_T(x'',t)f_L(x',t)\, dx''\, dx'\right>,
\end{multline*}
and the binary interaction dynamics (\ref{eq:binarystar1}) and (\ref{eq:bininterliar}) can be written as
\begin{equation}\label{eq:rescalebin}
    x_* - x = \varepsilon P(x) (x' - x) + \Theta_1^{\varepsilon}D(x),
\end{equation}
where $\Theta_1^{\varepsilon}$ is a random variable with zero mean and variance $\varepsilon\zeta$, and
\begin{equation}\label{eq:rescaleliar}
    x''_* - x'' = \varepsilon P(x'') \frac{\kappa  + P(x'')}{\kappa + \varepsilon P(x'')^2}(x_d -x'') + \Theta'^{\varepsilon}D(x'')
\end{equation}
where $\Theta'^{\varepsilon}$ is also a normal random variable with zero mean and variance $\varepsilon\zeta$. This means that in the parameter range $\varepsilon\ll1$, the change in opinion due to binary interactions is very small.

\subsection{Fokker-Planck equation}\label{sec:fokker-planck-easyP}

In order to recover the limit as $\varepsilon\to 0$, we consider the second order Taylor expansion of $\varphi$ around $x$,
\begin{eqnarray*}
    \varphi(x_*) - \varphi(x) &=& (x_* - x)\varphi'(x) + \frac{1}{2}(x_* - x)^2\varphi''(\tilde{x}),\\
    \varphi(x''_*) - \varphi(x'') &=& (x''_* - x'')\varphi'(x'') + \frac{1}{2}(x_*'' - x'')^2\varphi''(\hat{x}),
\end{eqnarray*}
where 
$$ \tilde{x} = \theta x_* + (1-\theta)x, \quad \hat{x} = \hat{\theta} x''_* + (1-\hat{\theta})x'',$$
for some $\theta, \hat{\theta}\in [0,1]$. Inserting these expansions into the interaction integrals (\ref{eq:bolts-set-rhs}), we obtain
\begin{equation*}
    (Q_T(f_T,f_T),\varphi)=\frac{1}{c_T\varepsilon}\left<\int_{\mathcal{I}^2}\left((x_*-x)\varphi'(x) + \frac{1}{2}(x_* - x)^2\varphi''(x)\right) f_T(x,t)f_T(x',t)\, dx\, dx'\right> + R_1(\varepsilon),
\end{equation*}
and
\begin{equation*}
    (Q_L(f_T,f_L),\varphi)=\frac{1}{c_L\varepsilon}\left<\int_{\mathcal{I}^2}\left((x''_*-x'')\varphi'(x'') + \frac{1}{2}(x''_* - x'')^2\varphi''(x'')\right) f_T(x'',t)f_L(x',t)\, dx''\, dx'\right> + R_2(\varepsilon),
\end{equation*}
respectively, where 
$$ R_1(\varepsilon) = \frac{1}{2c_T\varepsilon}\left<\int_{\mathcal{I}^2}(x_*-x)^2(\varphi''(\tilde{x})-\varphi''(x))f_T(x,t)f_T(x',t)\, dx\, dx'\right>,$$
and
$$ R_2(\varepsilon) = \frac{1}{2c_L\varepsilon}\left<\int_{\mathcal{I}^2}(x''_*-x'')^2(\varphi''(\hat{x})-\varphi''(x''))f_T(x'',t)f_L(x',t)\, dx''\, dx'\right>.$$
Then, by equations (\ref{eq:rescalebin}) and (\ref{eq:rescaleliar}), we can write the interaction integrals as
\begin{equation}\label{eq:interactint}
    (Q_T(f_T,f_T),\varphi) = \frac{1}{c_T}\int_{\mathcal{I}^2}\left[P(x)(x'-x)\varphi'(x) + \frac{\zeta}{2}D(x)^2\varphi''(x)\right]f_T(x,t)f_T(x',t)\, dx\, dx' + R_1(\varepsilon)+ \mathcal{O}(\varepsilon),
\end{equation}
and
\begin{multline}\label{eq:interactint2}
    (Q_L(f_T,f_L),\varphi) = \frac{\rho}{c_L}\int_{\mathcal{I}}\left[P(x'')\frac{\kappa +P(x'')}{\kappa + \varepsilon P(x'')^2}(x_d - x'')\varphi'(x'')+ \frac{\zeta}{2}D(x'')^2\varphi''(x'')\right]f_T(x'',t)\, dx'\\ + R_2(\varepsilon) + \mathcal{O}(\varepsilon),
\end{multline}
where we have used the fact that $\Theta_1^{\varepsilon}, \Theta'^{\varepsilon}$ have mean 0 and variance $\varepsilon\zeta$.
\proposition{Assuming the H\"{o}lder continuity of $\varphi(x)$, i.e., $\varphi\in \mathcal{F}_{2+\delta}(\mathcal{I})$ given by equation (\ref{eq:holder-continuity}), with order $0<\delta\leq 1$, $\kappa>0, x_d\in\mathcal{I}$ and $0\leq P(x)\leq1$, the remainder terms $R_1(\varepsilon)\to 0$ and $R_2(\varepsilon)\to 0$ as $\varepsilon\to 0$.}\label{prop:R2tozero}\normalfont

\textbf{Proof}: By the same argument as Toscani \cite{toscani2006kinetic}, the first remainder term $R_1(\varepsilon)\to 0$. Now for $R_2(\varepsilon)$. By the triangle inequality, we have
\begin{equation*}
    |R_2(\varepsilon)|\leq \frac{1}{2c_L\varepsilon}\int_{\mathcal{I}^2}\left<\left|(x''_* - x'')^2(\varphi''(\hat{x}) - \varphi''(x''))\right|\right>f_T(x'',t)f_L(x',t)\, dx''\, dx'.
\end{equation*}
Since $\varphi\in\mathcal{F}_{2+\delta}(\mathcal{I})$ and $|\hat{x} - x''| = \hat{\theta}|x''_* - x''|$, we have
\begin{equation*}
    |\varphi''(\hat{x}) - \varphi''(x'')|\leq ||\varphi''||_{\delta}|\hat{x} - x''|^{\delta}\leq\hat{\theta}||\varphi''||_{\delta}|x''_* - x''|^{\delta},
\end{equation*}
where 
\begin{equation}\label{eq:holder-continuity}
    ||\varphi''||_{\delta} = \sup_{x',x\in\mathcal{I}, x'\neq x}\frac{|\varphi''(x') - \varphi''(x)|}{|x' - x|^{\delta}}<\infty,
\end{equation}
for some $0<\delta \leq 1$. The H\"{o}lder inequality, gives us
\begin{equation*}
    |R_2(\varepsilon)|\leq \frac{1}{2c_L\varepsilon}\int_{\mathcal{I}^2}\left<\left|(x''_* - x'')^{4}\right|\right>^{1/2}\left<\left|(\varphi''(\hat{x}) - \varphi''(x''))^{2}\right|\right>^{1/2}f_T(x'',t)f_L(x',t)\, dx''\, dx',
\end{equation*}
and then by H\"{o}lder continuity,
\begin{equation*}
    |R_2(\varepsilon)|\leq \frac{\hat{\theta}||\varphi''||_{\delta}}{2c_L\varepsilon}\int_{\mathcal{I}^2}\left<\left|x''_* - x''\right|^{4 + 2\delta}\right>^{1/2} f_T(x'',t)f_L(x',t)\, dx''\, dx'.
\end{equation*}
From equation (\ref{eq:rescaleliar}), we have
\begin{align*}
    |x_*'' - x''|^{4 + 2\delta} &= \left| \varepsilon P(x'') \frac{\kappa + P(x'')}{\kappa + \varepsilon P(x'')^2}(x_d - x'') + \Theta'^{\varepsilon}D(x'')\right|^{4+2\delta}\\ &\leq \left| \varepsilon P(x'') \frac{\kappa + P(x'')}{\kappa + \varepsilon P(x'')^2}(x_d - x'')\right|^{4+2\delta} + |\Theta'^{\varepsilon}D(x'')|^{4+2\delta}
\end{align*}
by the triangle inequality. The first term on the right hand side is maximised by taking $P(x'') = 1$ and $|x_d - x''| = 2$, hence
$$ |x_*'' - x''|^{4 + 2\delta}\leq \left|2\varepsilon\frac{\kappa + 1}{\kappa + \varepsilon}\right|^{4+2\varepsilon} + |\Theta'^{\varepsilon}D(x'')|^{4+2\delta}.$$
Substituting this inequality into the integral, we have
\begin{equation*}
    |R_2(\varepsilon)| \leq \frac{\rho\hat{\theta}||\varphi''||_{\delta}}{2c_L\varepsilon}\int_{\mathcal{I}}\left<\left(2\varepsilon\frac{\kappa + 1}{\kappa + \varepsilon}\right)^{4+2\delta}+(|\Theta'^{\varepsilon}D(x'')|^2)^{2+\delta}\right>^{1/2}f_T(x'',t)\, dx'',
\end{equation*}
and then by linearity of expectation, we have
$$ |R_2(\varepsilon)|\leq \frac{\rho\hat{\theta}||\varphi''||_{\delta}}{2c_L\varepsilon}\int_{\mathcal{I}}\varepsilon^{1+\delta/2}\left(\varepsilon^{2+\delta}\left(2\frac{\kappa + 1}{\kappa + \varepsilon}\right)^{4 + 2\delta} + (\zeta D(x'')^2)^{2+\delta}\right)^{1/2}f_T(x'',t)\, dx''.$$
Hence, as $\varepsilon\to 0$, we have that $|R_2(\varepsilon)|\to 0$. $\blacksquare$

Taking the limit of the interaction integrals (\ref{eq:interactint}) and (\ref{eq:interactint2}) as $\varepsilon\to 0$ and applying Proposition \ref{prop:R2tozero}, we have
\begin{multline*}
    \frac{d}{dt}\int_{\mathcal{I}}\varphi(x)f(x,t)\, dx = \frac{1}{c_T}\int_{\mathcal{I}^2}\left[P(x)(x'-x)\varphi'(x)+ \frac{\zeta}{2}D(x)^2\varphi''(x)\right]f_T(x,t)f_T(x',t)\, dx\,dx'\\ + \frac{\rho}{c_L}\int_{\mathcal{I}}\left[P(x) \left(1 + \frac{P(x)}{\kappa}\right)(x_d-x)\varphi'(x) + \frac{\zeta}{2}D(x)^2\varphi''(x)\right]f_T(x,t)\, dx.
\end{multline*}
Integrating by parts twice, we have the integral equation
\begin{multline*}
    \int_{\mathcal{I}}\frac{\partial}{\partial t}(f_T(x,t))\varphi(x)\, dx = -\frac{1}{c_T}\int_{\mathcal{I}^2}\frac{\partial}{\partial x}\left(\left(P(x)(x' -x)\right)f_T(x,t)f_T(x',t)\right)\varphi(x)\, dx\, dx' \\- \frac{\rho}{c_L}\int_{\mathcal{I}}\frac{\partial}{\partial x}\left(P(x)\left(1 + \frac{P(x)}{\kappa}\right)(x_d - x)f_T(x,t)\right)\varphi(x)\, dx + \left(\frac{1}{c_T} + \frac{\rho}{c_L}\right)\int_{\mathcal{I}}\frac{\zeta}{2}\frac{\partial^2}{\partial x^2}(D(x)^2f_T(x,t))\varphi(x)\, dx.
\end{multline*}
Then, since $\varphi(x)$ is an arbitrary test function,we apply the variational lemma and arrive at the Fokker-Planck equation:
\begin{equation}\label{eq:fokplanck}
    \frac{\partial}{\partial t}f_T(x,t) + \frac{1}{c_T}\frac{\partial }{\partial x}\left(\mathcal{H}[f_T](x)f_T(x,t)\right) +\frac{\rho}{c_L}\frac{\partial }{\partial x}\left(\mathcal{K}[f_T](x)f_T(x,t)\right) = \left(\frac{1}{c_T} +\frac{\rho}{c_L}\right)\frac{\zeta}{2}\frac{\partial^2}{\partial x^2}(D(x)^2f_T(x,t)),
\end{equation}
where
\begin{eqnarray*}
    \mathcal{H}[f_T](x) &=& \int_{\mathcal{I}}P(x)(x'- x)f_T(x',t)\, dx' = P(x)(m_T(t) -x),\\
    \mathcal{K}[f_T](x) &=& P(x)\left(1+ \frac{P(x)}{\kappa}\right)(x_d-x).
\end{eqnarray*}
Equation (\ref{eq:fokplanck}) is complemented with no-flux boundary conditions and an initial condition $f_T(x,0) = f_0(x)$. The first drift term, $\mathcal{H}[f_T]$, represents drift toward the mean opinion of the truth tellers. The second term, $\mathcal{K}[f_T]$ represents drift towards the liars' opinion, $x_d$. For very small regularisation parameter $\kappa$, we see that this drift term will be very large -- this is consistent with the unregularised state where the liar can lie in any way they like and therefore bring consensus toward their goal opinion very quickly. In contrast, when $\kappa$ is very large, the liar can only show their true opinion. In this case, there is still a drift toward the liar's opinion (since all agents are connected and the liar's opinion does not change) but the liar's influence is comparable to any ordinary agent. Note that the term $\mathcal{K}[f_T]$ carries second order dependence on $P(x)$. Intuitively, we may explain this dependence by the fact that the lie told by agent 1 depends on the interaction function $P(x)$ and furthermore, in communicating, the lie must be passed through the interaction function. Hence the contribution from the liar is second order in $P(x)$. The prefactor $\rho/c_L$ determines the strength in communication of the liar's opinion compared to the overall strength of communication between truth-telling agents. We wish to take $\rho$ small to represent that the liar is only one person in a very large population size. When $\rho=0$, we recover the standard Fokker-Planck equation for an interacting system of agents obtained by Toscani in \cite{toscani2006kinetic}.

\subsection{Stationary solutions}

We next consider steady-state solutions to the Fokker-Planck equation (\ref{eq:fokplanck}) in the case $P(x) = \mathcal{P}$ for $0\leq \mathcal{P}\leq 1$ a constant. For $\mathcal{P}>0$ we have shown that $m_T(t)\to x_d$ exponentially by equation (\ref{eq:boltzmann-mean-quasi-invariant}) so at the steady state we have,
\begin{equation*}
    \frac{1}{c_T}\frac{\partial}{\partial x}\left(\mathcal{P}(x_d-x)f_T(x)\right) + \frac{\rho}{c_L}\frac{\partial }{\partial x}\left(\mathcal{P}\left(1 + \frac{\mathcal{P}}{\kappa}\right)(x_d-x)f_T(x)\right) = \left(\frac{1}{c_T} + \frac{\rho}{c_L}\right)\frac{\zeta}{2}\frac{\partial^2}{\partial x^2}\left(D(x)^2f_T(x)\right).
\end{equation*}
Integrating with respect to $x$ and applying the no-flux boundary condition, we have
\begin{equation}\label{eq:FP-constantP}
    \left(\frac{1}{c_T}+ \frac{\rho}{c_L}\left(1+ \frac{\mathcal{P}}{\kappa}\right)\right)\mathcal{P}(x_d-x)f_T(x) = \left(\frac{1}{c_T}+ \frac{\rho}{c_L}\right)\frac{\zeta}{2}\frac{\partial}{\partial x}\left(D(x)^2f_T(x)\right).
\end{equation}
This equation has the solution,
\begin{equation*}
    f_{\infty}^{A}(x) = \frac{C}{D(x)^2}\exp\left(2A\int^x \frac{x_d - s}{D(s)^2}\, ds\right),
\end{equation*}
where $C$ is such that $\int_\mathcal{I}f_T(x) = 1$ and
\begin{equation}\label{eq:constant-stationary-solution}
    A = \frac{\mathcal{P}}{\zeta} + \frac{\mathcal{P}^2}{\kappa\zeta}\left(\frac{\rho}{c_L/c_T + \rho}\right)>0.
\end{equation}
In order for Proposition \ref{prop:bin-inter-conditions} to be satisfied, it is necessary that $D(\pm1) = 0$. Hence, $D(x) = 1-x^2$ is a simple choice of diffusion function satisfying conditions of Proposition \ref{prop:bin-inter-conditions}. When $D(x) = 1-x^2$, (\ref{eq:FP-constantP}) has the exact solution
\begin{equation*}
    f_{\infty}^{A}(x) = \frac{C}{(1-x^2)^2}\left(\frac{1+x}{1-x}\right)^{A\frac{x_d}{2}}\exp\left(-A\left(\frac{1-x_dx}{1-x^2}\right)\right).
\end{equation*}
In the expression for $A$, (\ref{eq:constant-stationary-solution}), we see the ratio of rates, $c_L/c_T$ appearing. If $1/c_L\gg1/c_T$ (i.e., the liar interacts with the populace much more than the populace interacts with one another), then the second term becomes very large, increasing the slope of the peak of the distribution, centred about $x_d$. In contrast, if $1/c_L\ll1/c_T$, then the size of $\rho$ becomes very important in determining how well the liar influences the populace.

If $\kappa\to\infty$ then $A\to \frac{\mathcal{P}}{\zeta}$, hence the regularisation term dominates and we expect to recover an uncontrolled system. Let $f_{\infty}^{\mathcal{P},\zeta}$ be the corresponding density, then
\begin{equation}\label{eq:stead-state-P-const}
    f_{\infty}^{\mathcal{P},\zeta}(x) = \frac{C}{(1-x^2)^2}\left(\frac{1+x}{1-x}\right)^{\frac{\mathcal{P}x_d}{2\zeta}}\exp\left(-\frac{\mathcal{P}(1-x_dx)}{\zeta(1-x^2)}\right).
\end{equation}
Note that this solution is independent of the rates $c_L,c_T$ and the relative weighting of the liar, $\rho$. When $\kappa\to\infty$, the liar communicates in the same way as a truth-telling individual except that the liar does not change their opinion in time. Hence the rates with which truth-tellers interact with each other versus interact with the liar have no bearing on the steady-state distribution.

\begin{figure}[htb!]
    \centering
    \begin{subfigure}[b]{0.45\textwidth}
        \centering
        \includegraphics[width=\textwidth]{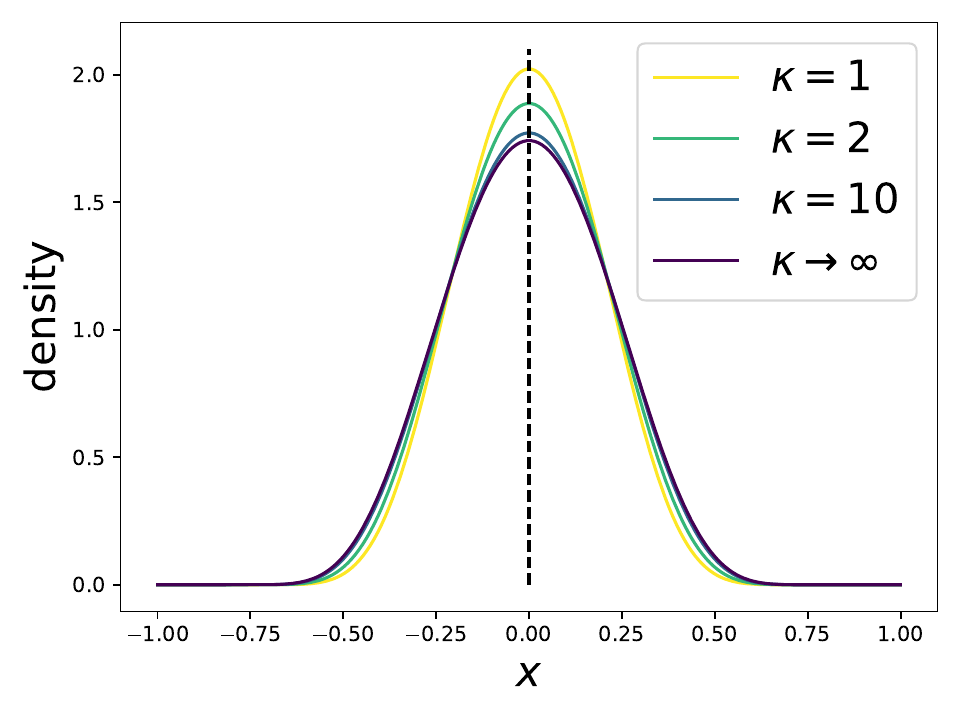}
        \caption{$x_d = 0$.}
    \end{subfigure}
    \begin{subfigure}[b]{0.45\textwidth}
        \centering
        \includegraphics[width=\textwidth]{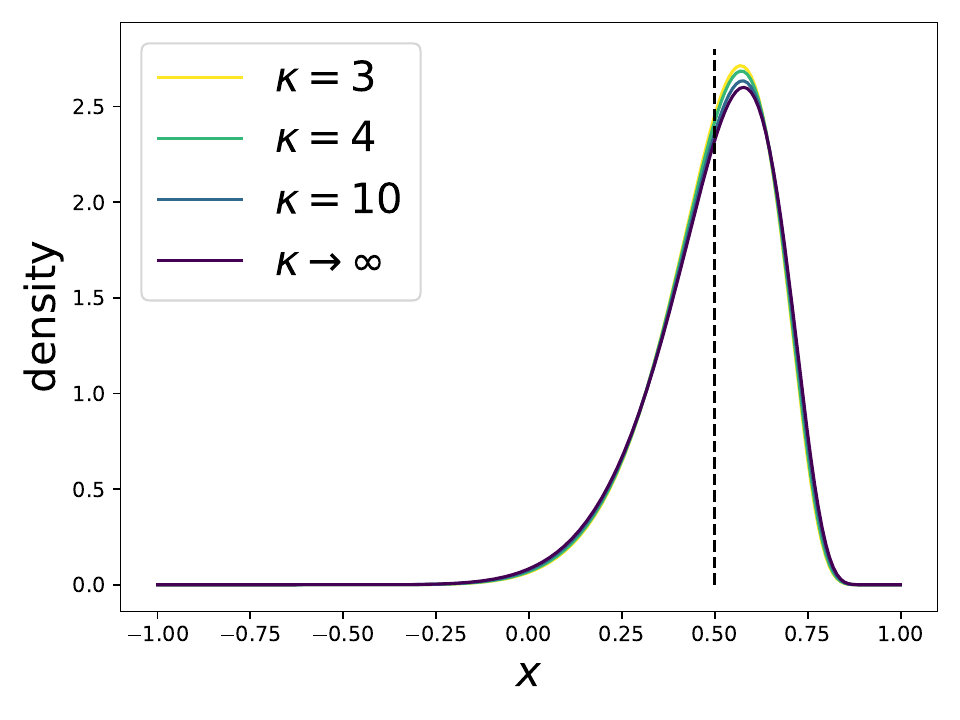}
        \caption{$x_d = 0.5$}
    \end{subfigure}
    \caption{Steady state profiles of distributions $f_T(x,t)$ in the case where $P(x) \equiv \mathcal{P} \equiv 1$ and $D(x)=1-x^2$. We take $\zeta = 0.1, c_T = c_L = 1$ and $\rho = 0.5$. Figure \textbf{(a)} shows profiles when the goal of the liar is $x_d = 0$ and \textbf{(b)} when it corresponds to $x_d = 0.5$, represented by a dashed line in each case. The values of $\kappa$ used ensure that Proposition \ref{prop:bound-on-reg} is satisfied.}
    \label{fig:steady-state-constP}
\end{figure}

Figure \ref{fig:steady-state-constP} displays some of these simple steady state solutions. We choose our values of $\kappa$ in accordance with Proposition \ref{prop:bound-on-reg}. The profile of the steady-state solution when $x_d=0$ has a peak at $x = 0$ and is Gaussian in quality. The peak of the steady-state solution when $x_d=0.5$ lies slightly above $x_d$ itself. This is due to the fact that the diffusion term $D(x)$ smoothly approaches zero at the boundary.

\begin{figure}[htb!]
    \centering
    \begin{subfigure}[b]{0.45\textwidth}
        \centering
        \includegraphics[width=\textwidth]{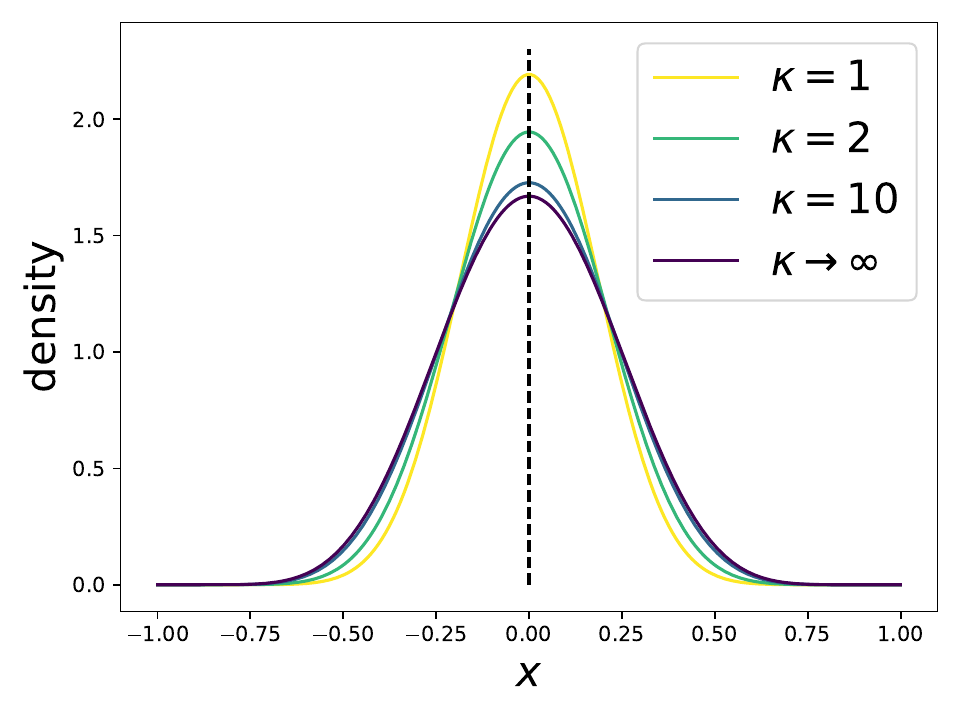}
        \caption{$x_d = 0$.}
    \end{subfigure}
    \begin{subfigure}[b]{0.45\textwidth}
        \centering
        \includegraphics[width=\textwidth]{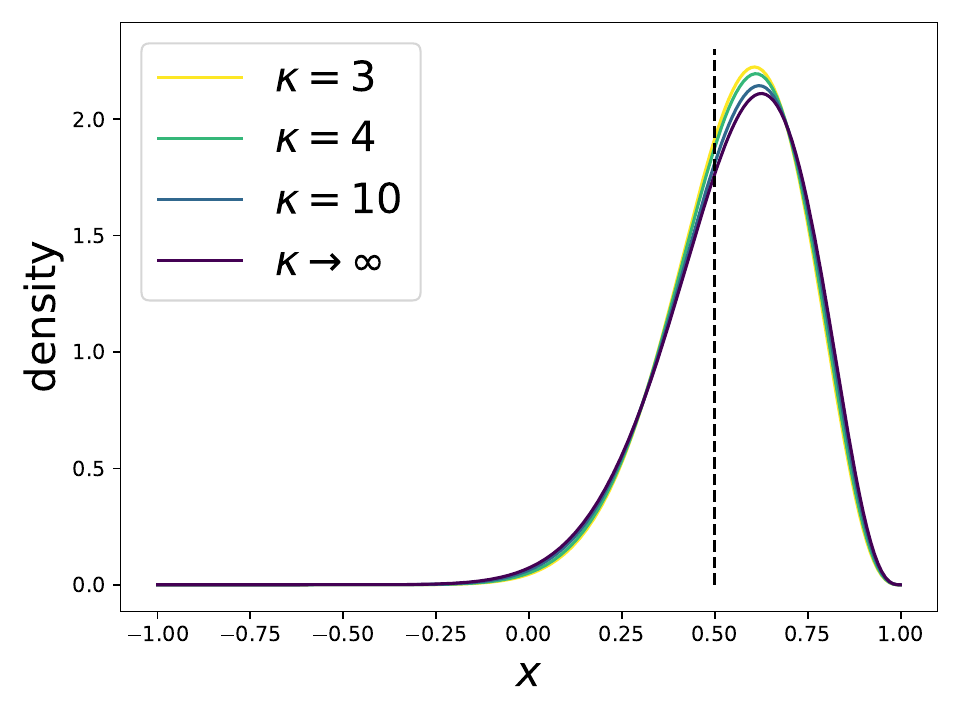}
        \caption{$x_d = 0.5$}
    \end{subfigure}
    \caption{Steady state profiles of distributions $f_T(x,t)$ in the case where $P(x) =1-x^2$ and $D(x)=1-x^2$. We take $\zeta = 0.1, c_T = c_L = 1$ and $\rho = 0.5$. Figure \textbf{(a)} shows profiles when the goal opinion of the liar is $x_d = 0$ and \textbf{(b)} when it is $x_d = 0.5$, represented by a dashed line in each case.}
    \label{fig:steady-state-simpleP}
\end{figure}

Now consider general $P(x)$ where $P(x)\not\equiv 0$ so $m_T(t) \to x_d$. Then stationary solutions satisfy
\begin{equation*}
    \left(\frac{1}{c_T} + \frac{\rho}{c_L}\left(1+ \frac{P(x)}{\kappa}\right)\right)P(x)(x_d - x)f_T(x) = \left(\frac{1}{c_T} + \frac{\rho}{c_L}\right)\frac{\zeta}{2}\frac{\partial}{\partial x}\left(D(x)^2f_T(x)\right).
\end{equation*}
We are able to find an exact solution for $f(x)$ when $P(x) = 1-x^2$ and $D(x) = 1-x^2$,
\begin{equation*}
    f_{\infty}^{\kappa,\zeta}(x) = \frac{C}{(1-x^2)^{2-1/\zeta}}\left(\frac{1+x}{1-x}\right)^{\frac{x_d}{\zeta}}\exp\left(-\frac{1}{Bc_L\kappa\zeta}\left(x^2 - 2x_d x\right)\right),
\end{equation*}
where
$$ B = \left(\frac{1}{c_T}+ \frac{\rho}{c_L}\right).$$
As $\kappa\to\infty$, the exponential term goes to one so the limiting distribution is simply,
\begin{equation*}
    f_{\infty}^{\zeta}(x) = \frac{C}{(1-x^2)^{2-1/\zeta}}\left(\frac{1+x}{1-x}\right)^{\frac{x_d}{\zeta}}.
\end{equation*}
Again, this expression carries no dependence on $c_T, c_L$ and $\rho$. Figure \ref{fig:steady-state-simpleP} displays some examples of these steady states in the cases $x_d=0$ and $x_d = 0.5$. We observe that for low values of the rescaled regularisation parameter $\kappa$, the peak of the distribution is slightly higher than in the $P(x) = 1$ case (see Figure \ref{fig:steady-state-constP}). However, for larger values of $\kappa$ the peak is lower. This reflects a trade-off between the liar's influence and the diffusion in the model. For small $\kappa$, the liar has a more dramatic and impactful effect on the truth-tellers whereas when $\kappa$ is large the liar cannot lie at all so causes less rate of change to the opinions of the truth-tellers and diffusion effects become more important.

\vspace{-0.25cm}
\section{Bounded confidence control}\label{sec:bounded-conf}

Throughout Section \ref{sec:boltzmann}, we assume that $P(x_i, x_j) = P(x_i)$. A more realistic choice of interaction function is the so-called \textit{bounded confidence} model \cite{krause2000discrete}, where $P(x_i, x_j) = P(x_j - x_i)$ is a symmetric function of the distance between opinions $x_i$ and $x_j$ given by equation (\ref{eq:bounded-con-introduciton}). The bounded confidence interaction kernel is also commonly written as
\begin{equation}\label{eq:bounded-conf-interfun}
    P(x_j - x_i) = \chi\left(|x_j - x_i|\leq R\right),
\end{equation}
where $R\in(0,|\mathcal{I}|]$ and $\chi$ is the indicator function.

\subsection{Microscopic model}\label{sec:bounded-conf-micro}

We will consider the interaction function (\ref{eq:bounded-conf-interfun}) and the resulting control in the microscopic model (\ref{eq:xidyn}) when the liar acts according to the classic regularisation control (\ref{eq:SR-control}) described in Section \ref{sec:Micro-model-std-reg}. 

Central to the argument is the fact that the liar either communicates with an agent $i$ by presenting an opinion that is inside the confidence interval $[x_i-R, x_i+R]$, where the interaction function $P(x_j - x_i) = 1$, or it is too expensive to communicate with agent $i$, in which case there is no benefit to lying to agent $i$ at all. Following this reasoning, our first step is to calculate the optimal lie in the case where $P(x_j - x_i) \equiv1$. By equation (\ref{eq:SR-control}), we have that if the liar can communicate with an agent of opinion $x_i^n$, $y_i^n$ is given by a projection of $\tilde{y}_i^n$, where
\begin{equation}\label{eq:bounded-con-control}
    \tilde{y}_i^n  = x_d - \frac{\Delta t}{\nu N}(x_i^{n+1} - x_d),
\end{equation}
onto $[x_i^n-R, x_i^n + R]\cup \{x_d\}$. Here, if $y_i^n\in[x^n_i-R, x^n_i+R]$ then
\begin{equation}\label{eq:xn1-withlie}
   x_i^{n+1} = x_i^n + \frac{\Delta t}{N}\left((y_i^n - x_i^n)+ \sum_{j=2}^NP(x_j^n - x_i^n)(x_j^n - x_i^n)\right), 
\end{equation}
and otherwise
\begin{equation}\label{eq:xn1-nolie}
    x_i^{n+1} = x_i^n + \frac{\Delta t}{N}\left( \sum_{j=2}^NP(x_j^n - x_i^n)(x_j^n - x_i^n)\right).
\end{equation}
Next we need to define the projection. Intuitively, we say that the liar will lie if there is a greater gain from influencing agent $i$ than there is cost from telling the lie.

We will discretise the cost function (\ref{eq:std-reg-cost}) by the trapezoidal rule, $\int_a^bf(t) \approx \frac{1}{2}(b-a)(f(a) +f(b))$, giving
\begin{equation*}
    \frac{\Delta t}{2}\left(\frac{1}{N}\sum_{i=1}^N\frac{1}{2}\left((x_i^n - x_d)^2 + (x_i^{n+1} - x_d)^2\right) + \frac{2}{N}\sum_{i=2}^N\frac{\nu}{2}(y_i^n - x_d)^2\right).
\end{equation*}
Since each $x_i^{n+1}$ is independent of $x_j^{n+1}$ for $j\neq i$, we can consider elements of the sum as independent. Hence, comparing the contribution to the cost function with/without lying, we have the following condition: $y_i^n = \mathbb{P}_{[x_i^n-R, x_i^n+R]}(\tilde{y}_i^n)$ if
\begin{equation}\label{eq:projection-definition}
    (x_i^{n+1} - x_d)^2 + 2\nu\left(\mathbb{P}_{[x_i^n-R, x_i^n+R]}(\tilde{y}_i^n) - x_d\right)^2< (\tilde{x}_i^{n+1} - x_d)^2,
\end{equation}
where $x_i^{n+1}$ is given by equation (\ref{eq:xn1-withlie}) and $\tilde{x}_i^{n+1}$ is given by equation (\ref{eq:xn1-nolie}). Otherwise, $y_i^n = x_d$.

\begin{figure}[htb!]
    \centering
    \begin{subfigure}[b]{0.9\textwidth}
        \centering
        \includegraphics[width=\textwidth]{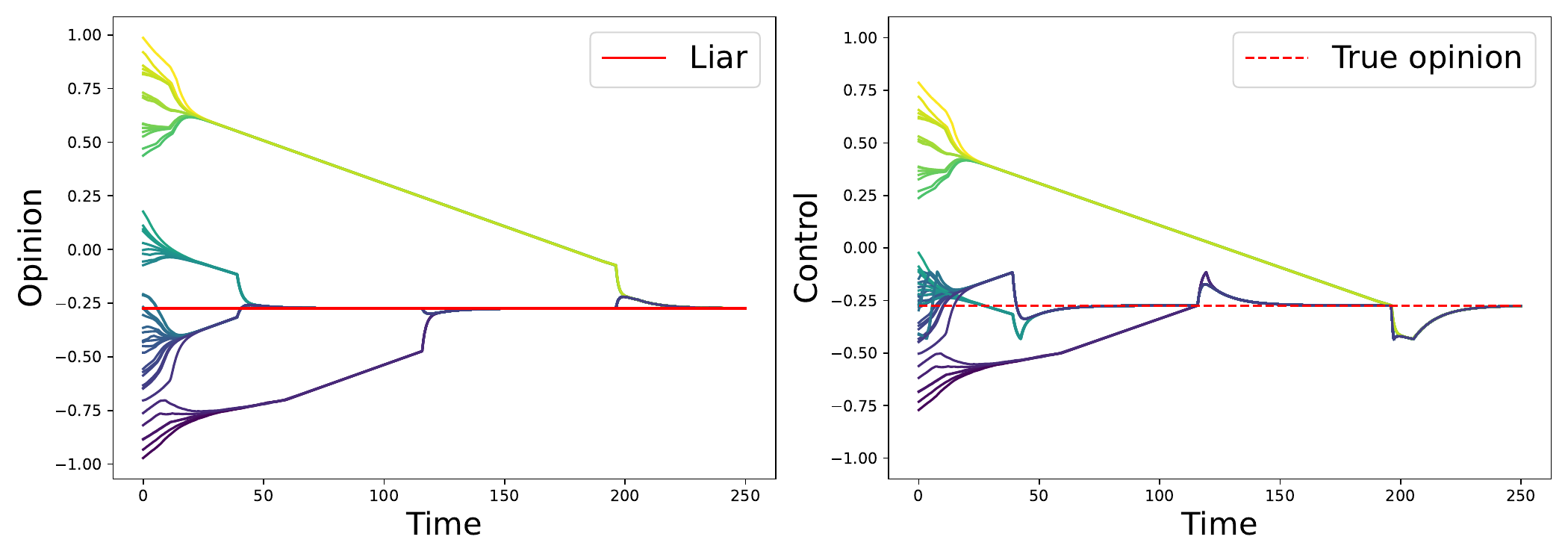}
        \caption{$\nu = 0.0001$}
        \label{fig:bounded-conf-nulessnutilde}
    \end{subfigure}
    \begin{subfigure}[b]{0.9\textwidth}
        \centering
        \includegraphics[width=\textwidth]{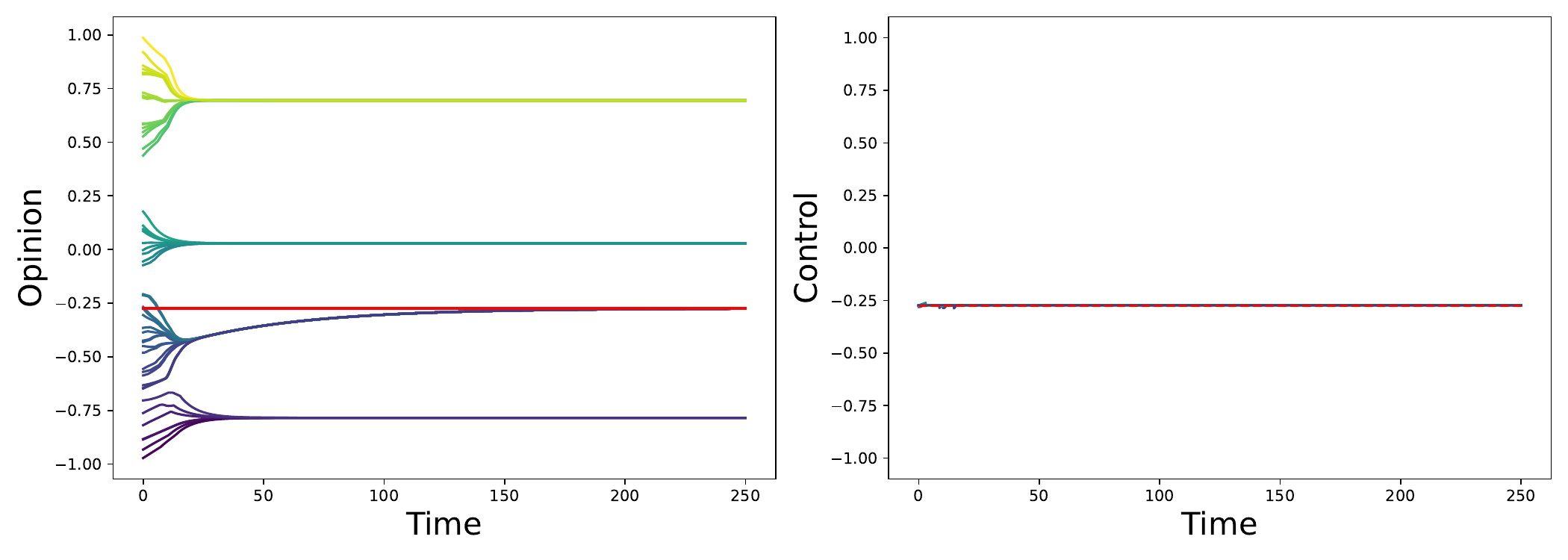}
        \caption{$\nu = 1$}
    \end{subfigure}
    \caption{Trajectory of agents (left) and the corresponding controls (right) in the case where the liar is penalised for deviating too far from their true opinion. Here we use a bounded confidence interaction kernel $P(x_j - x_i) = \chi(|x_j - x_i|\leq R)$ with $R=0.2$.}
    \label{fig:micro-model-bounded-conf}
\end{figure}

Essentially, if the optimal lie according to (\ref{eq:bounded-con-control}) is inside of agent $i$'s bounded confidence region, the liar tells this lie. If not, it is either optimal to project onto the boundary of the confidence region or tell the truth since $x_i^{n+1}$ will be unchanged by a lie outside of the confidence region and $(y_i - x_d)^2$ is minimised by taking $y_i = x_d$. Figure \ref{fig:micro-model-bounded-conf} displays an example of this control for a large and a small value of $\nu$. For small $\nu$, consensus at the liar's goal opinion is possible and the control strategy is to hug the edge of the bounded confidence interval of $x_i$. For larger $\nu$, consensus at $x_d$ may not be possible since it becomes too expensive to lie enough to place the liar in a position of influence for agents that are far away.

For a given $R\in(0,|\mathcal{I}|]$, we are interested in finding a range of $\nu$ such that consensus is possible. In other words, we are looking for some $\nu_{\text{crit}}>0$ such that for $\nu>\nu_{\text{crit}}$, consensus is not possible but for $\nu<\nu_{\text{crit}}$, consensus is guaranteed. A lower bound on $\nu_{\text{crit}}$ is to choose $\nu$ such that the liar is able to lie enough to reach the bounded confidence interval of every agent, we will let this value of $\nu$ be $\tilde{\nu}$.

\proposition{For all $0<\nu<\tilde{\nu}(R)$, consensus of the agents in the microscopic model (\ref{eq:xidyn}) at goal opinion $x_d\in\mathcal{I}$ is guaranteed when the control is given by equation (\ref{eq:bounded-con-control}) with the projection defined by (\ref{eq:projection-definition}). Here, 
\begin{equation}\label{eq:tildenu}
    \tilde{\nu} = \frac{\Delta t}{N}\frac{R|x_K^0 - x_d| - \frac{\Delta t}{N}R(\frac{1}{2}R-\text{sgn}(x_K^0-x_d)A_K^0)}{(R-|x_K^0 - x_d|)^2},
\end{equation}
where $K$ is the index of the agent who is initially the furthest away from the liar, i.e.,
\begin{equation}\label{eq:furthest-agent-bc}
     K = \arg\max_{j = 2, ..., N}\left\{|x_j^0 - x_d|\right\}.
\end{equation}
Furthermore, $A_K^0 = \sum_{j=2}^NP(x^0_j - x^0_K)(x^0_j - x^0_K)$ where $P(x_j - x_K)$ is the bounded confidence kernel with radius $R\in(0,|\mathcal{I}|]$ (\ref{eq:bounded-conf-interfun}).}\label{prop:nu-cond-boundedcon}\normalfont

\textbf{Proof}: We will assume that $R$ is such that $|x_K^0 - x_d|>R$ where $K$ is defined in equation (\ref{eq:furthest-agent-bc}), otherwise the result is trivial, see \cite{motsch2014heterophilious}. If the liar can communicate with agent $K$, then it can also communicate with all other agents so it is sufficient to find a condition such that the liar interacts with agent $K$. Since the model is contracting, agents can only move closer together at subsequent time-steps. Hence we only need to consider the dynamics at the very first time-step. In the rest of the proof, we will drop the superscript $0$'s so $x_K\equiv x_K^0$ etc. We want to find an upper bound for interaction so it is safe to say that the control will be
$$ y_K = x_K - \text{sgn}(x_K - x_d)R,$$
in order for the lie to be as close to the true opinion of $x_d$ as possible while still communicating with $x_K$. In order for this to be optimal, we also require that
\begin{equation}\label{eq:conditiononnu}
    (x_K^1 - x_d)^2 + 2\nu(y_K - x_d)^2<(\tilde{x}_K^1 - x_d)^2,
\end{equation}
where $x_K^1$ is the post-interaction opinion of agent $K$ when the liar tells lie $y_K$ and $\tilde{x}_K^1$ is the post-interaction opinion of agent $K$ when the liar tells their true opinion. Hence, we have that
\vspace{-0.2cm}
$$ x_K^1 = x_K + \frac{\Delta t}{N}\left(y_K - x_K + \sum_{j=2}^NP(x_j - x_K)(x_j - x_K)\right)$$
and
\vspace{-0.2cm}
$$ \tilde{x}_K^1 = x_K + \frac{\Delta t}{N}\sum_{j=2}^NP(x_j - x_K)(x_j - x_K).$$
Let $A_K = \sum_{j=2}^NP(x_j - x_K)(x_j - x_K)$. The inequality (\ref{eq:conditiononnu}) becomes
\vspace{-0.2cm}
$$ \left(x_K-x_d + \frac{\Delta t}{N}\left(y_K - x_K + A_K\right)\right)^2 + 2\nu(y_K - x_d)^2 < \left(x_K - x_d + \frac{\Delta t}{N}A_k\right)^2.$$
Rearranging this and cancelling terms gives us $\tilde{\nu}$ as in the proposition. $\blacksquare$

In order to complete our answer to the question of whether for given $R$ we can find a $\nu$ such that the liar can achieve consensus at $x_d$, we need to finally show that $\tilde{\nu}>0$ where $\tilde{\nu}$ is given by equation (\ref{eq:tildenu}). First, we note that $R|x_K^0-x_d|>R^2$ since we have that $|x_K^0 - x_d|>R$. Hence
$$ \tilde{\nu}> \frac{\Delta t}{N}\frac{R^2 - \frac{\Delta t}{N}(\frac{1}{2}R^2 - R\text{sgn}(x_K^0 - x_d)A_K^0)}{(R-|x_K^0 - x_d|)^2}.$$
We will consider only the case where $\text{sgn}(x_K^0 - x_d)=1$ (the other case follows similarly). By the definition of $K$, agent $K$ must be such that $x_K^0\geq x_i^0$ for all $i=1, ..., N$. Hence, we can write $A_K^0$ as
$$ A_K^0 = -\sum_{\{i: |x_K^0 - x_i^0|\leq R\}} |x_K^0 - x_i|\geq -NR.$$
Substituting this bound into the previous equation, we have that
$$\tilde{\nu}> \frac{\Delta t}{N(R - |x_K^0 - x_d|)^2}\left(1 - \frac{\Delta t}{2N} - \Delta t\right)>0$$
for 
$$ \Delta t< \frac{2N}{2N+1}.$$

\begin{figure}[htb!]
    \centering
    \includegraphics[width=0.9\textwidth]{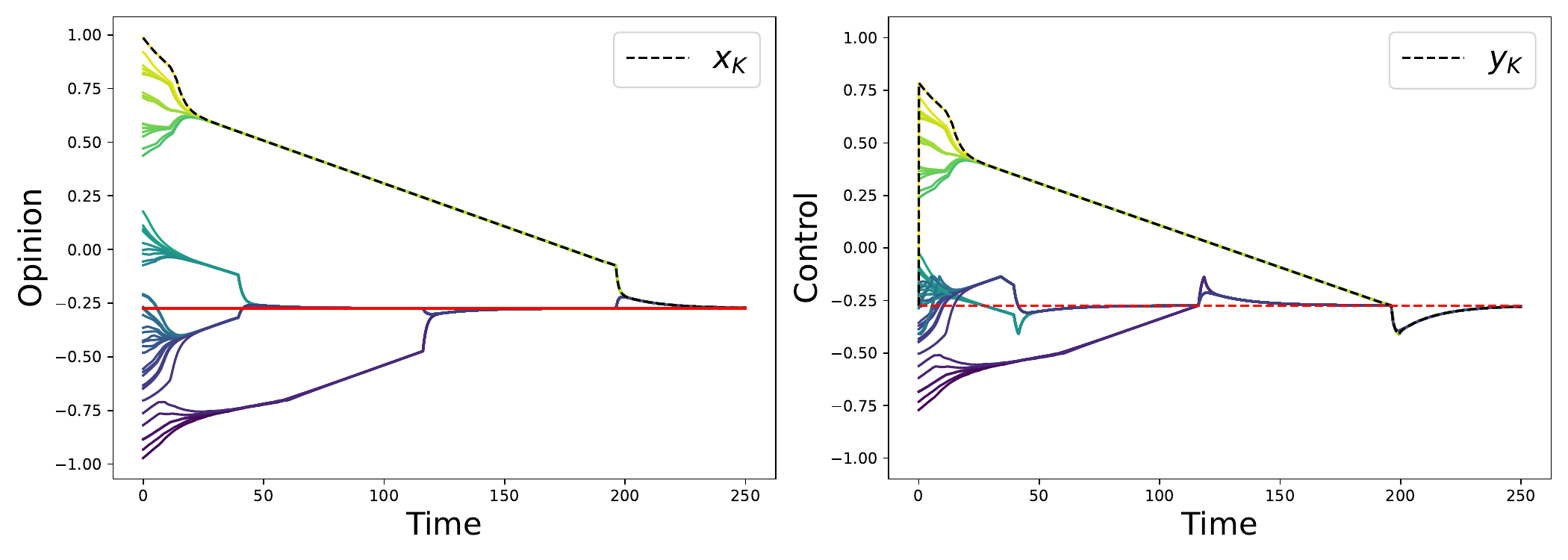}
    \caption{Trajectory of agents and the corresponding controls in the case where the liar cannot influence the most extreme individual initially. Here we use a bounded confidence interaction kernel $P(x_j - x_i) = \chi(|x_j - x_i|\leq R)$ with $R=0.2$ and we take $\nu = \tilde{\nu}\approx 0.000448$, computed according to equation (\ref{eq:tildenu}). The dashed line is agent $K$, the individual with initial opinion furthest from $x_d$.}
    \label{fig:micro-model-bounded-conf-nutilde}
\end{figure}

For $\nu<\tilde{\nu}$, the liar can reach any agent for all $t\geq 0$. For $\nu\geq \tilde{\nu}$, this is not the case as demonstrated in Figure \ref{fig:micro-model-bounded-conf-nutilde} where $\nu = \tilde{\nu}$. We see that consensus at $x_d$ is still possible. This is due to the fact that the most extreme agent $K$ is also influenced by the agents inside its confidence region who are all less extreme. After a short period of time, the liar is able to communicate with agent $K$ and influence them sufficiently to cause consensus at $x_d$. Obviously, if there are no agents within $R$ of $x_K^0$ then $\tilde{\nu} = \nu_{\text{crit}}$ however in general $\tilde{\nu}<\nu_{\text{crit}}$.

We can perform similar analysis as in the proof to Proposition \ref{prop:nu-cond-boundedcon} to show that $\tilde{\nu}(R)$ is an increasing function in $R$ for $\Delta t<\frac{2N}{2N+1}$.

\subsection{Kinetic model}\label{sec:bounded-con-kinetic}

We will now consider the limit of the system for large numbers of truth-telling agents $N\to \infty$. We define densities $f_T(x,t)$ and $f_L(x,t)$ as in Section \ref{sec:bin-interact} and have the identical Boltzmann equation, given by (\ref{eq:Boltzmann}).

We will consider binary interactions between truth-tellers, and interactions between truth-tellers and the liar as before. Given truth-telling agent opinions $x,x'$, we have post-interaction opinions given by
\begin{eqnarray*}
    x_{*} &=& x + \alpha P(x' - x)(x'-x) +  \Theta_1 D(x),\\
    x'_{*} &=& x' +\alpha P(x-x')(x-x') + \Theta_2D(x'),
\end{eqnarray*}
where $\alpha = \frac{\Delta t}{2}$, $P(x' - x)$ is the bounded confidence kernel (\ref{eq:bounded-conf-interfun}), $\Theta_1, \Theta_2$ are independent identically distributed normal random variables with mean 0 and variance $\sigma^2$ and $D(x)$ is a diffusion coefficient. Given a truth-telling agent with opinion $x''$, after interaction with the liar the truth-teller's opinion updates as 
$$ x''_{*} = x'' + \alpha P(y(x'')-x'')(y(x'')-x'') +\Theta'D(x''),$$
where $y(x'')$ is a function corresponding to the lie told to someone with opinion $x''$. We will continue to treat $y(x'')$ as a function and proceed with the Fokker-Planck derivation.

In the case of binary interactions, from equation (\ref{eq:SR-control}), we have the control function $y(x'')$ being given by
\begin{equation}
    y(x'') = x_d - \frac{\alpha}{\nu}(x'' - x_d+\alpha P(y(x'')-x'')(y(x'') - x''))\partial_{y}\left\{P(y(x'')-x'')(y(x'') - x'')\right\}.
\end{equation}
We proceed with the derivation and quasi-invariant limit as in Section \ref{sec:boltzmann}, finding the Fokker-Planck equation
\begin{equation}\label{eq:fok-plank-boundcon}
    \partial_tf_T(x,t) + \frac{1}{c_T}\partial_x(\mathcal{H}[f_T](x)f_T(x,t)) + \frac{\rho}{c_L}\partial_x(\mathcal{K}[f_T](x)f_T(x,t)) = \left(\frac{1}{c_T}+ \frac{\rho}{c_L}\right)\frac{\zeta}{2}\partial_{xx}(D(x)^2f_T(x,t)),
\end{equation}
where
\begin{subequations}
\begin{equation}\label{eq:convolution-term}
    \mathcal{H}[f_T](x) = \int_{\mathcal{I}}P(x'-x)(x' - x)f_T(x',t)\, dx',
\end{equation}
\begin{equation}
    \mathcal{K}[f_T](x) = P(y(x) - x)(y(x)-x),
\end{equation}
\end{subequations}
and $\kappa,c_L, c_T>0$ are rescaled regularisation and weight parameters, defined as in Section \ref{sec:quasi-invariant-limit}. This equation is supplemented with no-flux boundary conditions and an initial condition $f_T(x,0) =f_0(x)$. Under the quasi-invariant limit, the control must satisfy
\begin{equation}\label{eq:yeqn}
    y(x) = x_d - \frac{1}{\kappa}(x-x_d)\partial_y\left\{P(y(x) - x)(y(x)-x)\right\},
\end{equation}
projected onto $[-1,1]$.

\remark{We see that this is a general form for the Fokker-Planck equation with regularisation according to Section \ref{sec:Micro-model-std-reg}. Indeed, if we replace $P(y(x) -x)$ with $P(x)$ in equations (\ref{eq:fok-plank-boundcon}) and (\ref{eq:yeqn}), we return the Fokker-Planck equation discussed in Section \ref{sec:fokker-planck-easyP}.} \normalfont

Firstly, we will discuss a strategy to approximate $y(x)$ numerically. In order for the derivative $\partial_y\left\{P(y-x)\right\}$ to be defined, we need to take a smoothed version of the bounded confidence kernel (\ref{eq:bounded-conf-interfun}). As in \cite{nugent2024bridging,nugent2025opinion}, we take $r_1, r_2$ with $r_2>r_1>0$ and define $P(y-x)$ such that
\begin{equation}
    P(y-x) = \begin{cases}
        1 \quad &\text{for }|y-x|<r_1,\\
        \tilde{P}(y-x) \quad &\text{for }r_1\leq |y-x|\leq r_2,\\
        0 \quad &\text{for }|y-x|>r_2,
    \end{cases}
\end{equation}
where $\tilde{P}(y-x)$ is a function that smoothly interpolates between 1 and 0 in the region $r_1\leq |y-x|\leq r_2$ such that $P(y-x)$ is twice smoothly differentiable, i.e., $P\in C^2(\mathcal{I})$.

Then there are three distinct regions where our control takes different values:
\begin{equation}\label{eq:bounded-con-control-num}
    y(x) = \begin{cases}
        x_d- \frac{1}{\kappa}(x-x_d)\quad &\text{for }|y(x) - x|<r_1,\\
        x_d - \frac{1}{\kappa}(x-x_d)\partial_y\left\{\tilde{P}(y(x) - x)(y(x)-x)\right\} \quad &\text{for }r_1\leq |y(x) - x|\leq r_2,\\
        x_d\quad &\text{for }|y(x) - x|>r_2.
    \end{cases}
\end{equation}
The first case gives us a condition for $x$ very close to $x_d$:
$$ |x_d-x|<\frac{r_1\kappa}{\kappa+1}\implies y(x) = x_d - \frac{1}{\kappa}(x-x_d).$$
Given a discretisation of the opinion interval $\mathcal{I}$ of $L+1$ points, labelled $x_i$ for $i=0, ..., L$, we solve for $y(x_i)$. We start with the known region $|x_i - x_d|<\frac{r_1\kappa}{\kappa +1}$. We then use an inbuilt Newton solver for non-linear equations to find roots of the following equation
\begin{equation}\label{eq:discrete-numeric-solve}
    g(y;x_i) = \kappa(y - x_d) + (x_i - x_d)\partial_y\left\{P(y-x_i)(y-x_i)\right\}.
\end{equation}
We need to use a very fine discretisation for the results of the numerical solver to be accurate. This is due to the fact that there are multiple roots to $g(y;x_i)=0$ when $|x_i-x_d|>r_2$. Indeed, note that $y=x_d$ is always a root of $g(y;x_i)$ for $|x_i-x_d|>r_2$. The root finding algorithm is initialised using the optimal lie told at the closest discretisation point to $x_i$ in the direction of $x_d$. Details of the numerical method and an explanation of the fine discretisation can be found in Appendix \ref{app:y(x)-approximation}.

\begin{figure}[htb!]
    \centering
    \begin{subfigure}[b]{0.45\textwidth}
        \centering
        \includegraphics[width=\textwidth]{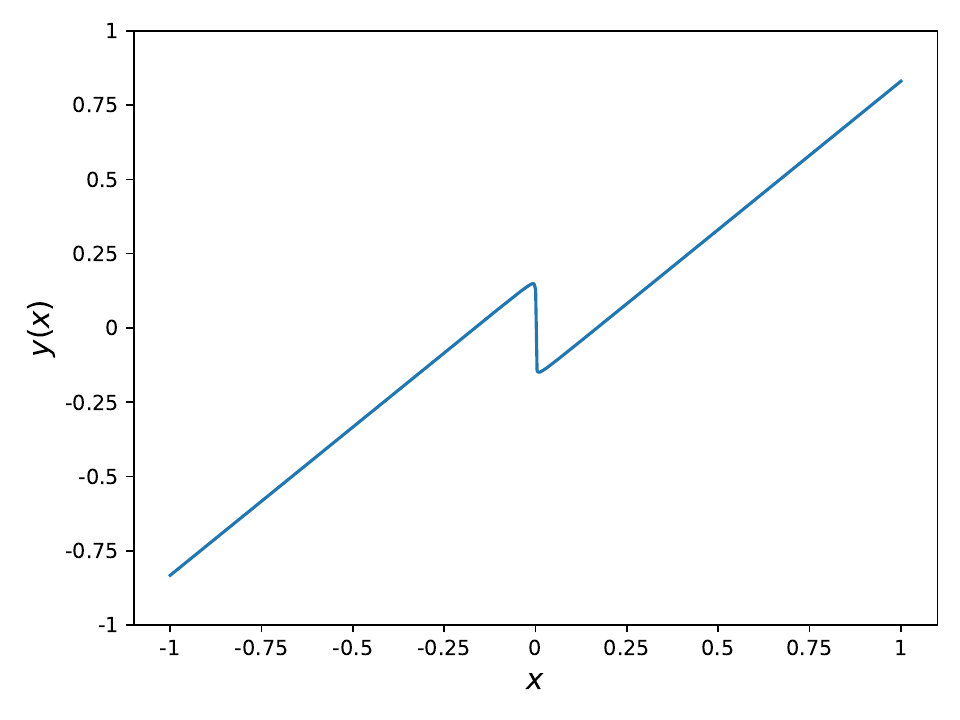}
        \caption{$\kappa = 0.01$, $r_1 = 0.1$, $r_2 = 0.2$}
        \label{fig:y(x)-control1}
    \end{subfigure}
    \begin{subfigure}[b]{0.45\textwidth}
        \centering
        \includegraphics[width=\textwidth]{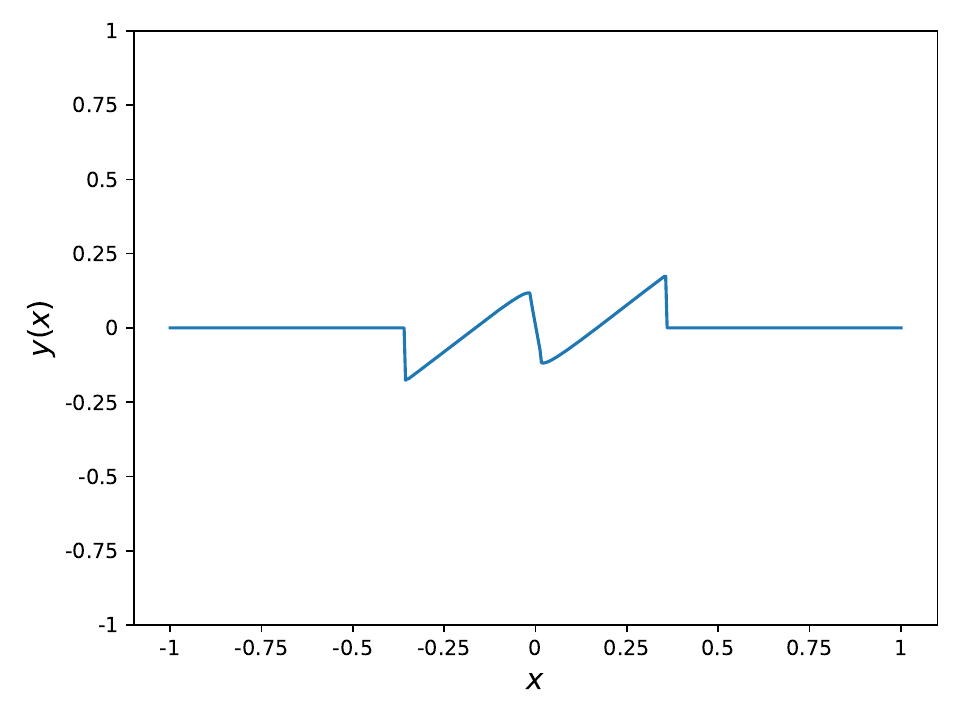}
        \caption{$\kappa = 0.15$, $r_1 = 0.1$, $r_2 = 0.2$}
        \label{subfig:y(x)-control-kappa0.15}
    \end{subfigure}
    \caption{Numerically computed $y(x)$ for different values of $\kappa$. Here we choose $x_d = 0$ and use the iterative solving method described in Appendix \ref{app:y(x)-approximation}.}
    \label{fig:y(x)-controlexamples}
\end{figure}
%\vspace{-0.5cm}

Figure \ref{fig:y(x)-controlexamples} shows the results of this scheme for different values of $r_1, r_2$ and $\kappa$. The scaling of $\kappa$ depends on the tolerance used for the numerical solver. We observe similar behaviour from the control as in the microscopic case. For small $\kappa$, the liar hugs the confidence region of agent $x$ in order to create the greatest rate of change of opinion $x$ in the direction of $x_d$. For larger $\kappa$, the liar decreases the amount by which they lie and can no longer communicate with agents of opinion that is very far from $x_d$. Furthermore, we note that in Figure \ref{subfig:y(x)-control-kappa0.15}, where $\kappa$ is comparatively large, there are local maxima and minima very close to $x_d$. This is due to the fact that within this region, the liar is able to communicate freely with an agent of opinion $x$ but can still lie without much cost in order to cause a greater increase/decrease to the opinion of a truth-teller with opinion $x$. Hence this strategy is more optimal than telling the truth.

\remark{Despite the appearance of Figure \ref{fig:y(x)-controlexamples}, $y(x)$ is typically not piecewise linear. Indeed, near to the peaks of the first local minima/maxima close to $x_d$ you can observe a slight curve. In the case where the interpolating polynomial $\tilde{P}(y-x) = 1-3t^2-2t^3$ for $t = \frac{|y-x|-r_1}{r_2-r_1}$, it can be shown that $y(x)$ is piecewise linear in the unregularised case $\kappa=0$ only.}\normalfont

\begin{figure}[htb!]
    \centering
    \begin{subfigure}[b]{0.45\textwidth}
        \centering
        \includegraphics[width=\textwidth]{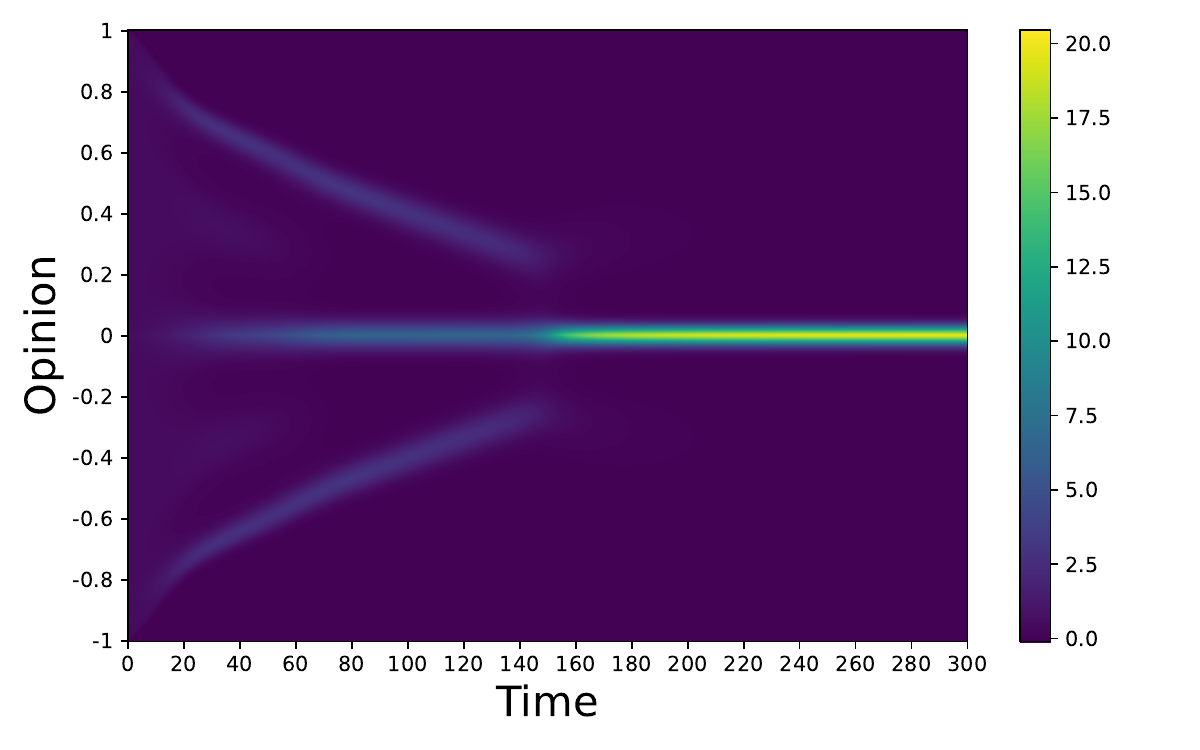}
        \caption{$\kappa = 0.01$}
    \end{subfigure}
    \begin{subfigure}[b]{0.45\textwidth}
        \centering
        \includegraphics[width=\textwidth]{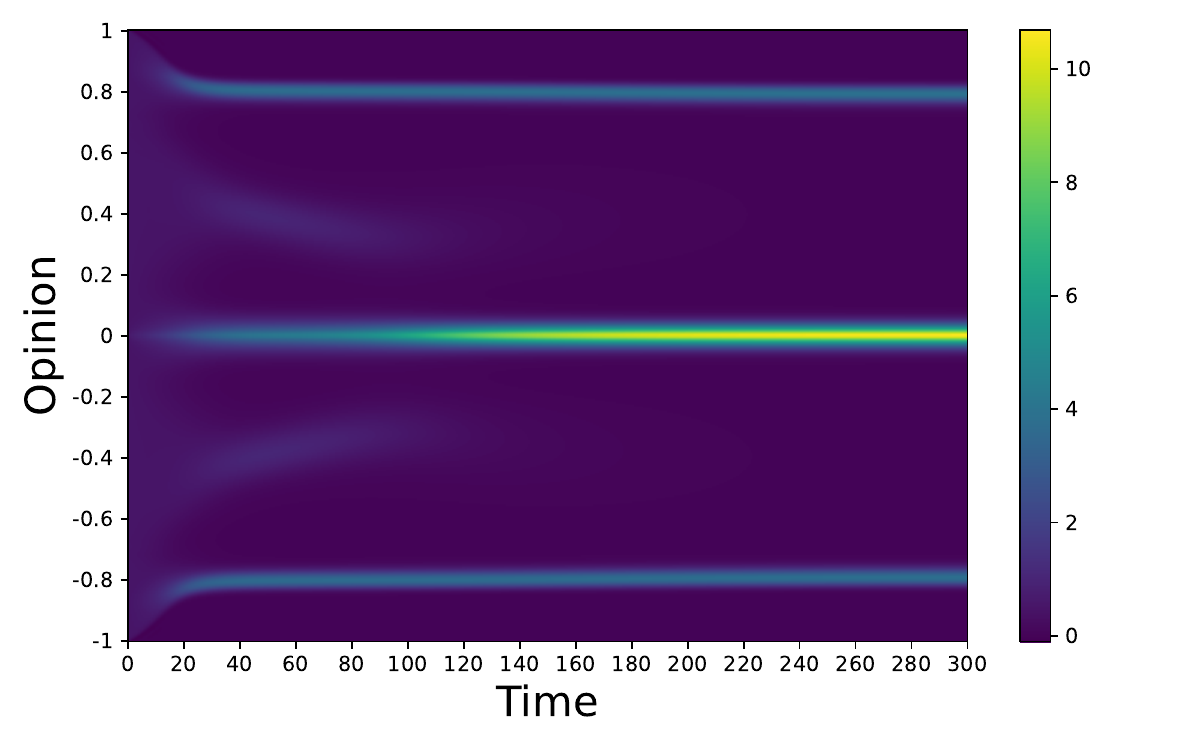}
        \caption{$\kappa = 0.1$}
    \end{subfigure}
    \begin{subfigure}[b]{0.45\textwidth}
        \centering
        \includegraphics[width=\textwidth]{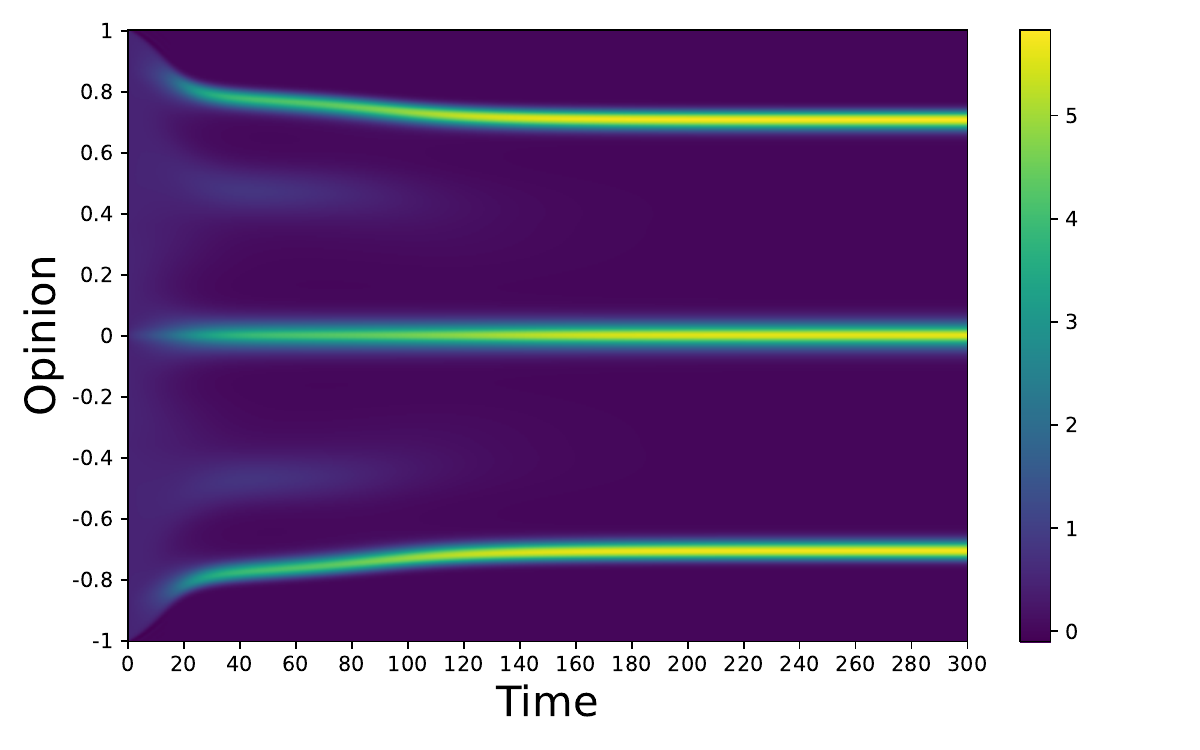}
        \caption{$\kappa = 0.15$}
    \end{subfigure}
    \begin{subfigure}[b]{0.45\textwidth}
        \centering
        \includegraphics[width=\textwidth]{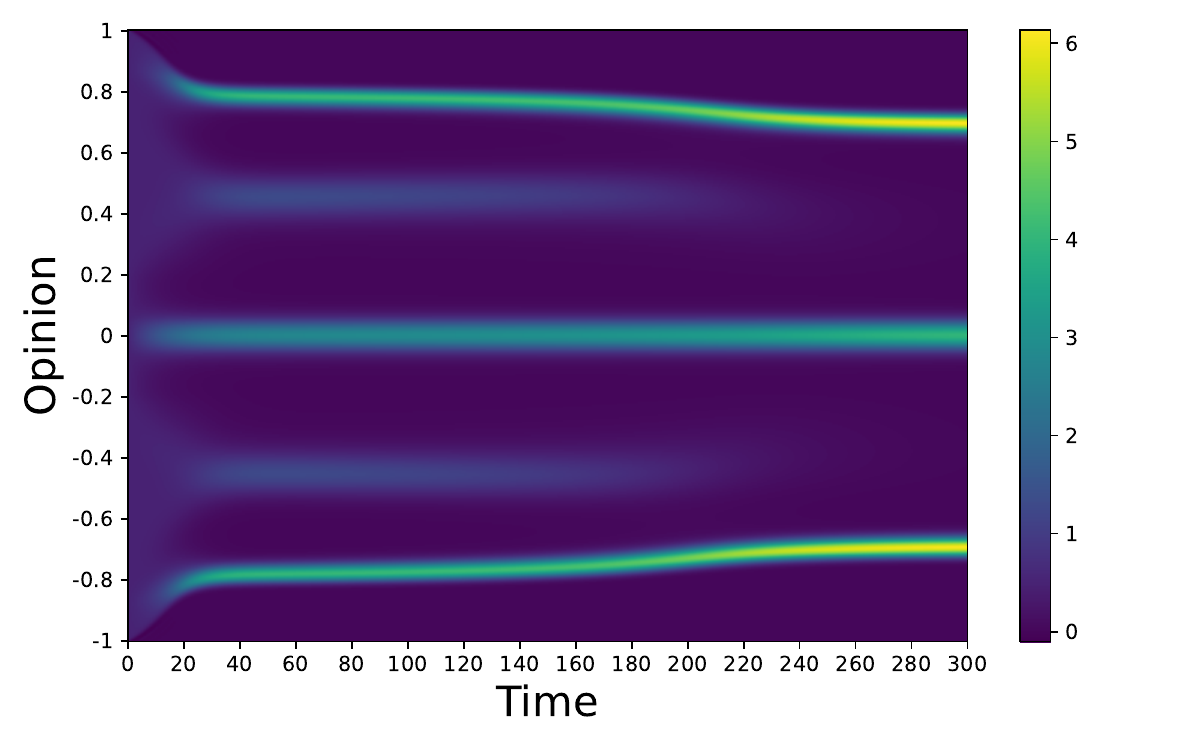}
        \caption{$\kappa = 1$}
    \end{subfigure}
    \caption{Solutions to the Fokker-Planck equation (\ref{eq:fok-plank-boundcon}) when $r_1 = 0.1$, $r_2 = 0.2$ and $x_d = 0$. Solutions are simulated using a finite volume scheme according to Appendix \ref{app:finite-volume} with $L=501$. The control $y(x)$ is computed numerically according to (\ref{eq:bounded-con-control-num}) with $2^7\times L$ grid points.}
    \label{fig:bounded-con-fin-volume}
\end{figure}

We can then incorporate our solution for $y(x)$ into a finite volume scheme to solve the Fokker-Planck equation (\ref{eq:fok-plank-boundcon}), the details of which can be found in Appendix \ref{app:finite-volume}. 

Figure \ref{fig:bounded-con-fin-volume} displays the evolution of the density $f_T(x,t)$ for different values of the regularisation parameter $\nu$. For very small $\nu$, we see that consensus at $x_d$ is achieved quickly. As $\nu$ increases, bands of consensus toward the more extreme ends of the opinion space grow stronger. Indeed, when $\nu = 1$, we have that the cluster of opinions around $x_d$ has lower density than the two clusters either side of $x_d$.

Furthermore, we conjecture that the function
\begin{equation*}
    V(t) = \left(m_T(t) - x_d\right)^2,
\end{equation*}
where $m_T(t) = \int_{\mathcal{I}}xf_T(x,t)\, dx$, is a Lyapunov function for our problem.
\begin{figure}[htb!]
    \centering
    \begin{subfigure}[b]{0.45\textwidth}
        \centering
        \includegraphics[width=\textwidth]{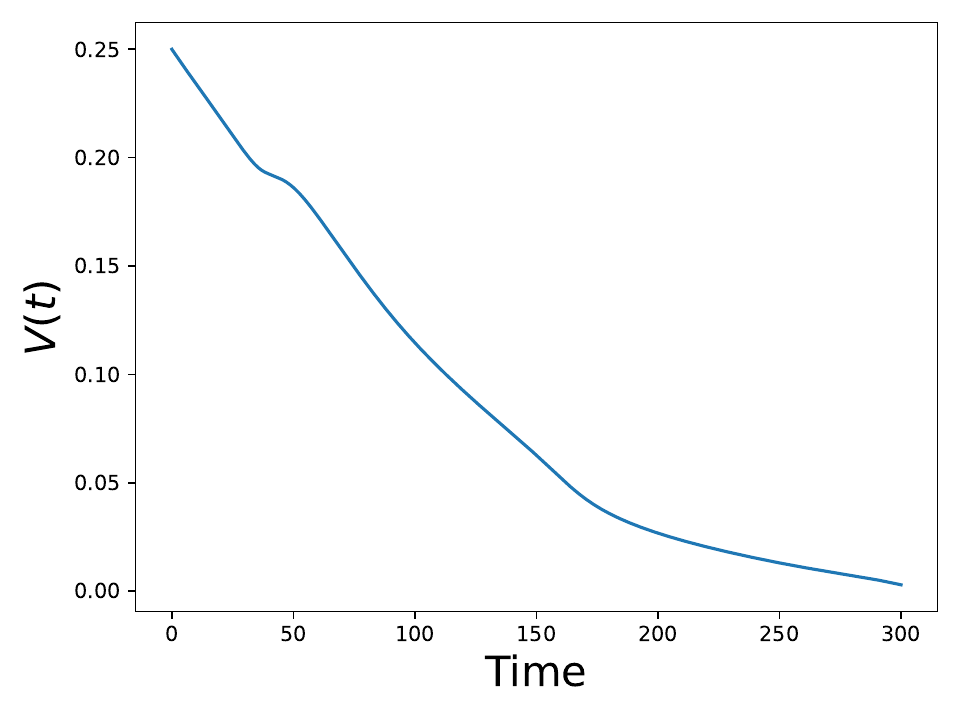}
        \caption{$\kappa = 0.01$}
    \end{subfigure}
    \begin{subfigure}[b]{0.45\textwidth}
        \centering
        \includegraphics[width=\textwidth]{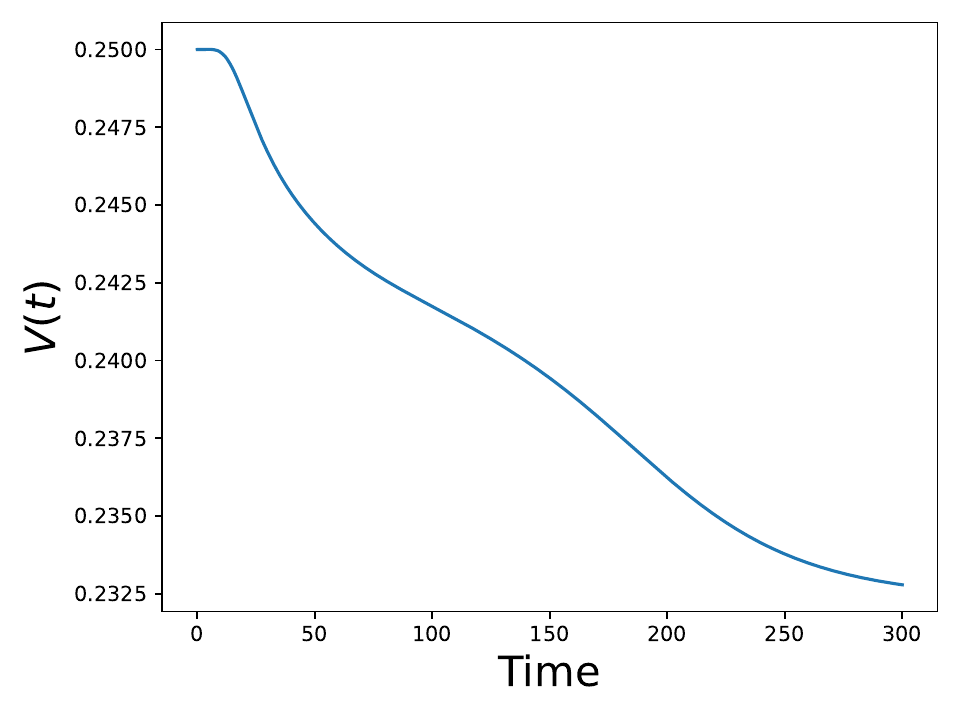}
        \caption{$\kappa = 1$}
    \end{subfigure}
    \caption{Plots showing the decay of $V(t)$ for a small and a large value of $\kappa$ respectively. Here, we take $x_d = -0.5$ and calculate $m_T(t)$ with numeric integration based on the samples of the function we collect from our finite volume discretisation.}
    \label{fig:lyapunov-func}
\end{figure}
Figure \ref{fig:lyapunov-func} demonstrates that for both large and small values of $\kappa$, $V(t)$ is found to be decreasing. Clearly, for small $\kappa$, the steady state of $V(t)$ will be a value much closer to 0 compared with larger $\kappa$.

\remark{Proving analytically that $V(t)$ is a Lyapunov function is challenging. Indeed, since $V(t)\geq0$ we need only show that $\dot{V}(t)\leq 0$. Now, $\dot{V}(t) = 2\frac{d}{dt}(m_T(t))(m_T(t)-x_d)$ where
$$\frac{d}{dt}m_T(t) = \frac{\rho}{c_L}\int_{\mathcal{I}}P(y(x)-x)(y(x)-x)f_T(x,t)\, dx.$$
Substituting equation (\ref{eq:yeqn}) for $y(x)$, we have
$$\frac{d}{dt}m_T(t) = -\frac{\rho}{c_L}\int_{\mathcal{I}}P(y-x)\left[1+\frac{1}{\kappa}P(y-x)+\frac{1}{\kappa}P'(y-x)(y-x)\right](x-x_d)f_T(x,t)\, dx,$$
where $y=y(x)$ and $P'(y-x) = \partial_y\{P(y-x)\}|_{(x,y(x))}$. The first two terms of the sum in the square bracket are bounded below by zero, however the third term is always negative. Indeed, if $x>x_d$ then $y<x$ (else the lie would not be optimal) and therefore $P'(y-x)>0$ by the oddness of $P'$. Hence $P'(y-x)(y-x)<0$ and similarly for $x<x_d$.
}\normalfont

\section{The effect of multiple liars}\label{sec:multiple-liars}

We now suppose that more than one population of liars exists, each one with a different strategy and goal opinion. We can describe the evolution of the kinetic density of the system through a Boltzmann approach by a similar method to Albi et al. in \cite{albi2014boltzmann}. Let $K>0$ be the number of groups of liars, each with some goal opinion $x_{d_p}$, relative volume $\rho_p$ and density $f_{L_p}$ for $p = 1, ..., K$ such that
$$ f_{L_p}(x,t) = \rho_p\delta(x - x_{d_p}), \quad \forall t\geq 0.$$
We suppose that a unique population of truth-tellers exists, and that every truth-teller interacts both with other truth-tellers and every family of liars $p=1, ..., K$. Then for a suitable test function $\varphi$, we have the following Boltzmann-type equation for the evolution of the density of truth-telling agents $f_T(x,t)$,
\begin{equation*}
    \frac{d}{dt}\int_{\mathcal{I}}\varphi(x)f_T(x,t)\, dx = (Q_T(f_T, f_T), \varphi) + \sum_{p=1}^K(Q_L(f_T, f_{L_p}), \varphi),
\end{equation*}
where $(Q_L(f_T, f_{L_p}), \varphi)$ is written analogously to equation (\ref{eq:boltz-rhs}).

When $P(x)=\mathcal{P}$ is constant, we can perform some similar analysis as in the single liar case in Section \ref{sec:boltzmann-eqnanal} to find the mean opinion of the truth-telling agents at large times. We will let each liar have interaction frequency $\eta_{L_p}$ and regularisation parameter $\nu_{p}$. Then, by taking $\varphi(x) = x$, similarly to equation (\ref{eq:mean-oneliar}), we have that the mean opinion of the truth-telling agents evolves according to
\begin{equation*}
    \frac{d}{dt}m_T(t) = \sum_{p=1}^K\left(\rho_p\eta_{L_p}\left(\frac{\nu_p\alpha \mathcal{P} +\alpha^2 \mathcal{P}^2}{\nu_p + \alpha^2\mathcal{P}^2}\right)(-m_T(t) + x_{d_p})\right).
\end{equation*}
We can solve this ODE exactly to obtain that the mean opinion at time $t\geq 0$ is 
\begin{equation*}
    m_T(t) = C \exp\left(-\sum_{p=1}^K\left(\rho_{p}\eta_{L_p}\frac{\nu_p\alpha \mathcal{P} + \alpha^2\mathcal{P}^2}{\nu_p + \alpha^2 \mathcal{P}^2}\right)t\right) + \frac{\sum_{p=1}^K\rho_p\eta_{L_p}\frac{\nu_p + \alpha \mathcal{P}}{\nu_p + \alpha^2 \mathcal{P}^2}x_{d_p}}{\sum_{p=1}^K\rho_p\eta_{L_p}\frac{\nu_p + \alpha \mathcal{P}}{\nu_p + \alpha^2 \mathcal{P}^2}},
\end{equation*}
for a constant $C$ determined by initial conditions. Note that in the particular case where $\nu_p = \nu$ for $p=1, ...,K$, we have that 
\begin{equation}\label{eq:Kliars-weightedmean}
    \lim_{t\to\infty}m_T(t) = \frac{\sum_{p=1}^K\rho_p\eta_{L_p}x_{d_p}}{\sum_{p=1}^K\rho_p\eta_{L_p}},
\end{equation}
i.e., the mean of the truth-tellers' opinions tends to a weighted average of the goal opinions of the liars. The weighting given to a liar is greater if they communicate more with the truth-tellers, either by having greater relative volume $\rho_p$ or greater interaction frequency $\eta_{L_p}$.

\begin{figure}[htb!]
    \centering
    \begin{subfigure}[b]{0.45\textwidth}
        \centering
        \includegraphics[width=\textwidth]{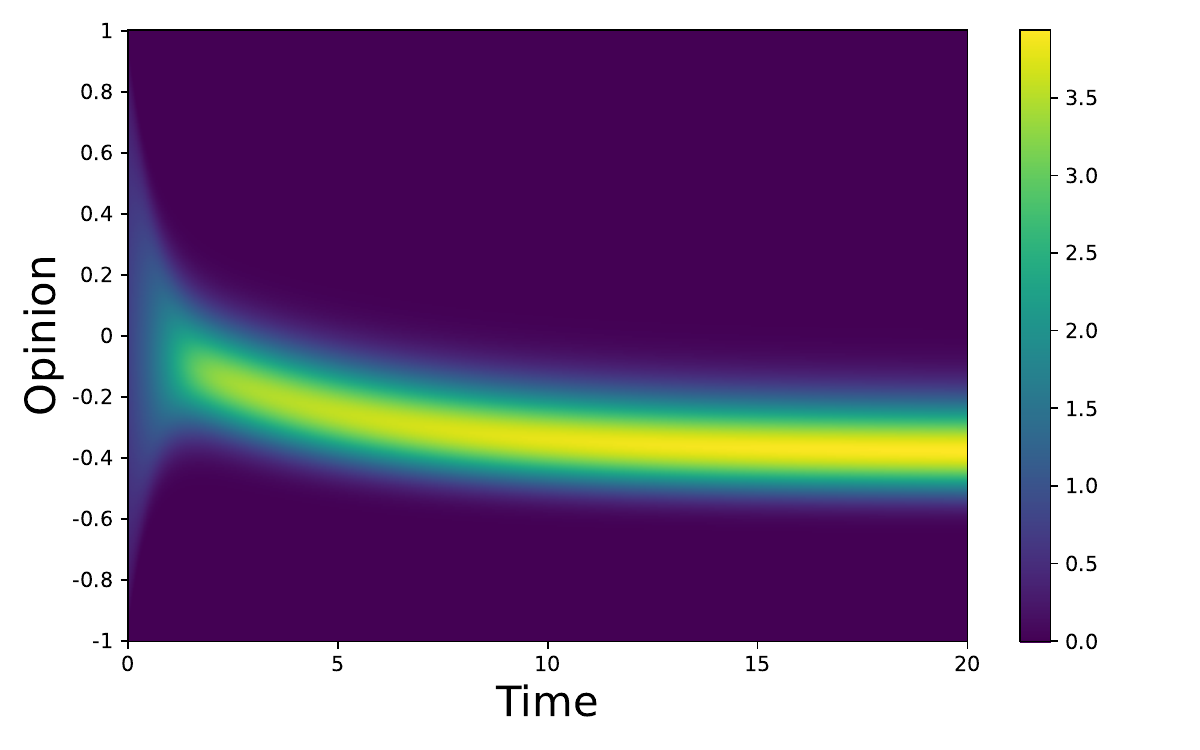}
        \caption{$f_T(0)$ uniform}
    \end{subfigure}
    \begin{subfigure}[b]{0.45\textwidth}
        \centering
        \includegraphics[width=\textwidth]{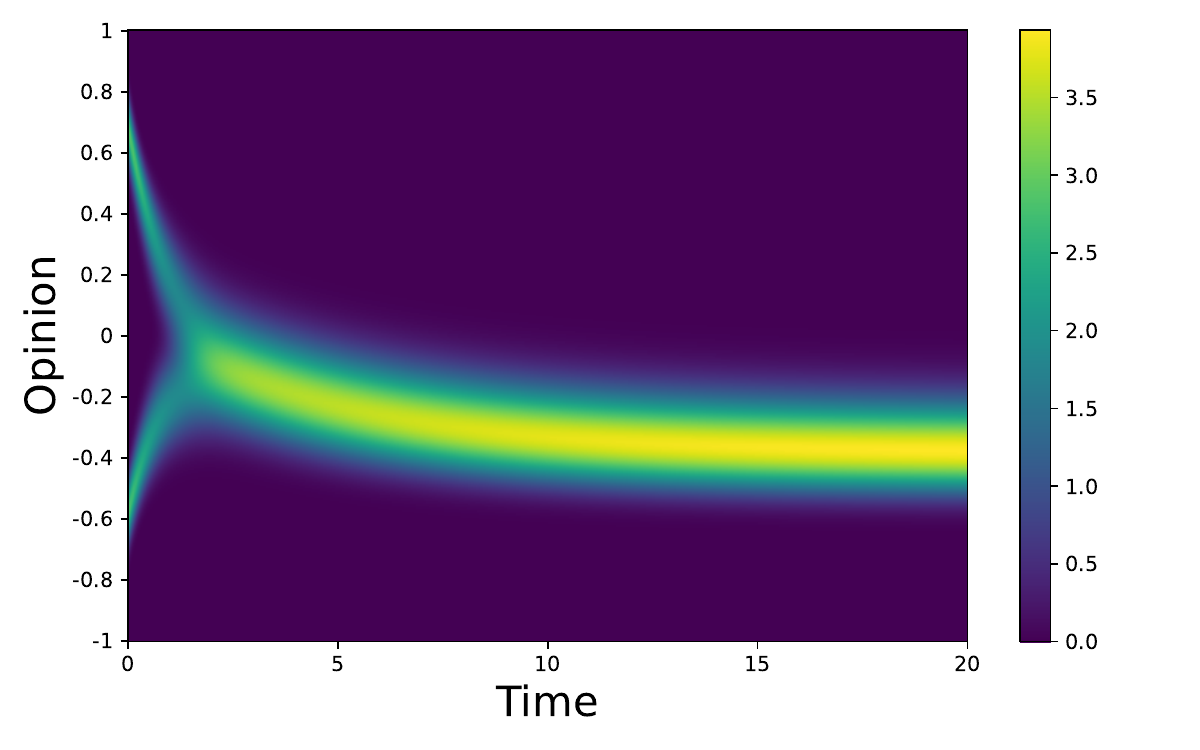}
        \caption{$f_T(0)$ non-uniform}
    \end{subfigure}
    \caption{Simulations for the truth-tellers' density $f_T(x,t)$ in the presence of two liars. The liars have goal opinions $x_{d_1} = 0.5$ and $x_{d_2} = -0.5$. We take relative volumes $\rho_{1} = 0.02$ and $\rho_{2} = 0.1$. All other parameters are set to be equal, $\kappa_1=\kappa_2 = 0.1$ and $c_{L_1} = c_{L_2} = 1$ and $P(x,\cdot) \equiv 1$. Both simulations are carried out by a finite volume solver for the corresponding Fokker-Planck equation.}
    \label{fig:monte-carlo-2liars}
\end{figure}

Figure \ref{fig:monte-carlo-2liars} displays simulations for the case where we have two liars, one with a goal opinion $x_{d_1} = 0.5$ and the other with a goal opinion $x_{d_2} = -0.5$. The first liar has relative volume $\rho_1 =0.02$, the second has relative volume $\rho_2 = 0.1$ and we take $\eta_{L_1} = \eta_{L_2}$. We see that in this simple case where $P(x) = 1$, the steady state solution reached is a compromise between these two goal opinions, weighted by their relative influences. By equation (\ref{eq:Kliars-weightedmean}), the mean of the steady state distribution is weighted by the relative volumes of the liars, i.e.,
$$\lim_{t\to\infty} m_T(t) = \frac{\rho_1}{\rho_1+\rho_2}x_{d_1} + \frac{\rho_2}{\rho_1 + \rho_2}x_{d_2}.$$
This is verified in Figure \ref{fig:monte-carlo-2liars}, where we have the expected long-time mean opinion, $-1/3$. We also consider the case where the initial condition on the truth-teller's opinions is non-uniform and clustered. Here, 
$$ f_T(x,0) = Z\left(e^{-100(0.7-x)^2}+e^{-100(-0.6-x)^2}\right),$$
where $Z$ is a normalising constant. As expected, both solutions approach identical steady state distributions at large times. 

When we introduce a bounded confidence kernel, where agents only interact with those who have a sufficiently similar opinion to their own, we are able to have clusters forming. To explore this regime, we first consider an extension to the Fokker-Planck equation (\ref{eq:fok-plank-boundcon}) where the truth tellers are under the influence of multiple liars. Using the same notation as in Section \ref{sec:bounded-con-kinetic}, we are able to write the following Fokker-Planck equation
\begin{multline}\label{eq:fok-planck-mult-liars}
    \partial_tf_T(x,t) + \frac{1}{c_T}\partial_x(\mathcal{H}[f_t](x)f_T(x,t)) + \sum_{k=1}^K\frac{\rho_k}{c_{L_k}}\partial_x(\mathcal{K}_k[f_T](x)f_T(x,t)) \\= \left(\frac{1}{c_T}+\sum_{k=1}^K\frac{\rho_k}{c_{L_k}}\right)\frac{\zeta}{2}\partial_{xx}(D(x)^2f_T(x,t)),
\end{multline}
where
$$ \mathcal{K}_k[f_T](x)= P(y_k(x) -x)(y_k(x) - x),$$
$\mathcal{H}[f_T](x)$ is given by equation (\ref{eq:convolution-term}) and $y_k(x)$ is given by (\ref{eq:bounded-con-control-num}) with goal opinion $x_{L_k}$ and regularisation parameter $\kappa_{k}$.

\begin{figure}[htb!]
    \centering
    \begin{subfigure}[b]{0.45\textwidth}
        \centering
        \includegraphics[width=\textwidth]{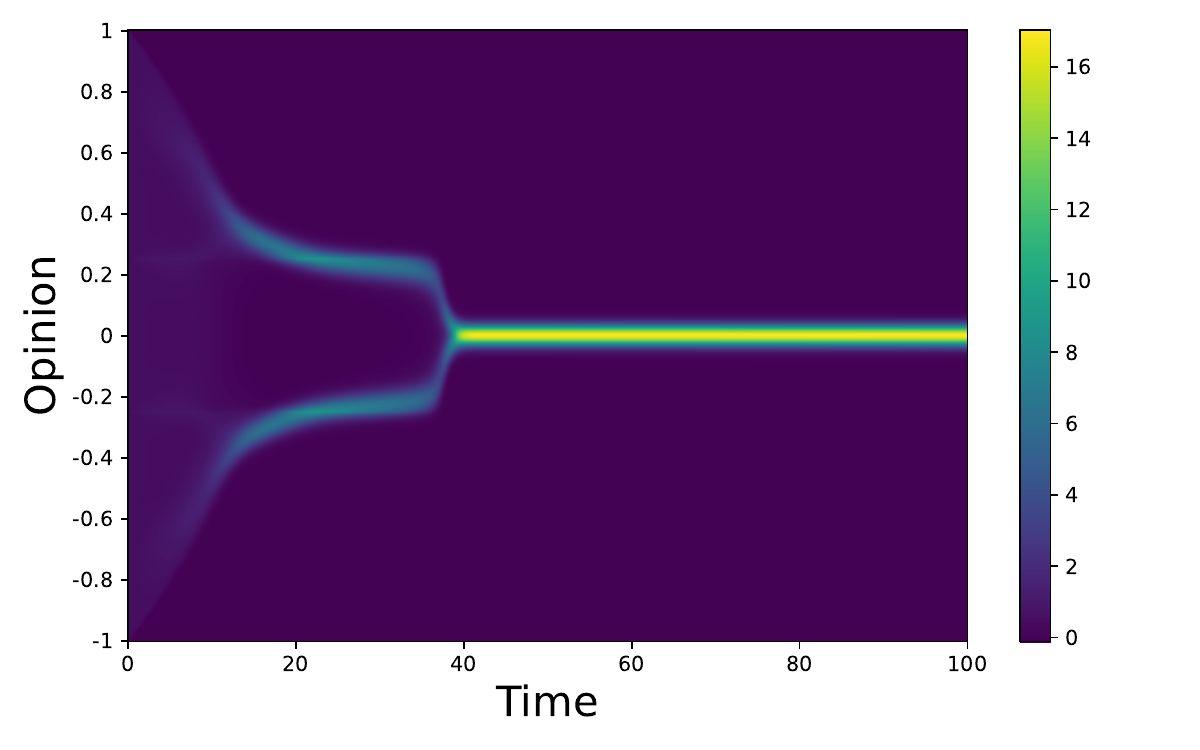}
        \caption{$\kappa_1 = \kappa_2 = 0.01$, $r_1=0.2$, $r_2 = 0.4$}
        \label{subfig:nuequal-twoliar-ev}
    \end{subfigure}
    \begin{subfigure}[b]{0.34\textwidth}
        \centering
        \includegraphics[width=\textwidth]{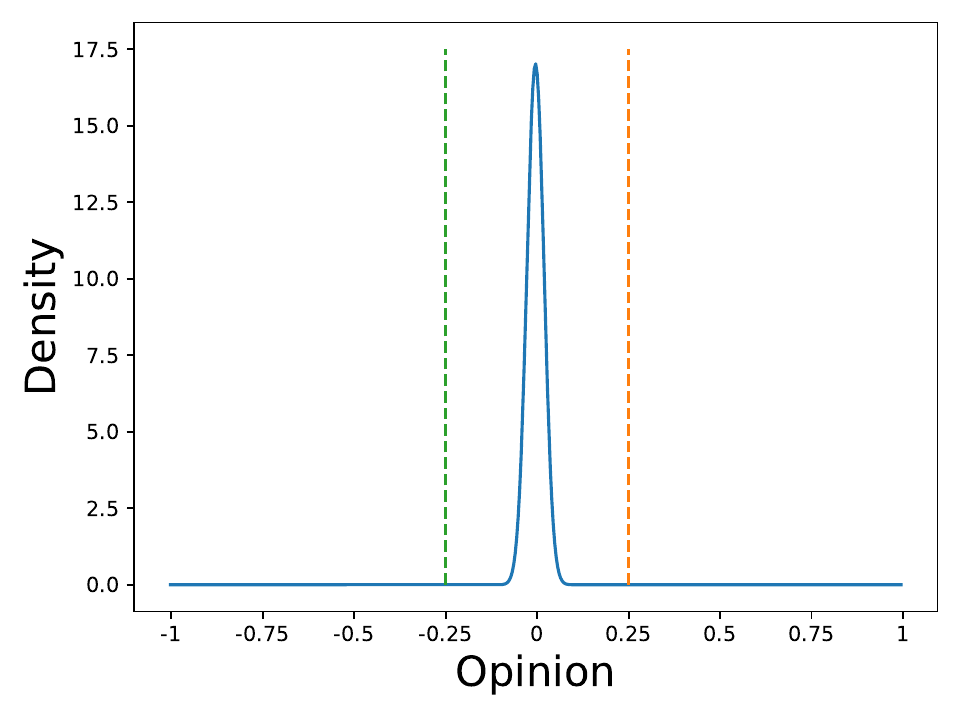}
        \caption{\centering Steady state $\kappa_1 = \kappa_2 = 0.01$, $r_1 = 0.2$, $r_2 = 0.4$}
        \label{subfig:nuequal-twoliar-final}
    \end{subfigure}
    \begin{subfigure}[b]{0.45\textwidth}
        \centering
        \includegraphics[width=\textwidth]{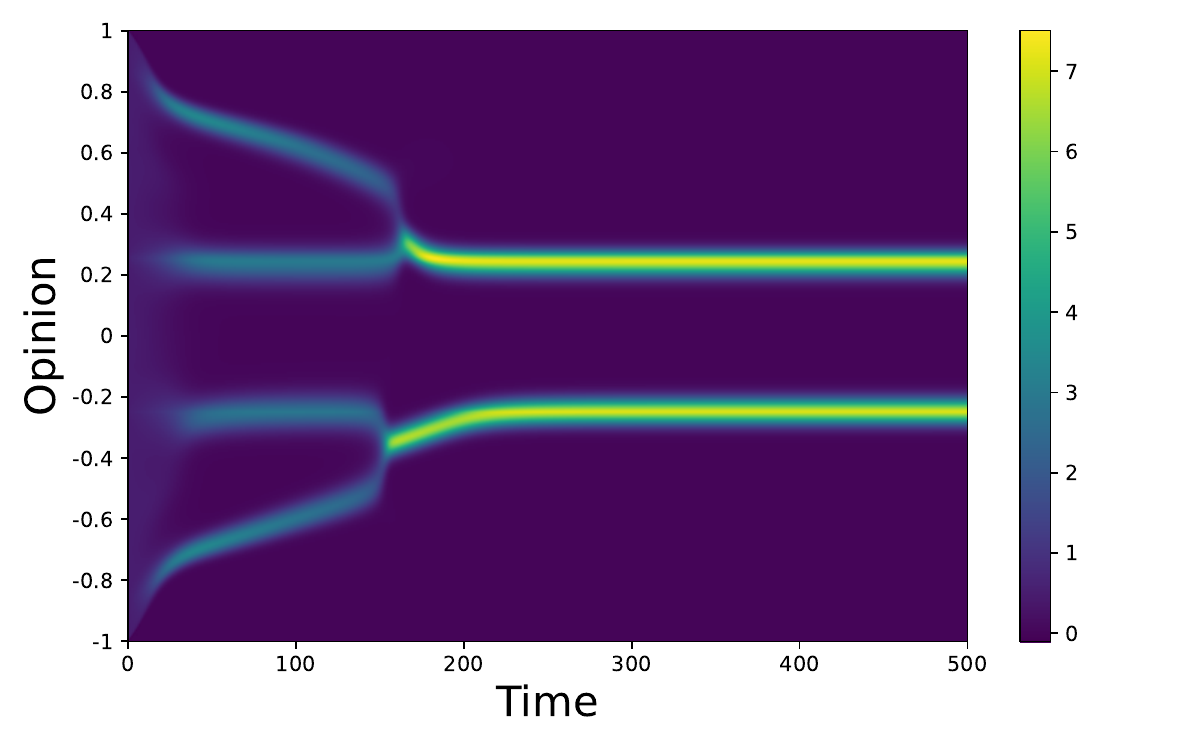}
        \caption{$\kappa_1 = 0.15$, $\kappa_2 = 0.01$, $r_1 = 0.1$, $r_2 = 0.2$}
    \end{subfigure}
    \begin{subfigure}[b]{0.34\textwidth}
        \centering
        \includegraphics[width=\textwidth]{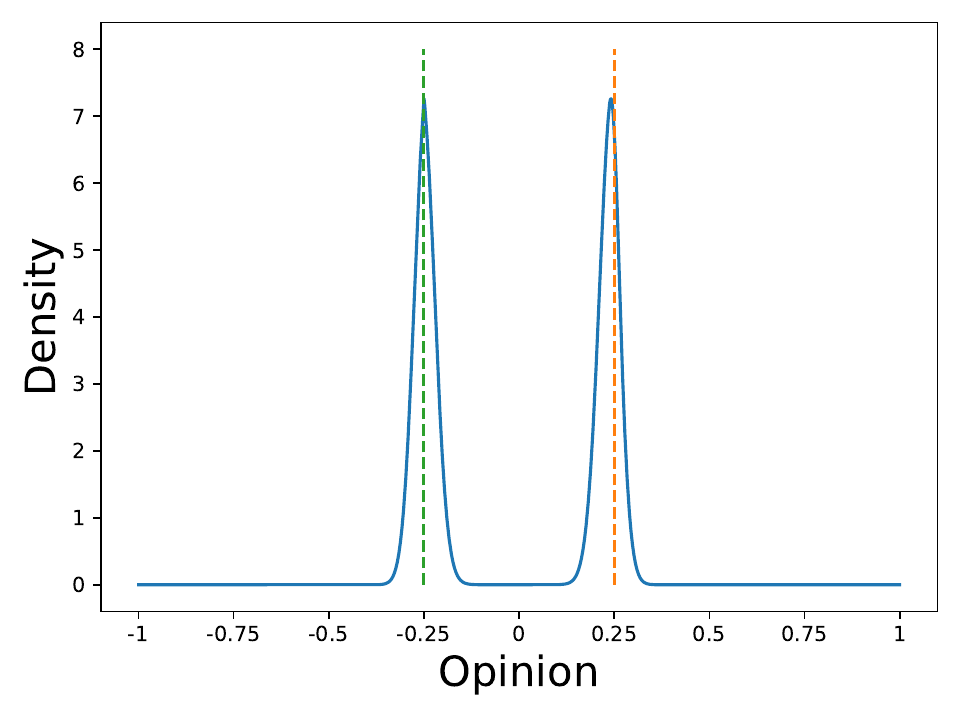}
        \caption{\centering Steady state $\kappa_1 =0.15$, $\kappa_2 = 0.01$, $r_1 = 0.1$, $r_2 = 0.2$}
    \end{subfigure}
    \begin{subfigure}[b]{0.45\textwidth}
        \centering
        \includegraphics[width=\textwidth]{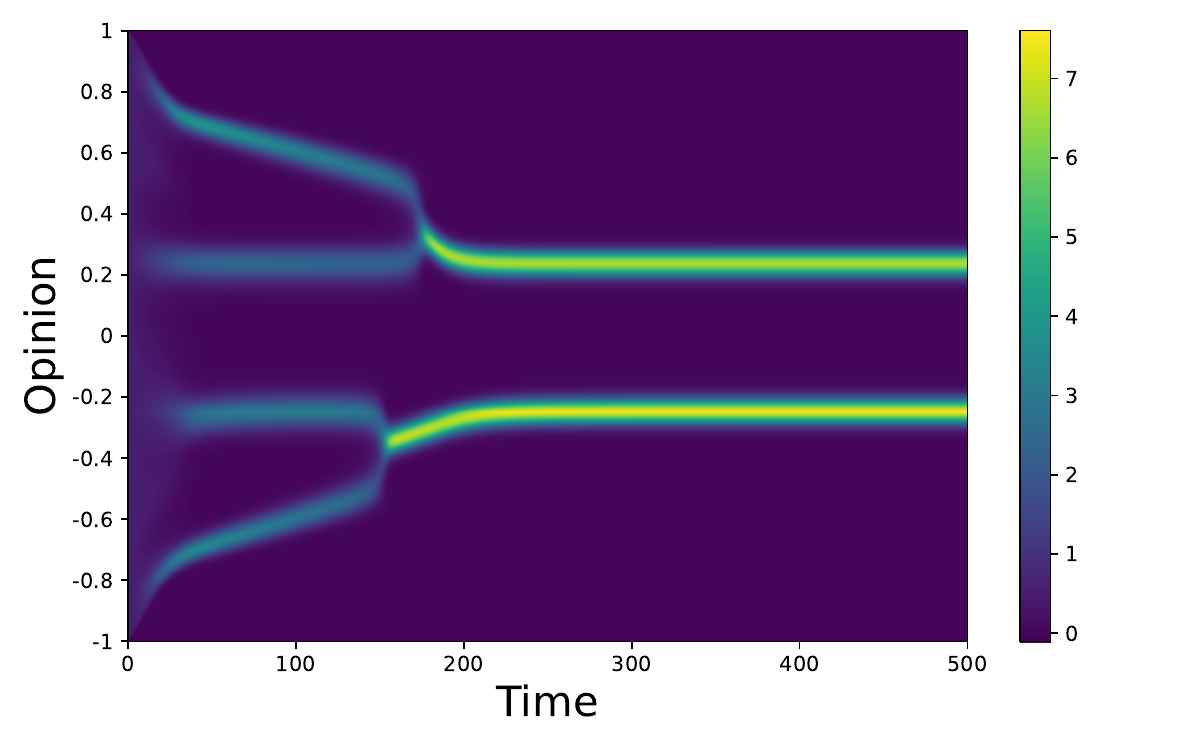}
        \caption{$\kappa_1 = 1$, $\kappa_2 = 0.01$, $r_1=0.1$, $r_2 = 0.2$}
    \end{subfigure}
    \begin{subfigure}[b]{0.34\textwidth}
        \centering
        \includegraphics[width=\textwidth]{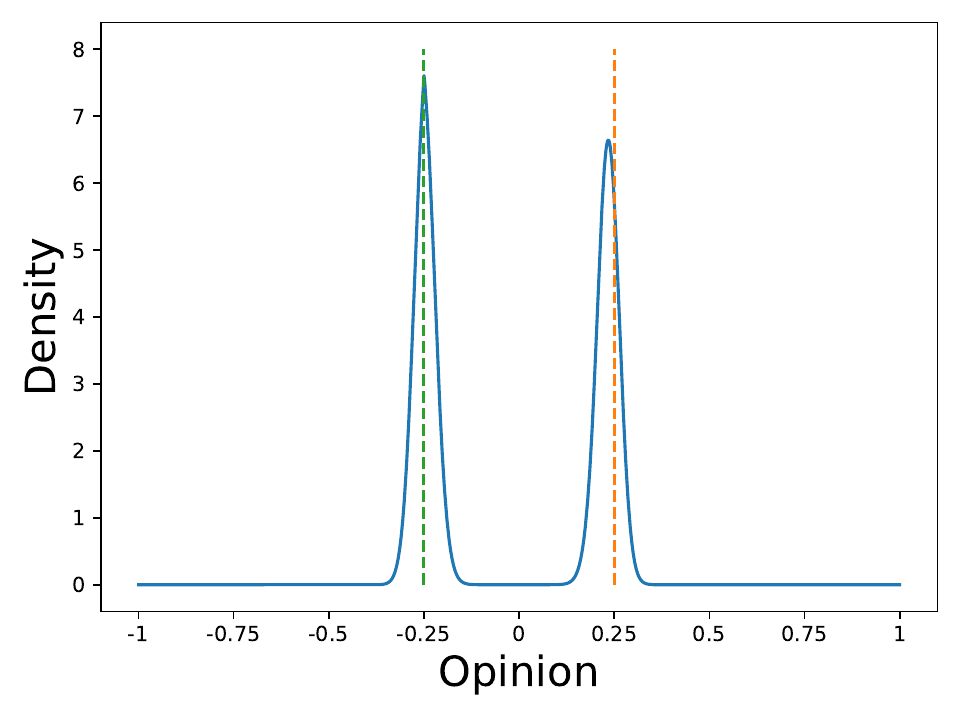}
        \caption{\centering Steady state $\kappa_1 =1$, $\kappa_2 = 0.01$, $r_1 = 0.1$, $r_2 = 0.2$}
    \end{subfigure}
    \caption{Numerical approximation to the solution of the Fokker-Planck equation (\ref{eq:fok-planck-mult-liars}) with a bounded confidence interaction kernel in the presence of two liars. Here we take $x_{d_1} = 0.25$, $x_{d_2} = -0.25$, $\rho_1 = \rho_2 = 0.05$ $\kappa_2 = 0.01$ and vary $\kappa_1$.}
    \label{fig:two-liars-bounded-conf}
\end{figure}
In Figure \ref{fig:two-liars-bounded-conf}, we display three simulations for the evolution of $f_T(x,t)$ in the presence of two liars in a bounded confidence regime. The discretisation method is similar to that outlined in Appendix \ref{app:finite-volume} but with extra terms corresponding to the influence of multiple liars. In Figures \ref{subfig:nuequal-twoliar-ev} and \ref{subfig:nuequal-twoliar-final}, we have that $\kappa_1=\kappa_2$ so both liars have the same amount of regularisation applied to their strategies. The truth-tellers reach consensus at the midpoint of $x_{d_1}, x_{d_2}$, i.e., at $x=0$. This is a \textit{compromise} and may only occur if $1/2|x_{d_1} - x_{d_2}|<r_2$. Increasing $\kappa_1$ decreases the amount by which the first liar is willing to lie. Interestingly a weaker liar can do better in the presence of another, stronger liar than it would do if it was a solo liar with the same regularisation parameter $\kappa_i$. The stronger liar can communicate with more agents than the weaker liar can and so accidentally brings some of the truth-telling agents into the confidence region of the weaker liar. In the corresponding steady state solutions, we see polarisation for $\kappa_1>\kappa_2$, where the cluster closest to the opinion $x_{d_1}$ is shifted slightly toward $x_{d_2}$ in a small compromise toward the stronger liar $L_2$. We see that for increased $\kappa_1$, the density close to $x_{d_1}$ at late times also decreases. 

\section{Conclusion}\label{sec:conclusions}

In this work, we have explored a novel approach to optimal consensus control with a strong real-world interpretation. We began by presenting the microscopic problem for liars and a method for seeking instantaneous solutions through a receding horizon optimal control strategy. Using this method, we presented several regularisation strategies for the liar, each inspired by some social convention or risk. We demonstrated how the model could be extended to a large-population limit and performed analysis on the corresponding Boltzmann and Fokker-Planck equations along with numerical simulation.

An alternative formulation for optimal consensus control by lying would be to assume that we have $K$ agents who are allowed to lie but that they can only express one apparent opinion to the whole population at each time-step. With this model, there are some interesting questions that could be addressed. For example, for the bounded confidence model (discussed in Section \ref{sec:bounded-conf}) with a certain confidence radius $R$, how many liars is it necessary to have such that consensus is reached at some desired opinion? Following similar methods to \cite{burger2020instantaneous}, it would also be possible to create a mean-field description of this model. Another possible direction for further research would be to consider a control framework in which the goal of the liar is to \textit{prevent} consensus, following similar strategies to D{\"u}ring et al. in \cite{during2024breaking}. This question is important since in this paper we have only discussed the best way to manipulate and deceive, which is certainly a morally questionable pursuit. Furthermore, the question of how to model \textit{trust} in these types of formulations is of interest \cite{urena2019review}. In considering trust, we could potentially model how a liar could be caught lying and the subsequent consequences of exposing a liar on the dynamics of consensus.

The introduction of control through lying produces interesting mathematical questions as well as qualitatively applicable strategies and results. This formulation carries two main differences from previous formulations of optimal consensus control. The first is that the control comes from within the system; there is no external figure deciding what the goal opinion or methods of influence are. Through this contrast, control through lying presents an intuitive method for controlling a population without leaving the population itself. The simplicity of the system also means that setting up experiments aiming to capture some of the behaviour seen here might be fruitful. Experiments could verify some of the strategies observed in this paper and help to explain under what conditions different regularisation become important. A second difference to previous formulations we have presented here is the complex dependency on the derivative of the interaction function $P(x,y)$ on the control. While this makes the problem challenging, we also see that the strategies employed by the liar are sophisticated and intuitive.  

\section{Acknowledgements}

We would like to thank Andrew Nugent for his insightful comments on the draft version of this paper and his thoughts during early conversations about this work.

SG is supported by the Warwick Mathematics Institute Centre for Doctoral Training, and gratefully acknowledges funding from the University of Warwick and the UK Engineering and Physical Sciences Research Council (Grant number: EP/W524645/1).

MTW acknowledges support by the Royal Society International Exchange Grant (Grant number: IES/R3/213113).

\bibliography{bibliography}
\bibliographystyle{unsrt}

\newpage
\appendix

\section{Appendix}

\subsection{Optimality conditions for Sections \ref{sec:Micro-model-std-reg}, \ref{sec:Micro-model-consis-reg}, \ref{sec:Micro-model-timecon} and \ref{sec:Micro-model-var-timecon}}\label{app:optimality-conditions}

In each case, the Hamiltonian corresponding to the opinion dynamics and respective cost functions (\ref{eq:std-reg-cost}), (\ref{eq:cost-fun-consis-control}), (\ref{eq:TC-euler}) and (\ref{eq:var-time-con-prob}) is given by
\begin{multline*}
    \quad \quad \quad H(x, p, y) = \sum_{i=2}^Np_i\frac{1}{N}\left[ \sum_{j=2}^N P(x_i, x_j)(x_j-x_i) + P(x_i,y_i)(y_i-x_i)\right] \\-\left(\frac{1}{N}\sum_{i=1}^N\frac{1}{2}(x_i-x_d)^2 + \nu\Psi(x,y)\right),
\end{multline*}
where $\Psi$ is a function given for each case in Table \ref{tab:table2} and $p$ is a vector of Lagrange multipliers. Differentiating with respect to $y_i, x_i$ and $p_i$ gives us the following set of equations for optimality,
$$ y_i = \ell(y) + \frac{k}{\nu}p_i\partial_{y_i}\{P(x_i, y_i)(y_i-x_i)\}, \quad i=2, ..., N,$$
$$ \dot{p}_i(t) = \frac{1}{N}(x_i-x_d) - \frac{1}{N}\left[ \sum_{j=2}^N R_{ij} + p_i\partial_{x_i}\{P(x_i, y_i)(y_i-x_i)\}\right], \quad p_i(t^{n+1}) = 0,\quad i=2, ..., N,$$
$$\dot{x}_i(t) = \frac{1}{N}\left(\sum_{j=2}^N P(x_i,x_j)(x_j-x_i) + P(x_i, y_i)(y_i-x_i)\right), \quad x_i(t^n) = \bar{x}_i,\quad i=2, ..., N,$$
where $R_{ij}$ is given by equation (\ref{eq:R}), $\ell(y)$ is the limit as $\nu\to\infty$ for each regularisation and $k$ is a prefactor. The explicit forms of $\ell(y)$ and $k$ are detailed in Table \ref{tab:table2}. The equation for $\dot{p}_i(t)$ is the same as that in Section \ref{sec:Micro-model-no-reg} so we have that $p_i^n$ is given by equation (\ref{eq:pin-basic}).

\begin{table}[htb!]
  \begin{center}
    \begin{tabular}{l|c|c|c}
      \textbf{Section} & \ref{sec:Micro-model-std-reg} & \ref{sec:Micro-model-consis-reg} & \ref{sec:Micro-model-timecon} \\ 
      \hline
      $\Psi(x,y)$ & $\frac{1}{N}\sum_{i=2}^N\frac{1}{2}(y_i-x_d)^2$ & $\frac{1}{2N^2}\sum_{i=2}^N\sum_{j=2}^N\frac{1}{2}(y_i-y_j)^2$ & $\frac{1}{N}\sum_{i=2}^N\frac{1}{2}\left(\frac{y_i - y_i^{n-1}}{\Delta t}\right)^2$ \\ 
      $\ell(y)$ & $x_d$ & $\frac{1}{N-1}\sum_{j=2}^N y_j$ & $y_i^{n-1}$\\
      $k$ & 1 & $\frac{N}{N-1}$ & $\Delta t^2$\\
    \end{tabular}
    \begin{tabular}{l|c}
      \textbf{Section} & \ref{sec:Micro-model-var-timecon}\\ 
      \hline
      $\Psi(x,y)$ & $\frac{1}{N}\sum_{i=2}^N\frac{1}{2}\left(\frac{y_i - y_i^{n-1}}{\Delta t} - \frac{1}{N}\sum_{j=2}^N P(y_i^{n-1}, x_j^{n-1})(x_j^{n-1} - y_i^{n-1})\right)^2$ \\ 
      $\ell(y)$ & $y_i^{n-1}+\frac{\Delta t}{N}\sum_{j=2}^NP(y_i^{n-1},x_j^{n-1})(x_j^{n-1}-y_i^{n-1})$ \\
      $k$ & $\Delta t^2$ 
    \end{tabular}
    \caption{Functions $\Psi(x,y)$, $\ell(y)$ and constant $k$ according to each of the regularisations discussed in Sections \ref{sec:Micro-model-std-reg}, \ref{sec:Micro-model-consis-reg}, \ref{sec:Micro-model-timecon} and \ref{sec:Micro-model-var-timecon}. }
    \label{tab:table2}
  \end{center}
\end{table}

\newpage
\subsection{Non-linear model predictive control algorithm with particle swarm optimisation}\label{app:nmpc-algorithm}

\begin{algorithm}
    \caption{Nonlinear model predictive control}\label{alg:sparse-nummethod}
    \begin{algorithmic}[1]
    \State \textbf{Inputs:} Maximum number of iterations, $I_{\max}$. Number of particles, $K$. Tolerance $\varepsilon\geq0$. Coefficients $c_1,c_2$. 
    \For{$n=0, ..., M-1$}
        \If{$n=0$}
            \State Draw random initial particle inputs (positions and velocities) $y^{n+h}_k(0)\in\mathcal{I}^{N-1}, v^{n+h}_k(0)\in\mathbb{R}^{N-1}$ for $k=1, ..., K$ and $h = 0, ..., H$.
        \Else
        \State Draw initial particle inputs based on the control found on the previous iteration, i.e. $y_k^{n+h}(0) = y^{n-1} + \Delta_k^{n+h}$ where $\Delta_k^{n+h} \sim N(0,\sigma^2)$ are independent identically distributed random variables for some $\sigma^2$ for each $k=1, ..., K$ and $h=0,...,H$. Draw random initial velocities $v^{n+h}_k(0)\in\mathbb{R}^{N-1}$.
        \EndIf
        \State Set $\tilde{y}^{n+h}_k(0) = y^{n+h}_k(0)$ for $k = 1, ..., K$ and $h = 0, ..., H$. 
        \State Calculate $k_{\text{best}} = \arg\min_{k = 1, ..., K}\{C_n(y):y = (\tilde{y}^n_k(0), ..., \tilde{y}^{n+H}_k(0))\}$ and set $y^n_{\text{best}}(0) = y^n_{k_{\text{best}}}(0)$. Set $m = 0$.
        \While{$m<I_{\max}$}
            \If{$C_n((y_{k_{\text{best}}}^n(m), ..., y^{n+H}_{k_{\text{best}}}(m)))<\varepsilon$}
                \State Stop and set $y^n = y^n_{k_{\text{best}}}(m)$. Proceed to Step \ref{algline:nmpc-loopback}.
            \Else
                \For{$k = 1, ..., K$}
                    \For{$h = 0, ..., H$}
                        \State Draw random variables $r^{n+h}_{k,1}(m), r^{n+h}_{k,2}(m)\sim U([0,1]^{N-1})$.
                        \State Calculate $$v^{n+h}_k(m+1) = v^{n+h}_k(m) + c_1r^{n+h}_{k,1}(m)(\tilde{y}^{n+h}_k(m) - y^{n+h}_k(m)) + c_2r^{n+h}_{k,2}(m)(y^{n+h}_{\text{best}}(m) - y^{n+h}_k(m)).$$
                        \State Set $y^{n+h}_k(m+1) = y^{n+h}_k(m) + v^{n+h}_k(m+1)$.
                    \EndFor
                    \If{$C_n((y^n_k(m+1), ..., y^{n+H}_k(m+1)))<C_n((\tilde{y}^n_k(m), ..., \tilde{y}^{n+H}_k(m)))$}
                    \State $\tilde{y}^{n+h}_k(m+1) = y^{n+h}_k(m+1)$ for each $h = 0, ..., H$. 
                    \Else
                    \State $\tilde{y}^{n+h}_k(m+1) = \tilde{y}^{n+h}_k(m)$ for each $h = 0, ..., H$.
                    \EndIf
                \EndFor
                \State Update $k_{\text{best}} = \arg\min_{k = 1, ..., K}\{C_n(y):y = (\tilde{y}^n_k(m+1), ..., \tilde{y}^{n+H}_k(m+1))\}$.
                \State Set $y^{n+h}_{\text{best}}(m+1) = y^{n+h}_{k_{\text{best}}}(m+1)$ for $h = 0, ..., H$. 
                \State Set $m = m+1$.
            \EndIf
        \EndWhile
        \State Advance dynamics (\ref{eq:xidyn}) using $y^n_{\text{best}}(m+1)$ as our simulated approximation for $y^n$.\label{algline:nmpc-loopback}
    \EndFor

    \end{algorithmic}
\end{algorithm}

\subsection{Numerical approximation of $y(x)$}\label{app:y(x)-approximation}

Assume that we want an approximation for $y(x_i)$ at the points $i = 0, ..., L$. We seek $y(x)$ to input into a finite volume scheme so we take $L$ to be not too large in order to combat computational cost (in our simulations $L=501$). We choose a $k\in\mathbb{N}$ such that the number of points we resolve over for the discretisation of $y(x)$ is $2^kL+1$ (taking $k=7$ has shown to be sufficient in our case). Consider the function 
$$ g(y;x_i) = \kappa(y - x_d) + \partial_y\left\{P(y - x_i) (y-x_i)\right\}.$$
An optimal instantaneous control must satisfy $g(y;x_i) = 0$. The initial guess we give our numerical solver is very important since $f$ has non-unique solutions. Hence, we begin with the region $|x - x_d|<\frac{r_1\kappa}{\kappa+1}$, where the control is known and iterate our initial guesses from this region. We calculate $y(x_i)$ for $i =0, ..., L$ according to Algorithm \ref{alg:y(x)-approximation}.
\begin{algorithm}
    \caption{Strategy to approximate $y(x)$}\label{alg:y(x)-approximation}
    \hspace*{\algorithmicindent} \textbf{Input:} $L\in\mathbb{N}$, desired discretisation length. $k\in\mathbb{N}$, factor by which we increase the number of discretisation points. $\text{tol}$, tolerance for $y$ to be a root of $g(y,x_i)=0$.\\
    \hspace*{\algorithmicindent} \textbf{Output:} $y(x_i)$ for $i=0, ..., L$.
    \begin{algorithmic}[1]
        \State Discretise $x$ on the interval $\mathcal{I}$ with $2^kL-1$ points. Label each point $\tilde{x}_j$ for $j = 0, ..., 2^kL$.
        \State Find $j_* =\arg\min_{j\in[0, ..., 2^kL]}|\tilde{x}_j - x_d|$ where $x_d$ is the goal opinion of the liar.
        \State Define $X_1 = [\tilde{x}_0, ..., \tilde{x}_{j_*}]$ and $X_2 = [\tilde{x}_{j_*}, ..., \tilde{x}_{2^kL}]$.
        \For{$j = j_*, ..., 0$}
        \If{$|\tilde{x}_j -x_d|<\frac{r_1\kappa}{\kappa+1}$}
        \State $$ y(\tilde{x}_j) = x_d - \frac{1}{\kappa}(\tilde{x}_j-x_d).$$
        \Else
        \State\label{algstep:newton-solve} Solve $$g(y;\tilde{x}_j) = 0,$$ where $g(y;\tilde{x}_j)$ is given by equation (\ref{eq:discrete-numeric-solve}), with initial guess $y(\tilde{x}_{j+1})$, giving some $\bar{y}$. If $|g(\bar{y}; \tilde{x}_j)|< \text{tol}$, then $y(\tilde{x}_j) = \bar{y}$. Otherwise, $y(\tilde{x}_j) = x_d$.
        \EndIf
        \EndFor
        \State Repeat with indices $j=j_*, ..., 2^kL$ and using initial guess $y(\tilde{x}_{j-1})$ for the Newton solve in Step \ref{algstep:newton-solve}.
        \State Let $Y_1' = [y(\tilde{x}_0), ..., y(\tilde{x}_{j_*})]$ and $Y_2' = [y(\tilde{x}_{j_*+1}), ..., y(\tilde{x}_{2^kL})]$. Note that we need to remove the repeated value $y(\tilde{x}_{j*})$ from $Y_2$.
        \State We are only interested in every $2^k$th entry of $Y' = [Y_1', Y_2']$, so define $Y=[y(x_0), ..., y(x_L)]$ such that
        $$ Y_p = y(x_p) = Y'_{2^kp} = y(\tilde{x}_{2^kp})$$
        for $p=0, ..., L$.
        \State Return $Y$.
    \end{algorithmic}
\end{algorithm}

We can explain the need for a very small discretisation when calculating $y(x_i)$ by the fact that $g(y;x_i)$ has multiple roots when $|x_i-x_d|>r_2$. Indeed, note that $y=x_d$ is always a root of $g(y;x_i)$. An example of $g(y,x_i)$ for $x_i$ such that $|x_i-x_d|>r_2$ is given in Figure \ref{fig:y(x)-against-f(y,x_i)}. In this example, if we set $k=3$ in Algorithm \ref{alg:y(x)-approximation} and we observe that $y(x_i) = x_d$ while $y(x_{i+1}), y(x_{i-1})\neq x_d$. We see that there are two roots of $g(y;x_i) = 0$, one at $x_d$ as expected, and another that is achieved by lying.
\begin{figure}[htb!]
    \centering
    \includegraphics[width=0.6\textwidth]{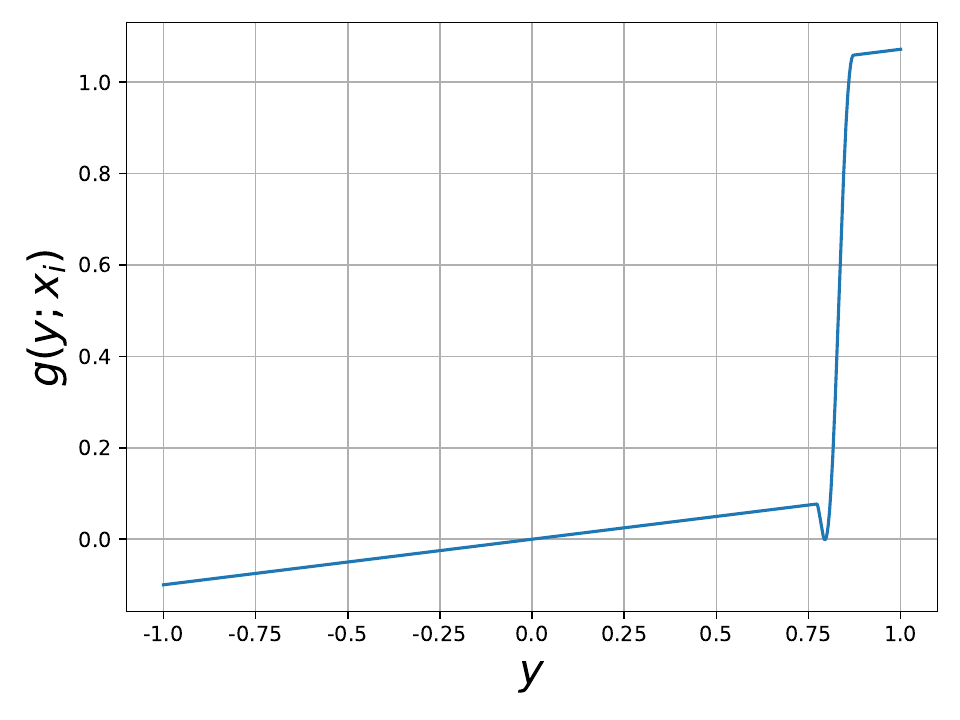}
    \caption{Example of $g(y;x_i)$ for $x_i = 0.972$ when $r_1=0.1, r_2 = 0.2$, $\kappa = 0.1$ and $x_d = 0$. We see that this function has two roots, one at $y = x_d$ and another where $r_1\leq |y-x|\leq r_2$ (the right-most root). Observe that the gradient close to the second root is very steep.}
    \label{fig:y(x)-against-f(y,x_i)}
\end{figure}

Figure \ref{fig:y(x)-against-f(y,x_i)} shows that the gradient around the root achieved by lying is very steep, meaning that we need to initialise our root finder very close to the true value of the root in order to reach it. This is why a very fine discretisation is necessary for accurate results. For particular values of $r_1, r_2$ and $\kappa$, $g(y;x_i)$ has three roots for $|x_i-x_d|>r_2$. One root of course will be $y=x_d$ but then two roots appear in the region $r_1\leq |y-x_i|\leq r_2$ with one close to the centre of the interval and one closer to $x_d$. These two roots illustrate a trade-off between a greater level of influence, $P(y-x)$, and a greater opinion shift, $y-x$.

\subsection{Finite volume scheme}\label{app:finite-volume}

We will employ a finite volume method for discretisation of the Fokker-Planck equation \cite{leveque2002finite}, similar to the one used in \cite{nugent2025opinion}. Recall the Fokker-Planck equation for general control $y(x)$ and $P(x,\cdot)$ (\ref{eq:fok-plank-boundcon}):
\begin{multline}\label{eq:general-fok-planck}
    \partial_tf_T(x,t) + \frac{1}{c_T}\partial_x\left(\left(\int_{\mathcal{I}}P(x,x')(x'-x)f_T(x',t)\, dx'\right)f_T(x,t)\right) \\+ \frac{\rho}{c_L}\partial_x\left(P(x,y)(y-x)f_T(x,t)\right) = \left(\frac{1}{c_T}+ \frac{\rho}{c_L}\right)\frac{\zeta}{2}\partial_{xx}(D(x)^2f_T(x,t)).
\end{multline}
We will take $x\in\mathcal{I} = [0,2]$ and discretise our interval by $L+1$ discrete points for $L\in\mathbb{N}$. Label these points by $x_i = i\Delta x$ for $i =0,...,L$ and where $\Delta x = 2/L$. We define the average value of our density over an interval $[x_i, x_{i+1}]$ at time $t\geq 0$ by
$$ u_i(t) = \frac{1}{\Delta x}\int_{x_{i-1}}^{x_{i}}f_T(x, t)\, dx, \quad t\geq 0, \quad i=1, ..., L.$$
We can write equation (\ref{eq:general-fok-planck}) as
\begin{equation*}
    \partial_tf_T(x,t) = \partial_xF(x,t),
\end{equation*}
where $F(x,t)$ is the flux, given by
\begin{multline}\label{eq:flux-term}
    F(x,t) = -\frac{1}{c_T}\left(\int_{\mathcal{I}}P(x,x')(x'-x)f_T(x',t)\, dx'\right)f_T(x,t) - \frac{\rho}{c_L}\left(P(x,y)(y-x)f_T(x,t)\right) \\+ \left(\frac{1}{c_T} + \frac{\rho}{c_L}\right)\frac{\zeta}{2}\partial_x(D(x)^2f_T(x,t)).
\end{multline}
Then, by our definition of $u_i(t)$, we have the differential equation
\begin{equation}\label{eq:fin-volume-idea}
    \frac{d}{dt}u_i(t) = \frac{1}{\Delta x}\int_{x_{i-1}}^{x_{i}}\partial_xF(x,t)\, dx = \frac{1}{\Delta x}(F(x_{i}, t) - F(x_{i-1},t)),
\end{equation}
for $i=1,...,L$, with initial conditions
$$ u_i(0) = \int_{x_{i-1}}^{x_{i}} f_T(x,0)\, dx.$$
By discretising the flux $F(x,t)$, we can approximate the right-hand-side of (\ref{eq:fin-volume-idea}) and then apply a time-stepping method to retrieve an approximation for $u_i(t)$. Additionally, the no-flux boundary conditions imply that
$$ F(x_0, t) = F(x_{L}, t) = 0, \quad \forall t\geq 0.$$
Define the constant matrix $\Pi\in\mathbb{R}^{L-1\times L}$ representing interactions by
$$ \Pi_{ij} = \int_{x_{j-1}}^{x_j}P(x_i,x')(x'-x_i)\, dx',$$
for $i=1, ..., L-1$ and $j = 1, ..., L$.

Next we define the matrix $F_1\in\mathbb{R}^{L-1\times L}$ to approximate the value of $f_T(x,t)$ on the midpoint of the finite volume interval $[x_{i-1}, x_i]$ for evaluating the convolution term in (\ref{eq:flux-term}),
$$ F_1 = -\frac{1}{c_T}\begin{pmatrix}
    1 &1&0&\dots&0&0\\
    0&1&1&\dots&0&0\\
    \vdots&\vdots&\vdots&\ddots&\vdots&\vdots\\
    0&0&0&\dots&1&0\\
    0&0&0&\dots&1&1
\end{pmatrix},$$
and a matrix to encode interaction with the liar $F_2\in\mathbb{R}^{L-1\times L}$, corresponding to the second term in (\ref{eq:flux-term}), given by
$$ F_2 = -\frac{\rho}{2c_L}\begin{pmatrix}
    \mathcal{Y}_1 &\mathcal{Y}_2&0&\dots&0&0\\
    0&\mathcal{Y}_2&\mathcal{Y}_3&\dots&0&0\\
    \vdots&\vdots&\vdots&\ddots&\vdots&\vdots\\
    0&0&0&\dots&\mathcal{Y}_{L-1}&0\\
    0&0&0&\dots&\mathcal{Y}_{L-1}&\mathcal{Y}_L
\end{pmatrix},$$
where 
$$ \mathcal{Y}_i = P(x_i, y(x_i))(y(x_i) - x_i)$$
for $i=1, ..., L$. Additionally, we need a matrix to approximate the spatial derivative appearing in the flux and define $F_3\in\mathbb{R}^{L-1\times L}$ as the standard forward difference matrix \cite{leveque2002finite}. Finally, we define the diffusion vector $d\in\mathbb{R}^L$ where
$$ d_i = D(x_i)^2, \quad i=1, ..., L.$$
To combine the fluxes according to (\ref{eq:fin-volume-idea}) and account for boundary conditions, we define the matrix $B\in \mathbb{R}^{L\times L-1}$ such that
$$ B = \begin{pmatrix}
    1 &0&0&\dots&0&0\\
    -1&1&0&\dots&0&0\\
    \vdots&\vdots&\vdots&\ddots&\vdots&\vdots\\
    0&0&0&\dots&-1&1\\
    0&0&0&\dots&0&-1
\end{pmatrix}.$$
Combining these discretisations, we have that an update step of $f_T(x,t)$ in time is given by
$$ u^{n+1} = u^n + \frac{\Delta t}{\Delta x}B\left((F_1u^{n})\odot(\Pi u^n) + F_2u^n + \left(\frac{1}{c_T}+\frac{\rho}{c_L}\right)\frac{\zeta}{2}F_3(d\odot u^n)\right)$$
where $\odot$ denotes a component-wise product and $u^n$ approximates $[u_1(t_n), ..., u_L(t_n)]^T$.

\end{document}

%% file: preamble.tex
\usepackage[utf8]{inputenc}
\usepackage[dvipsnames]{xcolor}
\usepackage{amsfonts}
\usepackage{amsthm}
\usepackage{amsmath}
\usepackage{parskip}
\usepackage{mathrsfs}
\usepackage{amssymb}
\usepackage{xtab}
\usepackage{setspace}
\usepackage{datetime}
\usepackage{svg}
\usepackage{csquotes}
\usepackage{fancyhdr}
\usepackage{float}
\usepackage{fontawesome}
\usepackage{comment}
\usepackage{enumitem}
\usepackage{mathalfa}
\usepackage{array}
\usepackage{titlesec}
\usepackage{mdframed}
\usepackage{listings}
\usepackage{booktabs}
\usepackage{dsfont}
\usepackage{textcomp}
\usepackage{graphicx}
\usepackage{wrapfig}
\usepackage{cases}
\usepackage{algorithm}
\usepackage{algpseudocode}

\usepackage[font=small,labelfont=bf]{caption}

\usepackage[utf8]{inputenc}
\usepackage[english]{babel}

\usepackage{placeins}
\usepackage{subcaption}
\usepackage{hyperref}
\usepackage{cleveref}

\usepackage{cite}

\usepackage{multirow}

\graphicspath{{Figures/}}

\allowdisplaybreaks

\usepackage{thmtools}
\usepackage{thm-restate}
\declaretheorem[name=Proposition,numberwithin=section]{proposition}

\declaretheorem[name=Remark,numberwithin=section]{remark}

\newcounter{subassumption}

\usepackage{circuitikz}
\usetikzlibrary{decorations.pathreplacing,calligraphy}

%% file: bibliography.bib
@article{albi2014boltzmann,
  title={Boltzmann-type control of opinion consensus through leaders},
  author={Albi, Giacomo and Pareschi, Lorenzo and Zanella, Mattia},
  journal={Philosophical Transactions of the Royal Society A: Mathematical, Physical and Engineering Sciences},
  volume={372},
  number={2028},
  pages={20140138},
  year={2014},
  publisher={The Royal Society Publishing}
}

@article{bailo2018optimal,
  title={Optimal consensus control of the \uppercase{C}ucker-\uppercase{S}male model},
  author={Bailo, Rafael and Bongini, Mattia and Carrillo, Jos{\'e} A and Kalise, Dante},
  journal={IFAC-PapersOnLine},
  volume={51},
  number={13},
  pages={1--6},
  year={2018},
  publisher={Elsevier}
}

@article{albi2014kinetic,
  title={Kinetic description of optimal control problems and applications to opinion consensus},
  author={Albi, Giacomo and Herty, Michael and Pareschi, Lorenzo},
  journal={Communications in Mathematical Sciences},
  year={2014}
}

@book{evans2005introduction,
  title={An introduction to mathematical optimal control theory},
  author={Evans, Lawrence Craig},
  year={2005},
  publisher={University of California}
}

@article{toscani2006kinetic,
  title={Kinetic models of opinion formation},
  author={Toscani, Giuseppe},
  year={2006},
  journal={Communications in Mathematical Sciences}
}

@book{leveque2002finite,
  title={Finite volume methods for hyperbolic problems},
  author={LeVeque, Randall J},
  volume={31},
  year={2002},
  publisher={Cambridge University Press}
}

@article{burger2020instantaneous,
  title={Instantaneous control of interacting particle systems in the mean-field limit},
  author={Burger, Martin and Pinnau, Ren{\'e} and Totzeck, Claudia and Tse, Oliver and Roth, Andreas},
  journal={Journal of Computational Physics},
  volume={405},
  pages={109181},
  year={2020},
  publisher={Elsevier}
}

@book{pareschi2013interacting,
  title={Interacting multiagent systems: kinetic equations and Monte Carlo methods},
  author={Pareschi, Lorenzo and Toscani, Giuseppe},
  year={2013},
  publisher={OUP Oxford}
}

@article{krause2000discrete,
  title={A discrete nonlinear and non-autonomous model of consensus formation},
  author={Krause, Ulrich and others},
  journal={Communications in difference equations},
  volume={2000},
  pages={227--236},
  year={2000}
}

@article{motsch2014heterophilious,
  title={Heterophilious dynamics enhances consensus},
  author={Motsch, Sebastien and Tadmor, Eitan},
  journal={SIAM review},
  volume={56},
  number={4},
  pages={577--621},
  year={2014},
  publisher={SIAM}
}

@article{mayne2000constrained,
  title={Constrained model predictive control: Stability and optimality},
  author={Mayne, David Q and Rawlings, James B and Rao, Christopher V and Scokaert, Pierre OM},
  journal={Automatica},
  year={2000},
  volume={36},
  number={6},
  pages={789--814},
  publisher={Elsevier}
}

@article{borghi2024kinetic,
  title={Kinetic models for optimization: a unified mathematical framework for metaheuristics},
  author={Borghi, Giacomo and Herty, Michael and Pareschi, Lorenzo},
  journal={arXiv preprint arXiv:2410.10369},
  year={2024}
}

@article{dixon2002political,
  title={Political skills or lying and manipulation? \uppercase{T}he choreography of the \uppercase{N}orthern \uppercase{I}reland peace process},
  author={Dixon, Paul},
  journal={Political studies},
  volume={50},
  number={4},
  pages={725--741},
  year={2002},
  publisher={SAGE Publications Sage UK: London, England}
}

@article{van2012learning,
  title={Learning to lie: Effects of practice on the cognitive cost of lying},
  author={Van Bockstaele, Bram and Verschuere, Bruno and Moens, Thomas and Suchotzki, Kristina and Debey, Evelyne and Spruyt, Adriaan},
  journal={Frontiers in psychology},
  volume={3},
  pages={526},
  year={2012},
  publisher={Frontiers Media SA}
}

@article{bond2004maintaining,
  title={Maintaining lies: The multiple-audience problem},
  author={Bond Jr, Charles F and Thomas, B Jason and Paulson, Ren{\'e} M},
  journal={Journal of Experimental Social Psychology},
  volume={40},
  number={1},
  pages={29--40},
  year={2004},
  publisher={Elsevier}
}

@article{allgeier1979waffle,
  title={The waffle phenomenon: Negative evaluations of those who shift attitudinally 1},
  author={Allgeier, AR and Byrne, Donn and Brooks, Barbara and Revnes, Diane},
  journal={Journal of Applied Social Psychology},
  volume={9},
  number={2},
  pages={170--182},
  year={1979},
  publisher={Wiley Online Library}
}

@article{nasr2024times,
  title={The times they are a-changin': An experimental assessment of the causes and consequences of sudden policy \uppercase{U}-turns},
  author={Nasr, Mohamed and Hoes, Emma},
  journal={European Journal of Political Research},
  volume={63},
  number={4},
  pages={1655--1673},
  year={2024},
  publisher={Wiley Online Library}
}

@article{nugent2024steering,
  title={Steering opinion dynamics through control of social networks},
  author={Nugent, Andrew and Gomes, Susana N and Wolfram, Marie-Therese},
  journal={Chaos: An Interdisciplinary Journal of Nonlinear Science},
  volume={34},
  number={7},
  year={2024},
  publisher={AIP Publishing}
}

@article{albi2017mean,
  title={Mean field control hierarchy},
  author={Albi, Giacomo and Choi, Young-Pil and Fornasier, Massimo and Kalise, Dante},
  journal={Applied Mathematics \& Optimization},
  volume={76},
  pages={93--135},
  year={2017},
  publisher={Springer}
}

@article{albi2017recent,
  title={Recent advances in opinion modeling: control and social influence},
  author={Albi, Giacomo and Pareschi, Lorenzo and Toscani, Giuseppe and Zanella, Mattia},
  journal={Active Particles, Volume 1: Advances in Theory, Models, and Applications},
  pages={49--98},
  year={2017},
  publisher={Springer}
}

@article{barrio2015dynamics,
  title={Dynamics of deceptive interactions in social networks},
  author={Barrio, Rafael A and Govezensky, Tzipe and Dunbar, Robin and Iniguez, Gerardo and Kaski, Kimmo},
  journal={Journal of the Royal Society Interface},
  volume={12},
  number={112},
  pages={20150798},
  year={2015},
  publisher={The Royal Society}
}

@article{cucker2007mathematics,
  title={On the mathematics of emergence},
  author={Cucker, Felipe and Smale, Steve},
  journal={Japanese Journal of Mathematics},
  volume={2},
  number={1},
  pages={197--227},
  year={2007},
  publisher={Springer}
}

@article{degroot1974reaching,
  title={Reaching a consensus},
  author={DeGroot, Morris H},
  journal={Journal of the American Statistical association},
  volume={69},
  number={345},
  pages={118--121},
  year={1974},
  publisher={Taylor \& Francis}
}

@article{during2024breaking,
  title={Breaking consensus in kinetic opinion formation models on graphons},
  author={D{\"u}ring, Bertram and Franceschi, Jonathan and Wolfram, Marie-Therese and Zanella, Mattia},
  journal={Journal of Nonlinear Science},
  volume={34},
  number={4},
  pages={79},
  year={2024},
  publisher={Springer}
}

@article{urena2019review,
  title={A review on trust propagation and opinion dynamics in social networks and group decision making frameworks},
  author={Urena, Raquel and Kou, Gang and Dong, Yucheng and Chiclana, Francisco and Herrera-Viedma, Enrique},
  journal={Information sciences},
  volume={478},
  pages={461--475},
  year={2019},
  publisher={Elsevier}
}

@article{tillmann2024cardinality,
  title={Cardinality minimization, constraints, and regularization: a survey},
  author={Tillmann, Andreas M and Bienstock, Daniel and Lodi, Andrea and Schwartz, Alexandra},
  journal={SIAM Review},
  volume={66},
  number={3},
  pages={403--477},
  year={2024},
  publisher={SIAM}
}

@article{kanzow2021sequential,
  title={Sequential optimality conditions for cardinality-constrained optimization problems with applications},
  author={Kanzow, Christian and Raharja, Andreas B and Schwartz, Alexandra},
  journal={Computational Optimization and Applications},
  volume={80},
  number={1},
  pages={185--211},
  year={2021},
  publisher={Springer}
}

@article{during2009boltzmann,
  title={Boltzmann and Fokker--Planck equations modelling opinion formation in the presence of strong leaders},
  author={D{\"u}ring, Bertram and Markowich, Peter and Pietschmann, Jan-Frederik and Wolfram, Marie-Therese},
  journal={Proceedings of the Royal Society A: Mathematical, Physical and Engineering Sciences},
  volume={465},
  number={2112},
  pages={3687--3708},
  year={2009},
  publisher={The Royal Society Publishing}
}

@book{livi2017nonequilibrium,
  title={Nonequilibrium statistical physics: a modern perspective},
  author={Livi, Roberto and Politi, Paolo},
  year={2017},
  publisher={Cambridge University Press}
}

@article{nugent2025opinion,
  title={Opinion Dynamics with Continuous Age Structure},
  author={Nugent, Andrew and Gomes, Susana N and Wolfram, Marie-Therese},
  journal={arXiv preprint arXiv:2503.04319},
  year={2025}
}

@article{nugent2024bridging,
  title={Bridging the gap between agent based models and continuous opinion dynamics},
  author={Nugent, Andrew and Gomes, Susana N and Wolfram, Marie-Therese},
  journal={Physica A: Statistical Mechanics and its Applications},
  volume={651},
  pages={129886},
  year={2024},
  publisher={Elsevier}
}

@article{lear2020grassmannian,
  title = {Grassmannian reduction of \uppercase{C}ucker-\uppercase{S}male systems and dynamical opinion games},
  journal = {Discrete and Continuous Dynamical Systems},
  volume = {41},
  number = {12},
  pages = {5765-5787},
  year = {2021},
  issn = {1078-0947},
  doi = {10.3934/dcds.2021095},
  url = {https://www.aimsciences.org/article/id/2baf0504-76c7-42b2-ac07-5180709e39ff},
  author = {Daniel Lear and David N. Reynolds and Roman Shvydkoy},
}

@book{shvydkoy2021dynamics,
  title={Dynamics and analysis of alignment models of collective behavior},
  author={Shvydkoy, Roman and others},
  year={2021},
  publisher={Springer}
}
